\documentclass[12pt]{amsart}
\usepackage{amssymb,amsmath,amsthm,amsfonts,mathrsfs,amsrefs,amscd,url}
\usepackage[hmargin=3cm,vmargin=3.5cm]{geometry}

\BibSpec{article}{%
    +{}  {\PrintAuthors}                {author}
    +{,} { \textit}                     {title}
    +{.} { }                            {part}
    +{:} { \textit}                     {subtitle}
    +{,} { \PrintContributions}         {contribution}
    +{.} { \PrintPartials}              {partial}
    +{,} { }                            {journal}
    +{}  { \textbf}                     {volume}
    +{}  { \PrintDatePV}                {date}
    +{,} { \issuetext}                  {number}
    +{,} { \eprintpages}                {pages}
    +{,} { }                            {status}
    +{,} { \url}                        {url}    % <---- ADDED
    +{,} { \PrintDOI}                   {doi}
    +{,} { available at \eprint}        {eprint}
    +{}  { \parenthesize}               {language}
    +{}  { \PrintTranslation}           {translation}
    +{;} { \PrintReprint}               {reprint}
    +{.} { }                            {note}
    +{.} {}                             {transition}
    +{}  {\SentenceSpace \PrintReviews} {review}
}

\usepackage[utf8]{inputenc}
\usepackage{indentfirst}
\usepackage[english]{babel}
\usepackage[unicode,verbose]{hyperref}
\usepackage[all]{xy}
\usepackage{pdfpages}
\usepackage{graphicx}
\usepackage{epstopdf}
\usepackage{braket}
\usepackage{commath}
\usepackage{mathtools}
\usepackage{empheq}
\usepackage{cancel}
\usepackage{float}
\usepackage[skip=15pt]{caption}
\newtheorem{theorem}{Theorem}
\newtheorem{corollary}{Corollary}[theorem]
\newtheorem{lemma}{Lemma}
\newtheorem{prop}{Proposition}
\theoremstyle{definition}
\newtheorem{definition}{Definition}
\newtheorem{remark}{Remark}
\newtheorem{example}{Example}
\newtheorem{conj}{Conjecture}
\newtheorem{problem}{Problem}

\numberwithin{prop}{section}
\numberwithin{conj}{section}
\numberwithin{example}{section}
\numberwithin{remark}{section}
\numberwithin{definition}{section}
\numberwithin{lemma}{section}
\numberwithin{corollary}{section}
\numberwithin{theorem}{section}
\numberwithin{problem}{section}
\numberwithin{problem}{section}

\def\R{\mathbb R}
\def\Q{\mathbb Q}
\def\Z{\mathbb Z}

\def\C{\mathbb C}
\def\S{\mathbb S}
 
\def\lra{\longrightarrow}
\newcommand{\SL}{\mathsf{SL}}
\newcommand{\SO}{\mathsf{SO}}
\newcommand{\GL}{\mathsf{GL}}
\newcommand{\CP}{\mathbb{C}P}
\newcommand{\Fl}{\mathsf{Fl}}
\newcommand{\brak}[1]{\langle #1 \rangle}
\newcommand{\und}[1]{\underline{#1}}
\newcommand{\eps}{\epsilon}
\newcommand{\Ker}{\mathsf{Ker}}
\newcommand{\Fr}{\mathsf{Fr}}
\newcommand{\mmH}{\mathrm{H}}
\newcommand{\mcH}{\mathcal{H}}
\newcommand{\slt}{\mathfrak{sl}_2}

\newcommand{\Bal}{\mathsf{Bal}}
\newcommand{\Hom}{\mathrm{Hom}}
\newcommand{\unde}{\und{\eps}}
\newcommand{\undep}{\und{\eps}^{\prime}}
\newcommand{\B}{\mathsf{B}}
\newcommand{\Seq}{\mathsf{Seq}}
\newcommand{\Tan}{\mathsf{Tan}}
\newcommand{\Ku}{\mathsf{Ku}}
\newcommand{\Rep}{\mathsf{Rep}}
\newcommand{\slthree}{\mathfrak{sl}_3}  % added 
\newcommand{\Tcob}{\mathsf{Tcob}}
\newcommand{\Webst}{\mathsf{W}} % category of usual webs 
\newcommand{\Web}{\mathsf{Web}}
\newcommand{\Uqthree}{U_q(\mathfrak{sl}_3)}
\newcommand{\Inv}{\mathsf{Inv}}

\newcommand{\undlambda}{\underline{\lambda}}
\newcommand{\mfg}{\mathfrak{g}}
\newcommand{\Dots}{\mathsf{Dot}}
\newcommand{\mbZ}{\mathbb{Z}}
\newcommand{\mcC}{\mathcal{C}}
\newcommand{\mcF}{\mathcal{F}}
\let\emptyset\varnothing

\counterwithin{figure}{section}

\title{Lectures on $\SL(3)$ foams and link homology}

\author{Mikhail Khovanov}
\address{M.K.: Department of Mathematics, Johns Hopkins University, 
3400 N. Charles Street, 
Baltimore, MD 21218, USA}
\email{khovanov@jhu.edu}
\author{Dmitry Solovyev}
\address{D.S.: Yau Mathematical Sciences Center, Tsinghua University, Beijing, China}
\email{dimsol42@gmail.com}

\makeatletter
\@namedef{subjclassname@2020}{%
  \textup{2020} Mathematics Subject Classification}


\date{July 23, 2025}

\begin{document}
\begin{abstract}
These notes are based on the three lectures that one of the authors gave at Tsinghua University in the summer of 2023 as part of the workshop on Geometric Representation Theory and Applications. They contain an introduction to the evaluation of $\SL(3)$ foams and the associated topological theory of trivalent planar graphs and foam cobordisms between them. A categorification of the Kuperberg quantum $\mathfrak{sl}_3$ web and link invariant and the Robert-Wagner $\SL(N)$ foam evaluation are reviewed as well.
\end{abstract}
\maketitle
\tableofcontents

\section*{Introduction}
Foams appeared for the first time in link homology in the categorification of the Kuperberg quantum link invariant~\cite{Kup} for the $A_2$ spider, which is also the Reshetikhin-Turaev quantum invariant for the Lie algebra $\mathfrak{sl}_3$ and its fundamental representation~\cite{RT}. The same invariant is a specialization of the HOMFLYPT two-variable polynomial to the values $(q^3,q)$ of the two parameters~\cites{HOMFLY,PT}. 

The foam approach proved subtle to generalize from $\mathfrak{sl}_3$ to $\mathfrak{sl}_n$, at first, and link homology for the fundamental representation of quantum $\mathfrak{sl}_n$ and its exterior powers was originally constructed using matrix factorizations~\cites{KhRz,Yn,Wu}. 
Since then, these and related link homology theories have been derived in many ways as various categorifications of the Reshetikhin-Turaev link and tangle invariants. In the Reshetikhin-Turaev theory~\cite{RT},~\cites{Tu,RR1,RR2}, invariants of tangles are viewed as maps between tensor products of representations of $U_q(\mathfrak{g})$. Categorification of tangle invariants is achieved by lifting tensor products of representations (or their subspaces of invariants) to Grothendieck groups of triangulated categories. The invariant  of a tangle becomes an exact functor between the categories. 

These triangulated categories and functors between them associated to tangles can be realized in an abundance of different ways. One  realization uses the highest weight category $\mathcal{O}$ of representations of $\SL(k)$, over all $k$. A categorification of $V^{\otimes k}$, for the fundamental representation $V$ of quantum $\mathfrak{sl}_2$ is achieved via maximal singular and parabolic blocks of this category, and this extends to tensor products of fundamental representations of quantum $\SL(N)$~\cites{BFK,S,Str,BrSt,MS}.  Closely related categories, functors and link homology invariants emerge from categories of coherent sheaves on suitable varieties, including convolution varieties for the affine Grassmannian~\cites{CK1, CK2} and quiver varieties~\cites{A,AN}. Mirror dual constructions realize these categories inside Fukaya-Floer categories for quiver varieties~\cites{AS,M}. 

In a tour-de-force, Webster~\cites{W1,W2} categorified Reshetikhin-Turaev invariants for links and tangles for any simple Lie algebra $\mfg$ and any coloring of the components of a tangle by irreducible representations of $U_q(\mfg)$. Key role in his construction is played by a categorification of arbitrary tensor products $V_{\lambda_1}\otimes \dots \otimes V_{\lambda_n}$ of irreducible $U_q(\mfg)$-modules to finite-dimensional algebras whose Grothendieck groups are naturally integral lattices in these tensor products. Action of tangles is given by complexes of bimodules, and composition of tangles is represented by the tensor product of bimodules. When the tangle is a link, the invariant is a bigraded homology group with the Euler characteristic given by the Reshetikhin-Turaev invariant of that labelled link. Unless representations that are used for the labels are miniscule, these homology groups seem to be nontrivial in infinitely-many bidegrees, but with the Laurent polynomial Euler characteristic in $\Z[q,q^{-1}]$. 

It is expected that Webster homology theories are functorial for decorated tangle cobordisms iff all labels of components are miniscule representations. 
 A related open problem is to find a modification of the Webster homology and tangle invariants for non-miniscule representations which is fully functorial under decorated tangle cobordisms.

In a different approach, Cautis~\cite{Cau} categorified  Reshetikhin-Turaev invariants for all representations of $\mathfrak{sl}_n$ by replacing generalized Jones-Wenzl projections by infinite twists, and using braid presentations of links. 

The functoriality problem of extending homology from links to link cobordisms is also open for the Cautis homology and several other homology theories, including for reduced triply-graded HOMFLYPT homology~\cite{KhRo3}, $\mathfrak{gl}_{-n}$ homology of Queffelec-Rose-Sartori~\cite{QRS}, Bodish-Elias-Rose spin link homology~\cite{BER} and categorifications of the colored Jones polynomial~\cites{Kh6,CoKr,BeWh,Roz,Me,StSu}.  

%Braid group action on a suitable symplectic manifold $(M,w)$ gives an induced action on the corresponding Fukaya category \cite{KhS}. Later, this idea combined with Floer theory led to the definition of the symplectic Khovanov homology \cite{SS}, which is conjectured to coincide with the original Khovanov homology. \MK{Check here, probably already proved.} Symplectic manifolds in the construction mentioned above are quiver varieties and carry hyperkahler structure. Hence, they are objects of both the Fukaya category and the category of coherent sheaves with respect to the rotated almost complex structure. Automorphism algebras in both of these categories are isomorphic to ordinary cohomology, which in some cases leads to the equivalence of these categories. Hence, it is expected that extension of Khovanov homology to tangles can be realized through derived categories of coherent sheaves on the considered quiver varieties. \DS{need to mention categories of representations of certain finite-dimensional algebras} There are other realizations, establishing connections with the geometric Satake correspondence, $(n,n)$-Springer fibers and other topics.\DS{citations}

\vspace{0.07in}

Foams explicitly reappeared in the field of categorification in~\cite{MSV} and later in~\cites{RWe,We} following a categorification of quantum groups and their representations and the discovery of categorified quantum Howe duality~\cites{LQR,QR}.

In 2017, Robert and Wagner~\cite{RW1} developed a combinatorial approach to the $\mathfrak{sl}(N)$ link homology for link components colored by arbitrary fundamental representations $\Lambda^k V$, for all $k$, where $V\cong \C^N$ is the $N$-dimensional fundamental representation, via the universal construction and their special foam evaluation formula. Their approach was extended to link and tangle cobordisms by Ehrig, Tubbenhauer, and Wedrich~\cite{ETW}. 

The present lecture notes are dedicated mostly to the $\SL(3)$ foams and their applications. 
Sections~\ref{sec-one} and~\ref{section2}  cover  unoriented $\SL(3)$ foams, where a simpler version of the Robert-Wagner formula is written down and applied to build a topological field theory for foams with boundary, following~\cite{KhR}. Unoriented $\SL(3)$ foam theory conveys most of the spirit of the Robert-Wagner evaluation and the associated topological theory in a simpler algebraic setup involving only 3 rather than $N$ variables and working in characteristic 2, where subtle signs from~\cite{RW1} are absent. 

Section~\ref{sec_oriented} explains oriented $\SL(3)$ webs and foams following~\cites{Kup,Kh3} and sketches how to build the associated homology theory of links that categorifies the Kuperberg quantum $\mathfrak{sl}_3$ invariant. Foam evaluation in the oriented $\SL(3)$ case reduces to computations in the cohomology rings of $\CP^2$ and the flag variety of $\C^3$.  

In Section~\ref{websbnd} we explain how to extend this link homology to tangles, mostly following~\cite{MPT} and a technically simpler extension of Khovanov homology to tangles and tangle cobordisms in~\cites{Kh3,Kh5}, see also~\cite{BN} for the operadic approach.  

In Section~\ref{sec_GLN} we review the Murakami-Ohtsuki-Yamada~\cite{MOY} diagrammatical calculus of intertwiners between tensor products of miniscule representations of quantum $\GL(N)$ and $\SL(N)$ and sketch the Robert-Wagner foam evaluation~\cite{RW1} in the general case of $\GL(N)$ foams, briefly listing further developments in the subject. 

Section~\ref{KMtheory} contains a short discussion of the Kronheimer-Mrowka homology theory for 3-orbifolds and its relation to the Four-Color Theorem and to  unoriented $\SL(3)$ foams.  

\vspace{0.07in}

The notion of \emph{categorification} was proposed in the visionary work of L.~Crane and I.~B.~Frenkel in 1994, approximately 30 years ago~\cite{CrFr}. Categorification of link invariants and representation structures has been an active, quickly developing and fascinating area of research for at least the last 25 years. Despite many developments, the original Crane-Frenkel problem~\cite{CrFr} of categorification of quantum $\mathfrak{sl}_2$ at a root of unity and constructing an associated 4D TQFT, a categorification of the Witten-Reshetikhin-Turaev 3-manifold invariants, remains open. 

\vspace{0.07in}

{\bf Acknowledgements:} This paper is submitted to a volume of Proceedings of the Steklov Institute of Mathematics in honor of the 60th anniversary of Dmitri Orlov. One of us (M.K.) have met Dmitri Orlov back in the undergraduate years at Turin's seminar in the Steklov Institute and enjoyed learning some algebraic geometry from him. Orlov's foundational work  on matrix factorizations brought forward their importance in mirror symmetry and algebraic geometry. From our research viewpoint, matrix factorizations were the first tool to lead to a categorification of $\SL(N)$ quantum link invariants, and more recently they appeared in some advanced developments in link homology~\cites{Ob,ObRoz}. Although these notes are based on other approaches to link homology, the deep connection between matrix factorizations and link homology is a major bridge between the two fields which, we hope, justifies this submission to the Proceedings. Happy Birthday, Dmitri!

\vspace{0.07in} 

M.K. would like to thank Tsinghua University and the Yau Center for Mathematical Sciences for the invitation to give talks which led to the present notes. M.K. was partially supported by NSF grants DMS-2204033, DMS-2446892 and Simons Collaboration Award 994328 while working on this paper. D.S. expresses gratitude for the support by the Xing Hua Scholarship program of Tsinghua University. The authors are grateful to Joshua Sussan for valuable feedback on an earlier version of this work.

\section{Unoriented \texorpdfstring{$\SL(3)$}{SL(3)} foams and their evaluation}\label{sec-one}
In the present section we introduce the category $\mathsf{Foam}$ of unoriented $\SL(3)$ foams. Then using universal construction and the evaluation formula from~\cite{KhR} we define a lax TQFT, which takes generic cross-sections of foams,  called planar webs, and assigns to them state spaces. State spaces of planar webs are graded $R$-modules, where $R$ is the ring of symmetric polynomials in three variables over a field $\Bbbk$ with $\textrm{char}\, \Bbbk=2$.

\subsection{Foams and Tait colorings}\hfill\vspace{0.02in}\\
A foam in the definition below can be more precisely called a \emph{closed unoriented $\SL(3)$ foam}. 

\begin{figure}%[H]
\begin{center}
	{\includegraphics[width=300pt]{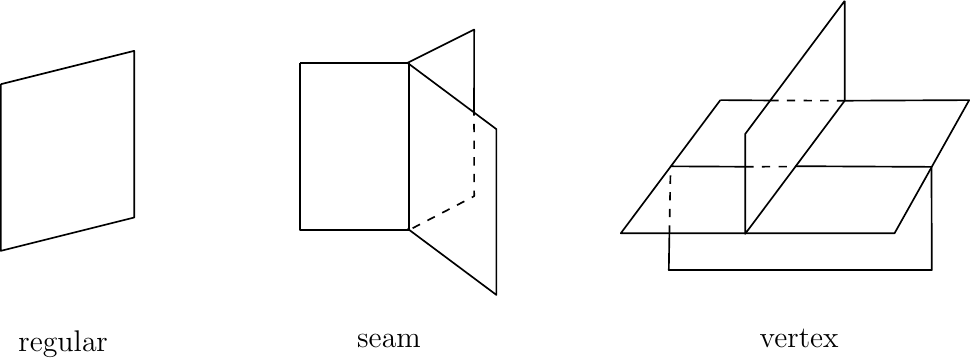}}
	\caption{Three types of points of a foam. From left to right: a regular point, a seam point, and a vertex.}
	\label{foam}
\end{center}
\end{figure}

\begin{definition}
	A foam $F\subset \mathbb{R}^3$ is a combinatorial two-dimensional compact $CW$-complex embedded in $ \mathbb{R}^3$, which can have points of the following three types, as depicted in Figure~\ref{foam}:
	\begin{itemize}
		\item a \textit{regular point}, with a neighborhood homeomorphic to a plane,
		\item a \textit{seam point}, with a neighborhood homeomorphic to a triple of half-planes joined along a line,
		\item a \textit{vertex}, with a neighborhood homeomorphic to the cone over the complete graph $K_4$ on four vertices.
	\end{itemize}
A neighbourhood of a vertex can be visualized by attaching two half-planes to a plane, one from above and one from below, with the boundaries of the two half-planes intersecting at the vertex,  as shown in Figure~\ref{foam} on the right.
Seam points form seam intervals (which end in vertices) and seam circles. At a vertex four seam intervals meet. The same interval may enter a vertex from two directions. 
 
The link of a regular point on a foam $F$ is a circle, the link of a seam point is a graph with two vertices and three edges connecting them, the link of a vertex is the complete graph $K_4$ on four vertices. 
\end{definition}

Seams and vertices are the local singularities that are allowed in foams. The union of seams and vertices of a foam $F$ is called the \textit{singular graph} of $F$ and denoted $s(F)$. It is a four-valent graph, possibly with loops, and it  may also have circles disjoint from the rest of the graph. The latter are seams that close upon themselves without encountering any vertices. We will call the set of connected components of $F\setminus s(F)$ the set of \textit{facets} $f(F)$. 

\begin{remark}
Locally a foam looks like a spine of a $3$-manifold. Consider a triangulation $T$ of a $3$-manifold $M$ then take the Poincare dual $PD(T)$. The two-skeleton $PD(T)^2$ of the Poincare dual of $T$ is a spine of $M$. Singularities of a spine are seams and vertices, as depicted in Figure~\ref{foam}.\par
\end{remark}

\begin{example}\label{ex3}
	A closed surface $S$ embedded in $\mathbb{R}^3$ is a foam.  It has only regular points, and its singular graph is empty, $s(S)=\emptyset$. Note that $S$ is necessarily orientable, has even Euler characteristic, and may have more than one connected component. 
\end{example}

\begin{example}\label{ex2}
Consider the so-called theta-foam $\Theta$ shown in Figure~\ref{theta} which consists of three disks glued together along a singular seam. 
	\begin{figure}[H]
    \begin{center}
		{\includegraphics[width=70pt]{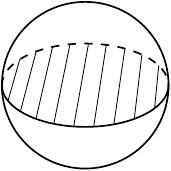}}
		\caption{The $\Theta$-foam has a singular seam which is a circle.}
		\label{theta}
    \end{center}
	\end{figure}
Its singular graph $s(\Theta)=\mathbb{S}^1$ consists of a verticeless circle, and the foam has three facets, which are all (open) two-disks.
\end{example}

\begin{remark}\label{rm_facets}
According to our definition a facet of $F$ is open in $F$. Sometimes by a facet we will mean its closure in $F$. We will either talk about \emph{open} or \emph{closed} facets, or it will be clear from the context whether an open or a closed facet is considered. 
\end{remark}

\begin{example}\label{ex1}
	Consider a foam $F$ given by  a two-torus standardly embedded in $\mathbb{R}^3$, with two disks glued in along the meridional and longitudinal cycles, respectively, as shown in Figure~\ref{ttd}.
	\begin{figure}[H]
    \begin{center}
		{\includegraphics[width=150pt]{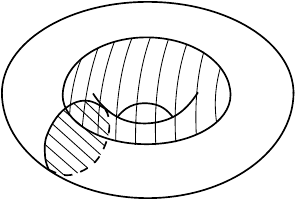}}
		\caption{Foam $F$ with three facets and a single singular vertex.}
		\label{ttd}
    \end{center}
	\end{figure}
Unlike previous examples, this foam  has a vertex point. The singular graph of this foam, $s(F)= S^1\vee S^1$, is a $4$-valent graph with a single vertex and two loops. The set of facets $f(F)$ consists of three open disks $D^2$.
\end{example}

\begin{example}\label{ex4}
	Example \ref{ex1} can be extended by replacing the meridional disk by $n$ parallel disks and the longitudional disk by $m$ parallel disks. The singular graph $s(F)$ then has $nm$ vertices connected by $2nm$ edges  (seams). It can be visualized as a square lattice on a torus with $nm$ squares. 
\end{example}

Now let us define a coloring of a foam.

\begin{definition}
	A \textit{Tait coloring} or an \emph{admissible coloring} of a foam $F$ is a map $c:f(F)\to\{1,2,3\}$  which is an assignment of a color or label in $\{1,2,3\}$ to each facet, such that facets along each seam are colored by three distinct colors. 
\end{definition}

\begin{figure}[H]
\begin{center}
	{\includegraphics[width=130pt]{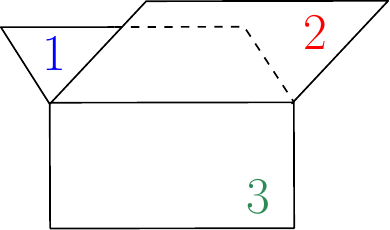}}
	\caption{Facets along each seam should be colored by three distinct colors.}
	\label{seamcolor}
    \end{center}
\end{figure}

In Example \ref{ex3}, the foam admits $3^m$ different Tait colorings, where $m$ is the number of connected components of the surface $S$, since there are three color choices for each component. $\Theta$-foam in Example \ref{ex2} foam has $6=3!$ different Tait colorings ($3$ color choices for the top facet, then $2$ choices for the second and $1$ for the third). Foam in Example $\ref{ex1}$ has no Tait colorings, because it has a seam which is approached by a single facet from two different directions. It is a good exercise to determine the number of Tait colorings of the foam in Example~\ref{ex4}, for all $n,m$. 
%Foam in  Example~\ref{ex4} with $n$ meridional and $m$ longitudional disks admits a Tait coloring iff $n,m$ are both even. In the latter case the foam has $6$ Tait colorings.   

Given a colored foam $(F,c)$, denote by $F_{ij}(c)$ the union of $i$- and $j$-colored \emph{closed} facets of $(F,c)$, where $1\leq i<j\leq 3$. (See Remark~\ref{rm_facets} on open vs. closed facets.) 
%Here we treat each facet as a closed subset of $F$ and $\mathbb{R}^3$. 
%If thinking of facets as open subsets of $F$, take their closures in the definition of $F_{ij}(c)$. 
We call $F_{ij}(c)$ \emph{the $ij$-bicolored surface} of foam $F$ and coloring $c$.  

\begin{prop}
	$F_{ij}(c)$ is an orientable surface in $\mathbb{R}^3$.
\end{prop}

\begin{figure}[H]
\begin{center}
	{\includegraphics[width=90pt]{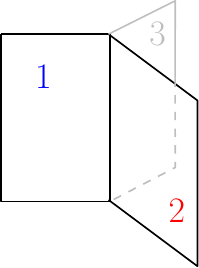}}
	\caption{Surface $F_{12}(c)$; facet colored $3$ is omitted.}
	\label{f1}
    \end{center}
\end{figure}

We need to check that $F_{ij}(c)$ is locally a surface, in a neighbourhood of each point. For regular points, that is, points on open facets colored by $i$ or $j$, that is clear. Next, consider a point on a seam, see Figure~\ref{f1}. Taking the union of two facets out of three shows that $F_{ij}(c)$ is a surface along each seam of $F$. 

\begin{figure}[H]
\begin{center}
	{\includegraphics[width=300pt]{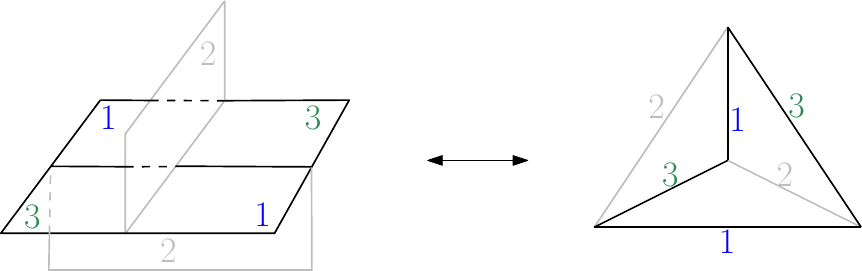}}
	\caption{An example of $F_{13}(c)$ and corresponding colored graph.}
	\label{f2}
\end{center}
\end{figure}

Next we examine the neighbourhood of a vertex in $F_{ij}(c)$. Picking distinct colors for two adjacent facets forces a unique coloring on the remaining four facets at the vertex, see Figure~\ref{f2}. In this coloring opposite facets (or opposite corners) carry the same color. Tait coloring of a neighbourhood of a vertex is unique up to permutation of colors.

Taking the union of facets colored $i$ or $j$ shows that $F_{ij}(c)$ is a surface in the neighbourhood of any vertex $v$ of $F$, and a suitable neighbourhood of $v$ in $F_{ij}(c)$ is homeomorphic to $\mathbb{R}^2$. Note that the bicolored surface $F_{ij}(c)$ contains all vertices of $F$. 

Another way to understand this property is by intersecting the foam with a small two-sphere placed at the center of the vertex. This intersection, called the link of the vertex, is a trivalent graph on the two-sphere, depicted in Figure~\ref{f2} on the right. It is the complete graph on four vertices, denoted $K_4$.
 
A Tait coloring of facets near $v$ correspond to a coloring of edges of $K_4$ so that at each vertex the colors are distinct. Such edge colorings of trivalent graphs are called Tait colorings. There is only one Tait coloring of $K_4$ up to a permutation of colors.
The union of edges with colors $i,j$ in a coloring of $K_4$ is a circle. A neighbourhood of $v$  in $F_{ij}(c)$ is given by the cone over this circle, which is another way to see that $F_{ij}(c)$ is a surface near each vertex of $F$. 

Since $F_{ij}(c)$ is a closed compact surface embedded in $\mathbb{R}^3$, it is orientable. 

Thus, for any foam $F$ and coloring $c$, the bicolored surface $F_{ij}(c)$ is a union of connected surfaces $\Sigma_g$ of various genera $g\ge 0$. Since  $F_{ij}(c)$ is closed and orientable, its  Euler characteristic $\chi_{ij}(c):=\chi(F_{ij}(c))$ is even.  Thus, for each Tait colored foam $(F,c)$ the three bicolored surfaces $F_{12}(c)$, $F_{13}(c)$ and $F_{23}(c)$ have even Euler characteristic.

\subsection{Foam evaluation from Tait colorings}
\begin{figure}[H]
\begin{center}
{\includegraphics[width=170pt]{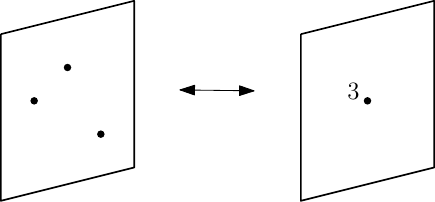}}
	\caption{Merging several dots on a facet into a single weighted dot.}
	\label{floating}
\end{center}
\end{figure}

\hfill\vspace{0.02in}\\
Let us generalize our foams by adding dots to facets. A dot can float in a facet, and several dots on the same facet may be denoted by a single dot carrying a positive integral weight, as depicted in Figure~\ref{floating}.
The dots can be thought of as observables of the theory we are about to construct.

\vspace{0.07in} 

Next we set up the algebraic side of the story. We work over a field $\Bbbk$ of characteristic two, $\textrm{char}\, \Bbbk=2$. 
The characteristic is restricted to $2$ for several reasons, including due to foams under consideration not endowed with a consistent choice of facet orientation. To a foam $F$ we will assign an invariant $\langle F \rangle\in R$ that lives in the ring $R$ of symmetric polynomials in three variables over $\Bbbk$:\begin{equation}\label{eq_ring_R}
R=\Bbbk[x_1,x_2,x_3]^{S_3}\cong \Bbbk\langle E_1, E_2, E_3\rangle,
\end{equation}
\begin{equation*}
E_1=x_1+x_2+x_3,\quad E_2=x_1 x_2+x_1 x_3+x_2 x_3,\quad E_3=x_1 x_2 x_3.
\end{equation*}
The three variables $x_1,x_2,x_3$ correspond to the three colors $1,2,3$. Let us define evaluation as follows
\begin{equation}\label{ev}
\langle F,c\rangle=\frac{x_1^{d_1(c)}x_2^{d_2(c)}x_3^{d_3(c)}}{(x_1+x_2)^{{\chi_{12}(c)}/{2}}(x_1+x_3)^{\chi_{13}(c)/{2}}(x_2+x_3)^{\chi_{23}(c)/{2}}},
\end{equation}
where
\begin{equation*}
d_i(c)=\#\{\text{dots on facets of color $i$}\},
\end{equation*}
\begin{equation*}
\chi_{ij}(c):=\chi(F_{ij}(c)).
\end{equation*}
Here $d_i(c)$ is the total number of dots on all facets colored $i$, and the Euler characteristic $\chi_{ij}(c)$ of the bicolored surface $F_{ij}(c)$ appears in an exponent in the denominator expression. 

Since we are working over a field of characteristic $2$, plus and minus signs in formula (\ref{ev}) do not matter. Replacing some plus signs with minuses will later reveal a connection to the Weyl character formula. Define the evaluation of a closed foam by 
\begin{equation}\label{eq_ev_F}
\langle F \rangle=\sum\limits_{c\in \mathrm{Tait}(F)} \langle F,c\rangle.
\end{equation}

\begin{theorem}\label{thm_in_R}
	For any foam $F$,
	\begin{equation*}
	\langle F \rangle\in R.
	\end{equation*}
\end{theorem}
This theorem states that $\langle F\rangle$ is a symmetric polynomial in $x_1,x_2,x_3$ rather than just a rational function with denominators $x_i+x_j$, $i<j$, see the definition of ring $R$ in \eqref{eq_ring_R}. That $\langle F\rangle$ is invariant under the permutation action of the symmetric group $S_3$ is immediate, since $S_3$ acts on the set of Tait colorings of $F$ by permuting the colorings, on rational functions in $x_1,x_2,x_3$ by permuting the indices, and the evaluation intertwines the two actions: $\langle F, \sigma(c)\rangle = \sigma \langle F,c\rangle$, $\sigma \in S_3$. 

\begin{figure}[H]
    \begin{center}
		{\includegraphics[width=100pt]{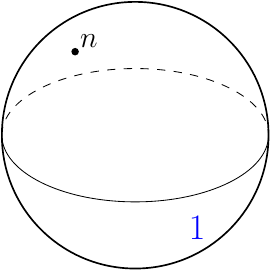}}
		\caption{The two-sphere foam. Label $1$ depics the coloring of the unique facet by that color.}
		\label{ex11}
    \end{center}
	\end{figure}
Let us look at some examples before sketching a proof of the theorem. 
\begin{example}\label{dottedsphere}
Consider the two-sphere foam $F=\mathbb{S}^2$ decorated by $n$ dots, see Figure \ref{ex11}. It has $3$ colorings. To evaluate this foam on the coloring $c_1$ by color $1$ we determine the bicolored surfaces:  
	
\begin{equation*}
	F_{12}(c_1)=\mathbb{S}^2,\quad \frac{\chi_{12}(c_1)}{2}=1,
\end{equation*}
\begin{equation*}
	F_{13}(c_1)=\mathbb{S}^2,\quad \frac{\chi_{13}(c_1)}{2}=1,
\end{equation*}
\begin{equation*}
	F_{23}(c_1)=\emptyset,\quad \frac{\chi_{23}(c_1)}{2}=0.
\end{equation*}
For instance, $F_{12}(c_1)$ is the union of facets colored $1$ or $2$. Since the unique facet of $F$ is colored $1$, this bicolored surface is the entire foam, ditto for $F_{13}(c_1)$. Surface $F_{23}(c_1)$ is empty and has Euler characteristic $0$. 
We see that the denominator term is $(x_1+x_2)(x_1+x_3)$ and 
\begin{equation*}
	\langle F,c_1\rangle = \frac{x_1^n}{(x_1+x_2)(x_1+x_3)}.
\end{equation*}
Summing over the three possible colorings of $F$, 
\begin{equation*}
	\langle F\rangle =\sum\limits_{\substack{i=1 \\\text{$j$, $k$ complementary to $i$}}}^{3}\frac{x_i^n}{(x_i+x_j)(x_i+x_k)}=h_{n-2}(x_1,x_2,x_3),	
\end{equation*}
where $h_k(x_1,x_2,x_3)$ is the complete symmetric function given by
\begin{equation*}
	h_k(x_1,x_2,x_3)\ :=\sum\limits_{a_1+a_2+a_3=k}x_1^{a_1}x_2^{a_2}x_3^{a_3}.
\end{equation*}
Polynomial $h_k$ is the character of the $k$-th symmetric power of the fundamental representation of $\SL(3)$, with coefficients modulo $2$. In particular, the two-sphere foam with no dots or one dot evaluates to $0$, with two dots --- to $1$, with three dots --- to $h_1=x_1+x_2+x_3$ and so on. Observe that for individual colorings $c$ the evaluation $\langle F,c\rangle$ is a rational function, while the sum is a symmetric polynomial in $x_1,x_2,x_3$.  
	%\begin{eqnarray*}
	%n=0\quad&\leftrightarrow&\quad \langle F \rangle = 0,\\
	%n=1\quad&\leftrightarrow&\quad \langle F \rangle = 0,\\
	%n=2\quad&\leftrightarrow&\quad \langle F \rangle = 1,\\
	%n=3\quad&\leftrightarrow&\quad \langle F \rangle = x_1+x_2+x_3,\\
	%&\ldots&
	%\end{eqnarray*}
	%and so on.
\end{example}

\begin{example}
Consider the theta-foam $\Theta$ from Example~\ref{ex2}, where the three facets are now decorated by $n_1$, $n_2$, $n_3$ dots respectively. This foam has six Tait colorings, which are in a bijection with permutations $\sigma\in S_3$, as depicted in Figure \ref{ex22}. These colorings can be denoted $c_{\sigma}$. The facets are colored by $\sigma(1),\sigma(2),\sigma(3)$, top to bottom. 
\begin{figure}[H]
\begin{center}
	{\includegraphics[width=120pt]{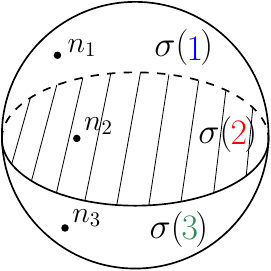}}
	\caption{Coloring $c_{\sigma}$ of a $\Theta$-foam.}
	\label{ex22}
\end{center}
\end{figure}
For any coloring $c_{\sigma}$ and any pair of colors $i,j$, $i<j$ the bicolored surface $F_{ij}(c_{\sigma})$ is a two-sphere. Consequently, the denominator term does not depend on the coloring, and  
\begin{equation*}
	F_{ij}(c_\sigma)= \mathbb{S}^2,\quad \frac{\chi_{ij}(c_\sigma)}{2}=1,
\end{equation*}
\begin{equation*}
	\langle F, c_\sigma \rangle= \frac{x^{n_1}_{\sigma(1)}x^{n_2}_{\sigma(2)}x^{n_3}_{\sigma(3)}}{(x_1+x_2)(x_2+x_3)(x_1+x_3)}.
\end{equation*}
Assuming $n_1\geq n_2\geq n_3$ and replacing some pluses by minuses (allowed in characteristic 2), we get
\begin{equation*}
	\langle F \rangle=\sum_{\sigma\in S_3} (-1)^{l(\sigma)}\frac{x^{n_1}_{\sigma(1)}x^{n_2}_{\sigma(2)}x^{n_3}_{\sigma(3)}}{(x_1-x_2)(x_2-x_3)(x_1-x_3)}=S_\lambda(x_1,x_2,x_3),
\end{equation*}
where $S_\lambda(x_1,x_2,x_3)$ is the Schur function, with coefficients reduced modulo $2$, for the partition $\lambda=(n_1-2,n_2-1,n_3)$. Schur function $S_{\lambda}$ is the character of an irreducible representation corresponding to the partition $\lambda$. In particular, if $n_1=n_2$ or $n_2=n_3$ then the partition is not defined and correspondingly $\langle F \rangle=0$.\footnote{There is a generalization of these foam evaluations to $\mathbb{Z}$ from characteristic $2$, where additionally one needs to keep track of the orientation of the facets. There is also a generalization to $N$ colors, see Section~\ref{sec_GLN}. Then, evaluation of suitable foams that generalize the $\Theta$-foam gives characters of irreducible representations of $\SL(N)$.}
\end{example}

We see that foam evaluation gives interesting functions already for the simplest foams. 

\vspace{0.07in}

Formula \eqref{ev} shows that the evaluation of a foam for a given coloring belongs to a larger ring, since 
$$
\langle F, c \rangle\in R'':=R'\left[\frac{1}{x_i+x_j}\right]_{i<j}, \ \ R':=\Bbbk[x_1,x_2,x_3], 
$$    
which is the ring $R'\supset R$ with $x_i+x_j$ for $i<j$ inverted. 
Theorem~\ref{thm_in_R} says that, unlike individual terms $\langle F, c \rangle$,  the sum \eqref{eq_ev_F} lies in the smaller ring $R$. The invariance of $\langle F, c \rangle$ under the $S_3$ action is immediate, see the remark after Theorem~\ref{thm_in_R}. Obstacles to $\langle F \rangle$ being a polynomial are the possible denominators in \eqref{ev}, which are powers of $x_i+x_j$, $i<j$. 

\vspace{0.07in} 

{\it Proof of Theorem~\ref{thm_in_R}.}
For a given coloring $c$, term $x_i+x_j$ can appear in the denominator only if some connected component of the bicolored surface 
$F_{ij}(c)$ is a two-sphere. 
A connected component of $F_{ij}(c)$ which is a two-torus contributes $1$ to the product $\langle F, c \rangle $, and higher genus components contribute powers of $x_i+x_j$ to the numerator. 

Suppose that a two-sphere $S$ is a connected component of $F_{ij}(c)$. Form another coloring $c^\prime$, where colors at all the other components are kept as is, color $k$ is kept as is, but the colors $i$ and $j$ on $S$ are swapped. Transforming coloring $c$ to the coloring $c^\prime$ is called the \textit{Kempe move} of $c$ on the two-sphere $S$.

\begin{example}
	In Figure \ref{kempe} we depict a two-sphere component of $F_{12}(c)$ for some coloring $c$. Coloring $c^\prime$ is obtained from $c$ by a Kempe move along that two-sphere. In this example, the intersection $S\cap s(F)$ of the two-sphere with the singular graph is a 4-valent graph on a sphere with 2 vertices and 4 edges. More generally, the intersection $S\cap s(F)$ is a 4-valent graph on $S$ that may also contain verticeless circles.  
	\begin{figure}[H]
    \begin{center}
		{\includegraphics[width=250pt]{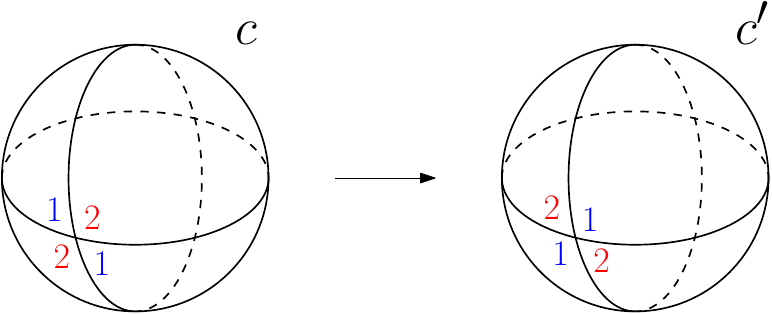}}
		\caption{Swapping colors $1,2$ on a two-sphere component of $F_{12}(c)$.}
		\label{kempe}
    \end{center}
	\end{figure}
\end{example}

The two-sphere $S$ results in the denominator $x_i+x_j$ in the evaluations $\langle F,c \rangle$ and $\langle F,c' \rangle$. The idea is to compare these two evaluations and show that $x_i+x_j$ in the denominator will cancel out in their sum $\langle F,c \rangle+\langle F,c^\prime \rangle$. Consider $i=1$, $j=2$ and assume that $S$ contains $n$, respectively $m$ dots on facets colored $1$, respectively $2$, see Figure~\ref{kempe}. 

Clearly, $\chi_{12}(c)=\chi_{12}(c^\prime)$, since $c$ and $c'$ have the same $12$-bicolored surface.

Furthermore, 
\begin{equation*}
	\chi_{13}(c)  =  \chi_{13}(c^\prime)+l, \ \ \ 
	\chi_{23}(c) =  \chi_{23}(c^\prime)-l 
	\end{equation*}
for some $l\in2\mathbb{Z}$. Writing down contributions from $S$ to the two evaluations, and contributions from the change in Euler characteristics of bicolored surfaces one obtains  
\begin{eqnarray*}
    \langle F,c \rangle &   =  &  \frac{1}{x_1+x_2}x_1^nx_2^m(x_2+x_3)^{\ell} \, y,\\
    \langle F,c' \rangle &   =  &  \frac{1}{x_1+x_2}x_2^nx_1^m(x_1+x_3)^{\ell} \, y, 
\end{eqnarray*}
for some rational function $y$, where the denominator coming from $S$ is shown as the first term in the products. 

The sum of the polynomials in the two terms is symmetric in $x_1,x_2$ and, since we work in characteristic $2$, it can be factored as follows:
$$
x_1^nx_2^m(x_2+x_3)^{\ell}  + 
x_2^nx_1^m(x_1+x_3)^{\ell} = (x_1+x_2)z,  
$$
for some polynomial $z$ in $x_1,x_2,x_3$. Factor $x_1+x_2$ cancels the same factor in the denominator in the sum $\langle F,c \rangle+\langle F,c^\prime \rangle$. 

All colorings $c\in \mathrm{Tait}(F)$ can now be partitioned into sets of order $2^k$, for various $k\ge 0$, where, for a coloring $c$ with $F_{12}(c)$ containing $k$ two-spheres, the set with $c$ contains the total of $2^k$ colorings obtained by Kempe moves on subsets of the set of $k$ spheres. 

The sum $\sum_{c'}\brak{F,c'}$ over all $c'$ in a given $(1,2)$-Kempe equivalence class has no $x_1+x_2$ terms in the denominator, by the above argument. Consequently, $\brak{F}$, being the sum of the above, over all $(1,2)$-Kempe equivalence classes, has no $x_1+x_2$ terms in the denominator. Replacing colors $1,2$ by other pair of colors $(1,3)$ and $(2,3)$ shows that $\brak{F}$ does not have $x_i+x_j$ denominator terms, for all $(i,j)$. Consequently, $\brak{F}$ is a polynomial in $x_1,x_2,x_3$. Note that the grouping of Tait colorings of $F$ into Kempe equivalence classes depends on the choice of the pair of colors $(i,j)$. 

This completes the proof of Theorem~\ref{thm_in_R}. \hfill 
$\square$

\begin{remark}
It is essential to consider foams embedded in $\mathbb{R}^3$ in our construction, since for foams not embedded in $\R^3$ bicolored surfaces $F_{ij}(c)$ are not necessarily orientable and the integrality condition $\chi_{ij}(c)\in 2\mathbb{Z}$ may fail. 
%{\it \color{red} Check Khovanov-Robert paper on whether the sum is integral if the foam is not embedded.}
\end{remark}

\begin{remark}
Reshetikhin-Turaev link invariants relate to invariants  of plane graphs $\Gamma\subset\R^2$ that represent convolutions of intertwiners and have state sum interpretations. Edges of $\Gamma$ are colored by elements of finite sets, which are in natural bijections with basis vectors of the corresponding representations. Each such coloring contributes to the quantum group evaluation of the graph via the product of contributions from vertices and local maxima and minima of the graph diagram. Summing over all colorings produces the Reshetikhin-Turaev invariant of a decorated planar graph (or of a closed planar network of intertwiners between tensor products of quantum group representations). 

Foam colorings are a simple model for going one dimension up and generalizing coloring models of planar graphs coming from quantum groups to coloring models for foams embedded in $\R^3$. It is an open problem to find more foam coloring models, beyond those in~\cites{KhR,RW1}. 
\end{remark}

\begin{remark}
    Foam evaluation formula~\eqref{eq_ev_F} was inspired by the Robert-Wagner evaluation formula~\cite{RW1} for $\GL(N)$-foams and the Kronheimer-Mrowka conjecture about the existence of consistent evaluation of unoriented $\SL(3)$ foams that naturally appear in their homology theory for 3-orbifolds~\cites{KM15,KM16,KM17}. 
\end{remark}

\subsection{Universal construction for foam evaluation}\label{subsec_uc}\hfill\vspace{0.02in}\\
A plane $T\subset\R^3$ intersects a foam $F$ {\it generically} if $T\cap F$ is a planar graph $\Gamma$ with no dots on its edges, and $N\cap F\cong \Gamma\times (-\epsilon,\epsilon)$ for a tubular neighborhood $N$ of $T$. Graph $\Gamma$ is necessarily trivalent, possibly with verticeless circles. We may also call such $\Gamma=T\cap F$ a \emph{generic cross-section} of $F$.

\begin{definition}
A foam $F$ with boundary is the intersection of a closed foam $F^\prime$ and $T \times [0, 1] \subset \R^3$ such that the planes $T\times \{0\}$ and $T\times \{1\}$ intersect $F^\prime$ generically.
\end{definition}

We view foam $F$ with boundary as a cobordism between graphs $\partial_j F = F \cap T \times \{j\}$ for $j = 0, 1$. The boundary of a foam is denoted $\partial F:=\partial_0 F\cup \partial_1 F$. In what follows we call a foam with boundary simply a foam. Any two planar trivalent graphs $\Gamma_0, \Gamma_1$ as above are the boundary graphs of some foam $F$, so the above definition is not restrictive. 
%(The above definition is a hack, but a more natural definition gives the same result, due to any planar web, as defined below, being the boundary of some foam.)  

Consider the general setup for the universal construction \cites{BHMV,Kh4} specialized to embedded foams. Suppose given a commutative ring $R$, and that to any closed foam $F\subset \R^3$ there is assigned an element $\alpha(F)\in R$ such that $\alpha(F)=\alpha (F^\prime)$ if foams $F\subset \R^3$ and $F'\subset \R^3$ are isotopic as embedded foams. Map $\alpha$ should be multiplicative under the disjoint union: 
\begin{equation*}
\alpha( F_1 \sqcup F_2)= \alpha( F_1)\alpha(F_2)
\end{equation*}
and, additionally, evaluate the empty foam to $1$: 
\begin{equation*}
\alpha( \emptyset) =1.
\end{equation*}
Evaluation (\ref{eq_ev_F}) considered above is an example of such a map.\par

\vspace{0.07in}

\begin{figure}[H]
\begin{center}
	{\includegraphics[width=300pt]{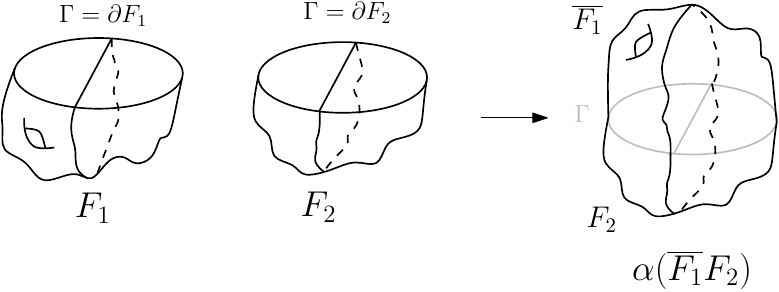}}
	\caption{$R$-valued bilinear form on $\alpha(\Gamma)$ is  obtained  by gluing two foams $F_1$ and $F_2$ along the common boundary web $\Gamma=\partial F_1=\partial F_2$ and evaluating the resulting closed foam.}
	\label{gluing}
\end{center}
\end{figure}

Define a \emph{planar web} or simply a \emph{web} to be a planar trivalent graph $\Gamma\subset \R^2$.
To a web $\Gamma$ we now assign a state space $\alpha(\Gamma)$, a module over $R$. 
First, consider all foams $F\in \R^2\times (-\infty,0]$ such that $\partial F=F\cap (\R^2\times \{0\})=\Gamma$, and let $\Fr(\Gamma)$ be the free $R$-module with a basis given by symbols $[F]$, over all such $F$. Define an $R$-bilinear form $(\,,\,)_{\Gamma}$:
\begin{eqnarray*}
	\Fr(\Gamma) \times \Fr(\Gamma) & \to & R\\
	([F_1],[F_2])_\Gamma & \mapsto & \alpha(\overline{F_1} F_2).
\end{eqnarray*}
The bilinear form on a pair of basis elements $[F_1],[F_2]$ is computed as follows. Foams $F_1,F_2$ share the common boundary $\Gamma$, and we glue the two foams along $\Gamma$ to get a closed foam in $\R^3$. That can be done, for instance, by reflecting $F_1$ about the plane $\R^2\times \{0\}$ to get the foam $\overline{F_1}$ in $\R^2\times [0,\infty)$ and then composing $\overline{F_1}$ with $F_2$ to form a closed foam $\overline{F_1}F_2$. This closed foam evaluates to an element of $R$ via $\alpha$, which is the value of the bilinear form on the given pair of vectors. 
The gluing of foams is schematically depicted in Figure \ref{gluing}.
Extending $R$-linearly, one obtains a symmetric bilinear form on $\Fr(\Gamma)$. Foams $F_1,F_2$ in this construction may have dots on them. Denote by $\Ker (\, , \,)_\Gamma\subset \Fr(\Gamma)$ the kernel of this bilinear form. \par

\vspace{0.07in} 

Next, for $\Gamma$ as above, define an $R$-module called \emph{the state space} of $\Gamma$:  
\begin{equation}\label{eq-alpha}
\alpha(\Gamma) = \Fr(\Gamma)/ \Ker (\, , \,)_\Gamma,
\end{equation}
An $R$-linear combination of foams $F_i$ with fixed boundary $\partial F_i=\Gamma$ is $0$ in the state space of $\Gamma$, 
\begin{equation}\label{eq_state}
\sum_i a_i [F_i]=0 \in \alpha(\Gamma),  \ a_i \in R,
\end{equation}
iff for any foam $F$ with boundary $\Gamma$ the following relation holds: 
\begin{equation}\label{eq_state_1}
\sum_i a_i\,\alpha(\overline{F} F_i)=0\in R.
\end{equation}
\par
That is, pairing with any foam $F$ with $\partial F = \Gamma$ and evaluating via $\alpha$ gives $0$. Relation \eqref{eq_state} is a strong condition since it requires \eqref{eq_state_1} to hold for any $F$ with $\partial F= \Gamma$. 

\vspace{0.07in} 

\begin{figure}[H]
\begin{center}
	{\includegraphics[width=200pt]{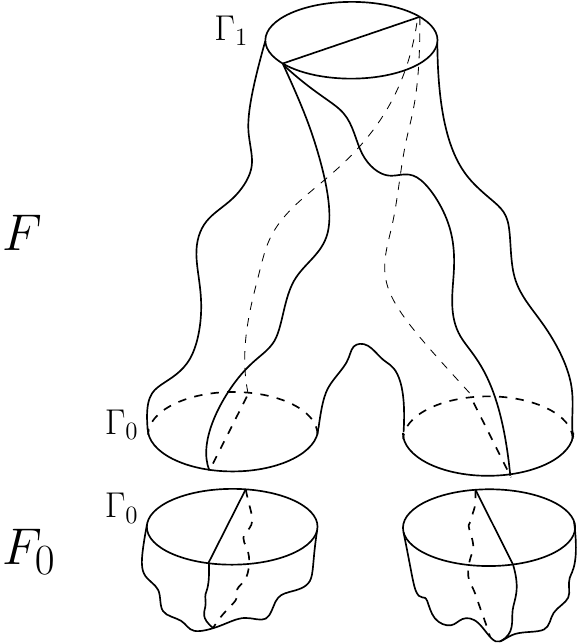}}
	\caption{Composing $F$ with various foams $F_0$ with boundary $\Gamma_0$ induces a map of state spaces $\alpha(F): \alpha(\Gamma_0)\lra \alpha(\Gamma_1)$. In this example $\Gamma_0$ has two connected components.}
	\label{funct}
\end{center}
\end{figure}

Assume now that $F\subset \R^2\times [0,1]$ is a foam with boundary, which is a cobordism between webs $\Gamma_0:=\partial_0 F$ and $\Gamma_1:=\partial_1 F$. Let $F_0\subset \R^2\times (-\infty,0]$ be a foam with boundary $\Gamma_0$. Note that $F_0$ represents a generator $[F_0]$ of $\alpha(\Gamma_0)$. 
Composing $F_0$ with $F$ induces a map of free $R$-modules
\[
\Fr(\Gamma_0)\lra \Fr(\Gamma_1), \ \ [F_0]\longmapsto [FF_0],
\]
with an example depicted in Figure~\ref{funct}.

It is immediate to see that $\Ker(\,,\,)_{\Gamma_0}$ is mapped to $\Ker(\,,\,)_{\Gamma_1}$. Consequently, 
there is an induced map on quotient spaces 
\begin{eqnarray}\label{eq_map_F}
    \alpha(F)\ :\ \alpha(\Gamma_0) & \to & \alpha(\Gamma_1)\\
    {[F_0]} & \mapsto & [FF_0].
\end{eqnarray}

These maps are consistent with composition of foams, and we obtain the following observation. 

\begin{prop}\label{prop_eval_extend}
    Multiplicative evaluation $\alpha$ as above extends to a functor 
    \begin{equation}\label{eq_func_alpha}
    \alpha \ : \ \mathsf{Foam} \lra R\text{-}\mathsf{mod} 
    \end{equation}
    from the category of foams with boundary to the category of $R$-modules. 
\end{prop}

Functor $\alpha$ takes a planar web $\Gamma$ to the $R$-module $\alpha(\Gamma)$ and a foam $F$ with boundary to the $R$-module homomorphism $\alpha(\partial_0F)\lra \alpha(\partial_1 F)$. 
\par

\vspace{0.07in}

Category $\mathsf{Foam}$ is braided monoidal, with the tensor product given by placing webs next to each other on the plane, and likewise for cobordisms between them. The functor $\alpha(\ast)$ is lax monoidal. In other words, there is a natural homomorphism 
\begin{eqnarray}\label{nathom}
	\alpha( \Gamma_1) \otimes_R \alpha( \Gamma_2)  &\to& \alpha( \Gamma_1 \sqcup \Gamma_2) \\
	{[F_1]  \otimes [F_2]} & \mapsto &  [F_1\sqcup F_2]\nonumber
\end{eqnarray}
from the tensor product of state spaces of webs $\Gamma_1,\Gamma_2$ to the state space of their disjoint union $\Gamma_1 \sqcup \Gamma_2$ given by taking a pair of foams with boundary $\Gamma_1$ and $\Gamma_2$ to their disjoint union. One of the axioms of a TQFT for foams is that homomorphism (\ref{nathom}) is an isomorphism for all $\Gamma_1,\Gamma_2$.

Let us specialize to the case when commutative ring $R$ is a field. The map (\ref{nathom}) is then injective for all multiplicative $\alpha$'s and webs $\Gamma_1,\Gamma_2$. However,  this map is usually not surjective. 
There may exist a foam $F$ with $\partial F= \Gamma_1\sqcup \Gamma_2$ such that $[F]\in \alpha(\Gamma_1\sqcup \Gamma_2)$ is not represented by an element of the tensor product $\alpha(\Gamma_1) \otimes_R \alpha( \Gamma_2)$, due to the latter being a proper subspace of the former. Consequently, the functor $\alpha$ in \eqref{eq_func_alpha} is only a \emph{lax monoidal functor} (maps \eqref{nathom} are not isomorphisms) and may be called a \emph{lax TQFT}.

\vspace{0.07in}

Clearly, there are many choices of multiplicative evaluations $\alpha$. Most of them produce uninteresting functors \eqref{eq_func_alpha}. One reason is that condition \eqref{eq_state} is very strong and, if $\alpha$ is chosen randomly, the kernel spaces $\Ker(\,,\,)_{\Gamma}$ will be trivial and $\alpha(\Gamma)$ will have infinite rank as an $R$-module.  

There are at least two natural conditions we can impose on state spaces to single out an interesting class of evaluations $\alpha$, by requiring that  

\begin{enumerate}
\item \label{cond_one} $\alpha(\Gamma)$ is a finitely-generated $R$-module for any planar web $\Gamma$, or  
\item  \label{cond_two} $\alpha(\Gamma)$ is a finitely-generated projective $R$-module for any $\Gamma$ and maps \eqref{nathom} are isomorphisms for all $\Gamma_1,\Gamma_2$. 
\end{enumerate}
Condition \eqref{cond_one} says that state spaces $\alpha(\Gamma)$ are small, in the sense of being finitely-generated $R$-modules, so that there are sufficiently many skein relations for the evaluation $\alpha$ and kernel subspaces $\Ker(\,,\,)_{\Gamma}$ are large. 
Condition \eqref{cond_two} is stronger than \eqref{cond_one} and says that $\alpha$ is a TQFT for foams.

\begin{remark} 
Dissection of any closed foam $F$ by a plane intersecting it generically gives rise to a planar trivalent graph $\Gamma$ (i.e., a planar web or simply a web). We can call this web {\it a cross-section} of $F$. To a cross-section $\Gamma$ there is assigned the $R$-module $\alpha(\Gamma)$,  {\it the state space} of $\Gamma$, see \eqref{eq-alpha}. Hence, one can think of going from an evaluation $\alpha$ of closed foams to state spaces $\alpha(\Gamma)$ for cross-sections $\Gamma$ and maps between them for foams with boundary as building the minimal model for what happens \emph{inside} the cobordisms given that the evaluation $\alpha$ is observed for closed foams (closed cobordisms).
\par
\end{remark} 

Universal construction for 3- and 4-manifolds when the ring $R$ is as large as possible (generated by the corresponding manifolds) was considered in~\cites{FKNSWW,CFW} and shown there to have markedly different features in dimensions $3$ and $4$. 

The universal construction can be applied in high generality to cobordism categories~\cite{Kh4} and, more generally, to an arbitrary monoidal category $\mathcal{C}$, see~\cite{IKO}.
One needs only a multiplicative evaluation $\alpha:\mathsf{End}_{\mathcal{C}}(\mathbf{1})\lra R$ from the commutative monoid of endomorphisms of the identity object $\mathbf{1}$ to a commutative ring $R$. 
Such an evaluation $\alpha$ gives rise to a lax monoidal functor $\alpha:\mathcal{C}\lra R$-$\mathsf{mod}$ from $\mathcal{C}$ to the category of $R$-modules. 

\vspace{0.07in} 

Let us now specialize from any multiplicative evaluation $\alpha$ of closed foams in $\R^3$ to a rather special evaluation $\brak{F}$ given by formulas \eqref{ev} and \eqref{eq_ev_F}, with the ring $R$ in \eqref{eq_ring_R}.

Denote by $\brak{\Gamma}:=\alpha(\Gamma)$ the state space of a planar web $\Gamma$ for this evaluation. We get a functor $\brak{\ast}:\mathsf{Foam}\lra R$-$\mathsf{mod}$ from the category of foams with boundary in $\R^2\times [0,1]$ to the category of $R$-modules. 

In the theory of foams dots on facets can sometimes be interpreted as first Chern classes of line bundles on flag varieties. In particular, it is natural to assign degree 2 to a dot and to generators $x_i$ of $R$, so that $\deg(x_i)=2$, and make $R$ into a graded ring. Next, extend the degree assignment to arbitrary foams as follows. 
Recall that $d_i(c)$, for a given coloring $c$ of $F$, is the number of dots on facets of color $i$. Also define $d(F)$ to be the total number of dots on facets of $F$, so that
\[
d(F) \ = \ d_1(c)+d_2(c)+d_3(c) 
\]
for any Tait coloring $c$ of $F$.

\begin{definition}\label{def_degree}
    Given a foam $F\subset \R^2\times (-\infty,0]$ such that $\partial F=\Gamma$, define the degree of the element $[F]\in \brak{\Gamma}$ by
    \begin{equation}\label{foam_degree}
        \textrm{deg}(F)=2d(F)-(\chi_{12}(c)+\chi_{23}(c)+\chi_{13}(c))
    \end{equation}
    for any Tait coloring $c$. 
\end{definition}
It is straightforward to check that $\deg(F)$ does not depend on the choice of a Tait coloring $c$ of $F$. Note also that if $F$ has no Tait colorings then $[F]=0$ in the state space $\brak{\partial F}$. 

\begin{prop}\label{prop_graded}
    $R$-module $\brak{\Gamma}$ is naturally graded, with degree of $[F]\in \brak{\Gamma}$ given by \eqref{foam_degree}. For a foam $F\subset \R^2\times [0,1]$ the induced map 
    \begin{equation}
    \brak{F}\ :\ \brak{\partial_0F}\lra \brak{\partial_1 F} 
    \end{equation}
    on state spaces has degree given by the same formula \eqref{foam_degree} but for foams in $\R^2\times[0,1]$.
\end{prop} 

\begin{corollary}\label{cor_graded}
    Functor $\brak{\ast}:\mathsf{Foam}\lra R$-$\mathsf{gmod}$
    takes values in the category of graded $R$-modules and homogeneous module morphisms. 
\end{corollary}

We record the following result.

\begin{prop}\label{prop_fg}%~\cite[Proposition~3.9]{KhR}
    The state space $\brak{\Gamma}$ is a finitely generated graded $R$-module, for any web $\Gamma$.
\end{prop}

The non-obvious part of this statement is that the graded module $\brak{\Gamma}$ is finitely generated, for which we refer to~\cite{KhR}*{Proposition~3.9}. 

\section{Skein relations and direct sum decompositions}\label{section2}
Recall that we defined a \emph{planar web} or a \emph{web} as a planar trivalent graph $\Gamma\subset \R^2$. More precisely, these can be called \emph{unoriented} $\SL(3)$ \emph{webs}. Foams with boundary appear as cobordisms between planar webs. In this section we sketch how to derive inductive direct sum decompositions for state spaces $\brak{\Gamma}$ of webs $\Gamma$ which have a region with at most four sides. We follow the strategy from~\cite{KhR}. State spaces $\brak{\Gamma}$ are graded $R$-modules, and direct sum decompositions are those of graded modules.  

\subsection{Local identities for foam evaluations and state space decomposition}\hfill\vspace{0.02in}\\
Let us derive a number of skein relations on foams, or, in other words, local relations on foam evaluations. These skein relations are then translated into direct sum decompositions of state spaces of planar webs. 

\begin{figure}[H]
	\begin{center}
		{\includegraphics[width=230pt]{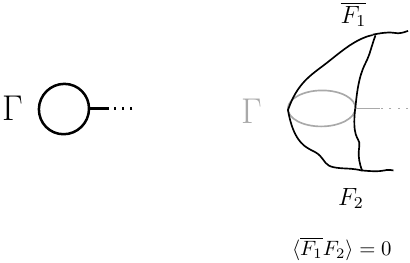}}
		\caption{The evaluation $\brak{\overline{F_1}F_2}$ is $0$ since this foam has no Tait colorings}
		\label{graphskein}
	\end{center}
\end{figure}

\begin{prop}\label{prop_trivial}
If web $\Gamma$ has a loop as in Figure \ref{graphskein}, then its state space is trivial, $\langle \Gamma \rangle=0$. 
\end{prop}

\begin{proof}
To understand the state space of $\Gamma$ we look at all foams $F_1$ such that $\partial F_1=\Gamma$, then we pair it up with $F_2$ such that $\partial F_2=\Gamma$ and compute the evaluation $\langle \overline{F_1}F_2 \rangle$. Foam $\overline{F_1}F_2$ is schematically depicted on the right side of Figure \ref{graphskein}.
Foam $\overline{F_1}F_2$ has a seam where the same facet comes from two sides. Consequently, $\overline{F_1}F_2$ has no Tait colorings and its evaluation $\langle \overline{F_1}F_2 \rangle=0$. The bilinear form on $\mathsf{Fr}(\Gamma)$ is identically zero and the state space $\langle \Gamma \rangle=0$. 
\end{proof}

%This proposition can be restated as vanishing of a state space. \MK{Discuss, perhaps don't even need this lemma, no difference with previous proposition.}
%\begin{lemma}\label{skein1}%[\cite{KhR}, Proposition 3.16]
%The following isomorphism holds. \MK{Add dots as in figure above?}
	%\begin{figure}[H]
		%\begin{center}
			%{\includegraphics[width=100pt]{skein1.pdf}}
		%\end{center}
	%\end{figure}
%\end{lemma}

\begin{example}\label{impex} 
Consider a foam $F_1$ and part of its facet shown in Figure~\ref{foam1} on the left. Form foam $F_2$, which adds to $F_1$ a bubble floating in that facet, with a dot added to the upper hemisphere, as shown in Figure~\ref{foam1} on the right. 
\begin{figure}[H]
	\begin{center}
		{\includegraphics[width=300pt]{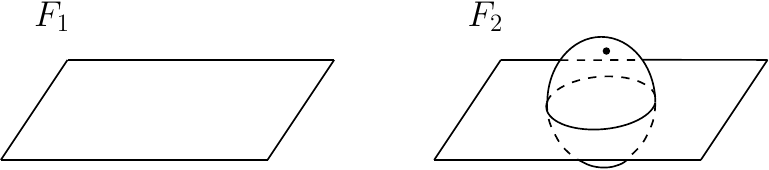}}
		\caption{Foams $F_1$ and $F_2$.}
		\label{foam1}
	\end{center}
\end{figure}
We claim that these two foams have the same evaluation, $\brak{F_1}=\brak{F_2}$. Indeed, pick a Tait coloring $c$ of the foam $F_1$, shown in Figure~\ref{foam2} on the left, with the facet colored by $i$.  There are two Tait colorings of $F_2$ that extend Tait coloring $c$ of $F_1$, shown in Figure~\ref{foam2} in the center and right and denoted $c_j$, $c_k$, the index referring to the color of the upper hemisphere. 
	\begin{figure}[H]
		\begin{center}
			{\includegraphics[width=400pt]{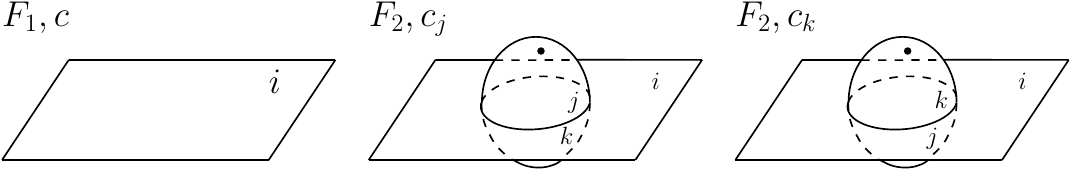}}
			\caption{Coloring $c$ of $F_1$ and its extensions to colorings $c_j,c_k$ of $F_2$.}
			\label{foam2}
		\end{center}
	\end{figure}
Let us compare $\langle F_1, c \rangle$ and $\langle F_2, c^\prime \rangle$, where $c^\prime$ stands for $c_k$ or $c_j$. Looking at the respective bicolored surfaces, we see that $F_{1,ij}(c)$ and $F_{2,ij}(c')$ are homeomorphic and likewise for $F_{1,ik}(c)$ and $F_{2,ik}(c')$. Surface $F_{2,jk}(c_j)$ is homeomorphic to the disjoint union of $F_{1,jk}(c)$ and the two-sphere:
\begin{equation*}
	F_{1,ij}(c)\cong F_{2,ij}(c_j), 
 \ \  F_{1,ik}(c)\cong F_{2,ik}(c_j), \ \ 
   F_{2,jk}(c_j)\cong F_{1,jk}(c)\sqcup \mathbb{S}^2.
\end{equation*}
We see that the bicolored surfaces for $c$ and $c_j$ contribute the same terms to $\brak{F_1,c}$ and $\brak{F_2,c_j}$, correspondingly, except for the extra term $(x_j+x_k)^{-1}$ contributed by the additional two-sphere in $F_{2,jk}(c_j)$ to the product for $\brak{F_2,c_j}$. The same computation extends to $c_k$ in place of $c_j$. Taking into account the additional dot in the upper hemisphere contributing to the numerator, we obtain
\begin{equation*}
	\langle F_2, c_j \rangle=\frac{x_j}{x_j+x_k}\langle F_1, c \rangle, \ \ 
	\langle F_2, c_k \rangle=\frac{x_k}{x_j+x_k}\langle F_1, c \rangle, 
\end{equation*}
so that
\begin{equation*}
    \brak{F_1,c} = \brak{F_2,c_j}+\brak{F_2,c_k}. 
\end{equation*}
As $c$ ranges over all colorings of $F_1$ with the given facet colored $i$, colorings $c_j,c_k$ describe the corresponding colorings of $F_2$. Summing over all colorings $c$ of $F_1$ results in the equality $\brak{F_1}=\brak{F_2}$. 
%\begin{eqnarray*}
%	\langle F_2 \rangle=\sum_{i=1}^3\big( \langle F_2, c_j \rangle+\langle F_2, c_k \rangle\big)=\sum_{i=1}^3\biggl(\frac{x_j}{x_j+x_k}\langle F_1, c \rangle+\frac{x_k}{x_j+x_k}\langle F_1, c \rangle\biggr)=\nonumber \\
%	=\sum_{i=1}^3\langle F_1, c \rangle=\sum_{c\in\mathrm{Tait}(F_1)}\langle F_1, c \rangle=\langle F_1 \rangle.
%\end{eqnarray*}
%Consequently, $\brak{F_1}=\brak{F_2}$. 
\end{example}

\begin{remark} \label{rmk_other_foams}
Consider foam $F_3$ defined as the foam $F_2$ but with no dots in either of the hemispheres, and foam $F_4$, which is the foam $F_2$ but with one dot in each hemisphere. Evaluation of both $F_3$ and $F_4$ is zero due to the ground field having characteristic $2$ and evaluations for individual colorings cancelling in pairs. %This implies that in Figure~\ref{skein2proof5} $\beta_2\alpha_1=0$ and $\beta_1\alpha_2=0$, which results in Lemma~\ref{skein2}. 
\end{remark} 

\begin{prop}%[\cite{KhR}, Proposition 2.23]
    ``Digon'' relation depicted in Figure \ref{skein2proof1} holds.
\begin{figure}[H]
	\begin{center}
		{\includegraphics[width=380pt]{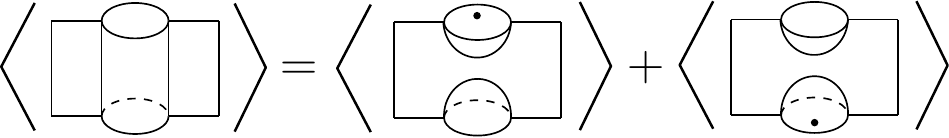}}
		\caption{An equality between foam evaluations, where the first foam on the right-hand side has a dot on the top facet away from the reader, while the second foam on the right-hand side has a dot on the bottom facet facing the reader.}
		\label{skein2proof1}
	\end{center}
\end{figure}
\end{prop}
\begin{proof}
We follow the idea of proof from Example \ref{impex}. Consider the foam $F$ on the left-hand side of Figure \ref{skein2proof1}, and pick a Tait coloring $c$ for it as depicted in Figure \ref{skein2proof2}.
\begin{figure}[H]  
	\begin{center}
		{\includegraphics[width=120pt]{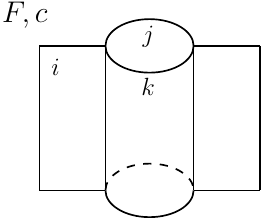}}
		\caption{Identity foam cobordism for the digon facet and its coloring $c$.}
		\label{skein2proof2}
	\end{center}
\end{figure}
There are two colorings $c^\prime$ of foams $G_1$ and $G_2$ on the right-hand side which extend the coloring $c$ of the foam $F$. We denote them  $c_j$ and $c_k$, with the index referring to the color of the back facet of the upper hemisphere, see Figure \ref{skein2proof3}.
\begin{figure}[H]
	\begin{center}
		{\includegraphics[width=250pt]{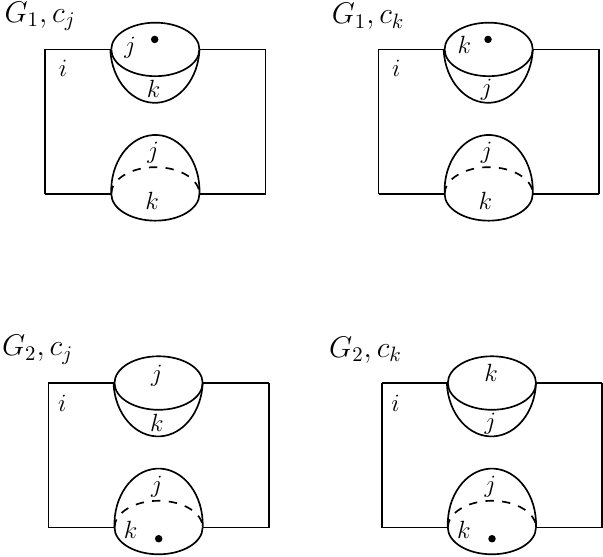}}
		\caption{Extending coloring $c$ of $F$ to colorings $c_j,c_k$ of $G_1,G_2$.}
		\label{skein2proof3}
	\end{center}
\end{figure}
Now let us compare $\langle F, c \rangle$ and $\langle G_1, c^\prime \rangle+ \langle G_2, c^\prime \rangle$. First consider intermediate foam $G_3$ with induced coloring $c^\prime$, as depicted in Figure \ref{skein2proof4}.
\begin{figure}[H]
	\begin{center}
		{\includegraphics[width=250pt]{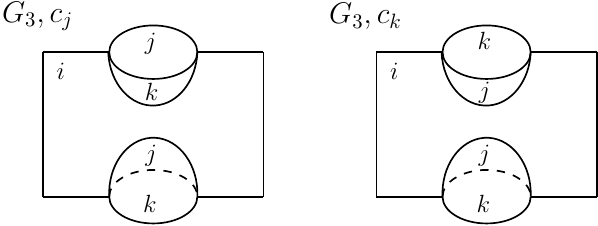}}
		\caption{Foam $G_3$ and its colorings $c_j,c_k$.}
		\label{skein2proof4}
	\end{center}
\end{figure}
Direct computation leads to the following relations
\begin{equation}\label{G131}
	\langle G_1, c_j \rangle+ \langle G_2, c_j \rangle=(x_j+x_k)\langle G_3, c_j \rangle,
\end{equation}
\begin{equation}\label{G132}
	\langle G_1, c_k \rangle+ \langle G_2, c_k \rangle=2x_k \langle G_3, c_k \rangle=0.
\end{equation}
Also note that
\begin{eqnarray*}
	\chi(F_{ij}(c))=\chi(G_{3,ij}(c_j)),\\
	\chi(F_{ik}(c))=\chi(G_{3,ik}(c_j)),
\end{eqnarray*}
but
\begin{equation*}
	\chi(F_{jk}(c))=0,\quad \chi(G_{3,jk}(c_j))=2.
\end{equation*}
 Hence,
\begin{equation*}
	\langle G_3, c_j \rangle=\frac{\langle F,c \rangle}{x_j+x_k},
\end{equation*}
and combined with (\ref{G131}\ref{G132}) we get
\begin{equation*}
    \langle F, c \rangle=\langle G_1, c_j \rangle+ \langle G_2, c_j \rangle +\langle G_1, c_k \rangle+ \langle G_2, c_k \rangle.
\end{equation*}
As $c$ ranges over all colorings $i\in\{1,2,3 \}$ of $F$, summation over $c$ corresponds to  summing over all colorings $c_j$, $c_k$ of $G_1$, $G_2$, which gives the relation depicted in Figure \ref{skein2proof1}.
%\begin{eqnarray*}
%	\langle G_1 \rangle+ \langle G_2 \rangle=\sum_{i=1}^3\big(\langle G_1, c_j \rangle+ \langle G_2, c_j \rangle +\langle G_1, c_k \rangle+ \langle G_2, c_k \rangle  \big)=\\
%	=\sum_{i=1}^3(x_j+x_k)\langle G_3, c_j \rangle=\sum_{i=1}^3\langle F, c \rangle=\langle F \rangle,
%\end{eqnarray*} 
\end{proof}
The equality implies that, irrespective of the structure of those foams outside the displayed region in Figure~\ref{skein2proof1}, evaluation of the foam on the left-hand side equals the evaluation of the foam on the right-hand side.\par

Cut the foams in that figure by a horizontal plane through the middle. 
The cross-cut of foams $G_1,G_2$ on the right hand side has a simpler structure than that of the foam $F$ on the left hand side of Figure~\ref{skein2proof1}. Together with the relations in Example~\ref{impex} and Remark~\ref{rmk_other_foams}, the above skein relation leads to a direct sum decomposition for the cross-section of $F$. 

%the following canonical decomposition of the state space of the graph from the middle section of the left-hand side of the Figure \ref{skein2proof1}:
\begin{lemma}\label{skein2}%[\cite{KhR}, Proposition 3.13]
	The following decomposition of state spaces holds
	\begin{figure}[H]
		\begin{center}
			{\includegraphics[width=320pt]{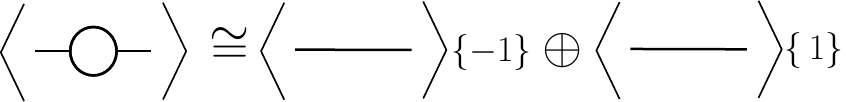}}
		\end{center}
	\end{figure}
\end{lemma}

The state space decompositions in Lemma \ref{skein2} and the following lemmas are those of graded $R$-modules. 

\begin{proof}
Consider state spaces of webs and maps between them induced by cobordisms, as depicted in Figure \ref{skein2proof5}.
\begin{figure}[H]
	\begin{center}
		{\includegraphics[width=330pt]{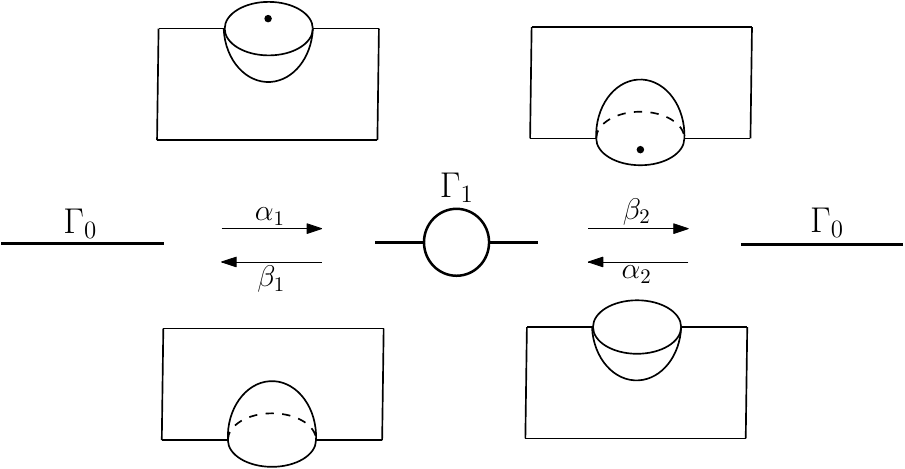}}
		\caption{Foams for the direct sum decomposition of a digon web.}
		\label{skein2proof5}
	\end{center}
\end{figure}

Cobordisms in Figure \ref{skein2proof5} are halves of the foams considered in Figure~\ref{skein2proof1}. They induce degree $\pm 1$ maps between state spaces for graphs $\Gamma_0$ and $\Gamma_1$. To make these maps degree $0$, we shift the state spaces as follows, where $\{n\}$ denotes the degree shift up by $n$: 
\begin{eqnarray*}
	\alpha_1:\langle \Gamma_0 \rangle\{1\} \to \langle \Gamma_1 \rangle,\quad \beta_1:\langle \Gamma_1 \rangle \to \langle \Gamma_0 \rangle\{1\},\\
	\alpha_2:\langle \Gamma_0 \rangle\{-1\}\to \langle \Gamma_1 \rangle,\quad \beta_2:\langle \Gamma_1 \rangle \to \langle \Gamma_0 \rangle\{-1\}.
\end{eqnarray*}
Relation in Figure~\ref{skein2proof1} then can be rewritten as
\begin{equation*}
	\textrm{id}_{\Gamma_1}=\alpha_1\beta_1+\alpha_2\beta_2,
\end{equation*}
which is part of the direct sum decomposition in Lemma~\ref{skein2}. Grading shifts in that decomposition match the degrees of maps $\alpha_i,\beta_i$ above. 

Maps in the direct sum decomposition are also given in Figure \ref{skein2proof6}.
\begin{figure}[H]
	\begin{center}
		{\includegraphics[width=350pt]{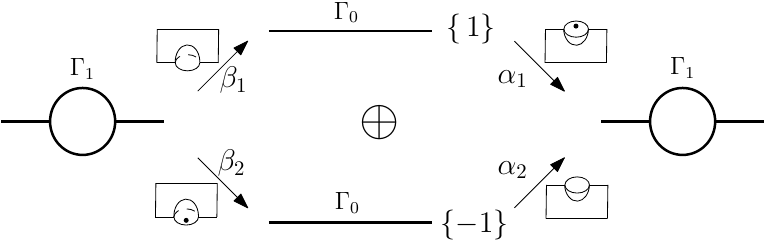}}
		\caption{}
		\label{skein2proof6}
	\end{center}
\end{figure}
Example~\ref{impex} implies that the following relations hold
\begin{equation*}
	\beta_1\alpha_1=\beta_2\alpha_2=\textrm{id}_{\Gamma_0}.
\end{equation*}
Indeed, $\beta_1\alpha_1$ is given by the foam $F_2$ in Figure~\ref{foam2}, which reduces to $F_1$ by the computation at the end of Example~\ref{impex}, implying $\beta_1\alpha_1=1$. Likewise, $\beta_2\alpha_2=1$. Moreover, Remark \ref{rmk_other_foams} implies that
\begin{equation*}
	\beta_2\alpha_1=\beta_1\alpha_2=0.
\end{equation*}
Combining all the above relations on cobordisms implies the direct sum decomposition in Lemma~\ref{skein2}.
\end{proof}

\begin{prop}%[\cite{KhR}, Proposition 2.22]
``Neck-cutting'' relation in Figure \ref{skein3proof1} holds.
\begin{figure}[H]
	\begin{center}
		{\includegraphics[width=330pt]{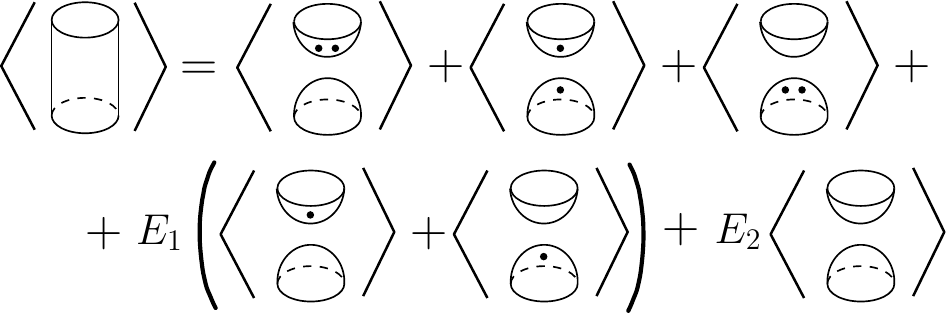}}
		\caption{A neck-cutting relation}
		\label{skein3proof1}
	\end{center}
\end{figure}
\end{prop}
\begin{proof}
Proof is similar to that of the identity in Figure~\ref{skein2proof1}. A Tait coloring of the foam $F_1$ on the left-hand side of the figure induces a coloring of each of the six terms on the right-hand side. The six terms are the same foam $F_2$ but with various number of dots on the two disks  shown (a ``cup'' and a ``cap''). In the coloring of $F_2$ induced from that of $F_1$ the cup and the cap carry the same color $i$. In the foam $F_1$ bicolored surfaces $F_{1,ij}$, for $j\not=i$, will have Euler characteristics that differ by 2 from those of surfaces $F_{2,ij}$ of the foam $F_2$. This difference will cancel with the sum of contributions from the dots and $E_1,E_2$ coefficients on the right-hand side. 

For the remaining colorings of $F_2$ the cup and the cap carry different colors $i,j$. One can then check that the right-hand side evaluation of the sum of six terms for this coloring adds up to 0. 
\end{proof}
This proposition, together with the evaluation of dotted two-spheres in Example~\ref{dottedsphere}, gives rise to the following direct sum decomposition.
\begin{lemma}\label{skein3}%[\cite{KhR}, Proposition 3.12]
	The following decomposition of state spaces holds
	\begin{figure}[H]
		\begin{center}
			{\includegraphics[width=330pt]{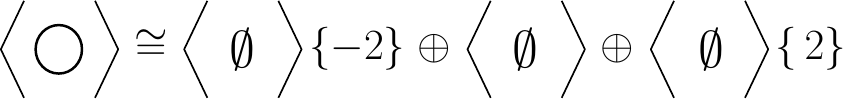}}
		\end{center}
	\end{figure}
\end{lemma}
\begin{proof}
The proof of this direct sum decomposition relies on the same idea as the one used for the proof of Lemma~\ref{skein2}. The relation in Figure \ref{skein3proof1} can be rewritten as
\begin{equation*}
    \textrm{id}_\Gamma=\alpha_1\beta_1+\alpha_2\beta_2+\alpha_3\beta_3,
\end{equation*}
as shown in Figure~\ref{skein3proof2}.
\begin{figure}[H]
	\begin{center}
		{\includegraphics[width=380pt]{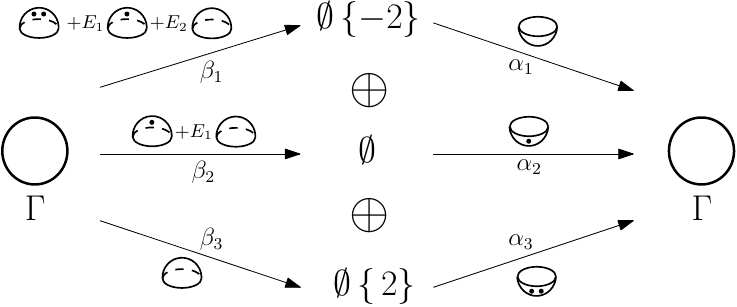}}
		\caption{Foam maps for the direct sum decomposition in Lemma~\ref{skein3}.}
		\label{skein3proof2}
	\end{center}
\end{figure}
The remaining relations for the direct sum decomposition: $\beta_i\alpha_j=\delta_{ij}$ follow from evaluations of dotted spheres in Example~\ref{dottedsphere}. 
\end{proof}
The following two relations are straightforward to verify. 
\begin{prop}%[\cite{KhR}, Proposition 2.25]
``Trivalent bubble'' relation depicted in Figure \ref{skein4proof1} holds.
\begin{figure}[H]
	\begin{center}
		{\includegraphics[width=270pt]{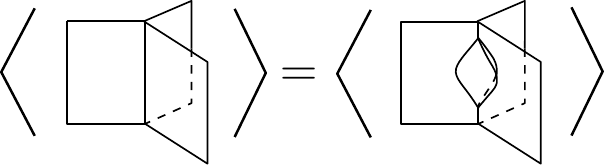}}
		\caption{Trivalent bubble relation.}
		\label{skein4proof1}
	\end{center}
\end{figure}
\end{prop}
\begin{prop}%[\cite{KhR}, Proposition 2.26]
``Vertices removal'' relation depicted in Figure \ref{skein4proof2} holds.
\begin{figure}[H]
	\begin{center}
		{\includegraphics[width=270pt]{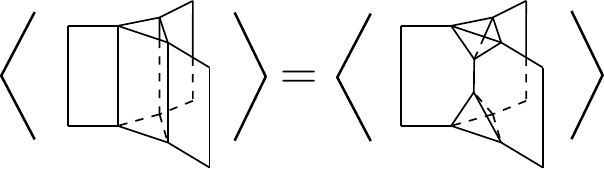}}
		\caption{Two vertices connected by an edge can be cancelled, resulting in a simpler foam with the same evaluation.}
		\label{skein4proof2}
	\end{center}
\end{figure}
\end{prop}
These, in turn, give rise to the following isomorphism
\begin{lemma}\label{skein4}%[\cite{KhR}, Proposition 3.15]
	The following webs have isomorphic state spaces 
	\begin{figure}[H]
		\begin{center}
			{\includegraphics[width=230pt]{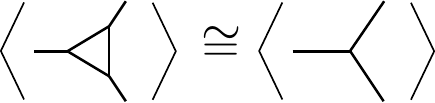}}
		\end{center}
	\end{figure}
\end{lemma}
\begin{proof}
The relations in Figures \ref{skein4proof1}, \ref{skein4proof2} can be written as
\begin{eqnarray*}
    \textrm{id}_{\Gamma_0}=\alpha\beta,\\
    \textrm{id}_{\Gamma_1}=\beta\alpha,
\end{eqnarray*}
where the foams for the maps $\alpha,\beta$ of state spaces are shown in Figure \ref{skein4proof3}.
\begin{figure}[H]
	\begin{center}
		{\includegraphics[width=330pt]{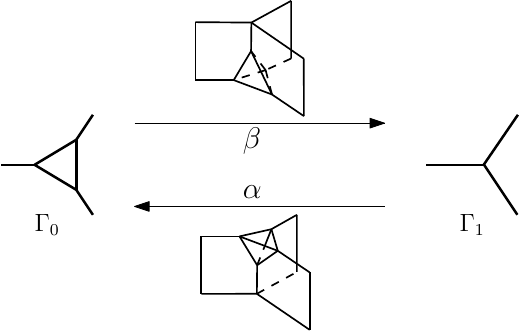}}
		\caption{Cobordisms that give mutually-inverse isomorphisms of state spaces for a web with a triangle and a contracted web.}
		\label{skein4proof3}
	\end{center}
\end{figure}
These identities follow by looking at possible colorings of foams $\alpha\beta$ and $\beta \alpha$. Colorings of $\alpha\beta$ are in a natural bijection with colorings of the identity foam $\Gamma_0\times [0,1]$, with respective bicolored surfaces homeomorphic and providing equal contributions to the evaluation functions. This implies the equality of maps of state spaces $\textrm{id}_{\Gamma_0}=\alpha\beta$. The same argument shows that $\textrm{id}_{\Gamma_1}=\beta\alpha$. 
\end{proof}

\begin{prop}%[\cite{KhR}, Proposition 2.24]
``Square decomposition'' relation depicted in Figure \ref{skein5proof1} holds.
\begin{figure}[H]
	\begin{center}
		{\includegraphics[width=400pt]{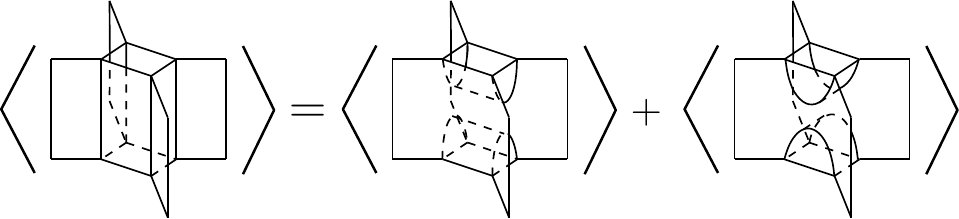}}
		\caption{A decomposition of the identity endomorphism of the square web.}
		\label{skein5proof1}
	\end{center}
\end{figure}
\end{prop}   
We refer to~\cite[Proposition 2.24]{KhR} for a proof. 
This relation, together with a relation in Remark~\ref{rmk_other_foams} and a relation to simplify a tube with two parallel disks leads to the following decomposition.
\begin{lemma}\label{skein5}%[\cite{KhR}, Proposition 3.14]
	The following direct sum decomposition of state spaces holds
	\begin{figure}[H]
		\begin{center}
			{\includegraphics[width=300pt]{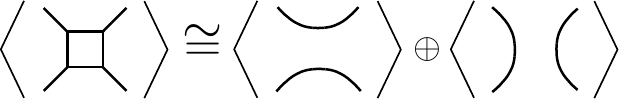}}
		\end{center}
	\end{figure}
\end{lemma}

\subsection{State spaces of planar webs, reducibility and Tait colorings}
\begin{figure}[H]
	\begin{center}
		{\includegraphics[width=50pt]{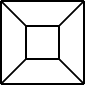}}
		\caption{This web can be reduced using skein relations \ref{skein2}, \ref{skein3} and \ref{skein5}.}
		\label{red1}
	\end{center}
\end{figure}

\hfill\vspace{0.02in}\\
We say that a web $\Gamma$ is {\it reducible} if it can be reduced to the empty web $\emptyset$ using the rules given by the relations in Proposition~\ref{prop_trivial} and Lemmas~\ref{skein2}, \ref{skein3}, \ref{skein4}, \ref{skein5}. These rules take a region of the web which has at most $n=4$ sides and reduce the web to one or several webs with fewer vertices and thus lower complexity. Proposition~\ref{prop_trivial} allows to discard a graph with a one-sided region ($n=1$). Lemma~\ref{skein2} reduces a web with a 2-sided region to a simpler web. Lemma~\ref{skein3} reduces a web with an innermost circle (which bounds a region with $0$ sides, in a sense) to the same web without the circle. Lemmas~\ref{skein4} and~\ref{skein5} reduce webs with a 3-sided, respectively 4-sided, region. Repeatedly applying these decompositions results in an inductive formula for the state space $\brak{\Gamma}$ of a reducible web $\Gamma$. In particular, the state space is a free graded module over the graded ring $R$ of symmetric polynomials in 3 variables. Figure~\ref{red1} shows an example of a reducible web.

\vspace{0.07in}

For a web $\Gamma$ we define a \emph{Tait coloring} $c$ as a map $c:e(\Gamma)\lra \{1,2,3\}$ from the set of edges of $\Gamma$ to the 3-element set of colors so that the colors of the three edges at each vertex of $\Gamma$ are distinct. If $\Gamma$ contains a circle, it can be colored by any of the 3 colors. Denote by $Tait(\Gamma)$ the set of Tait colorings of $\Gamma$ and let $t(\Gamma):=|Tait(\Gamma)|$ be the number of Tait colorings of $\Gamma$. A straightforward counting argument gives the following skein relations on the number of Tait colorings of webs:
\begin{figure}[H]
	\begin{center}
		{\includegraphics[width=350pt]{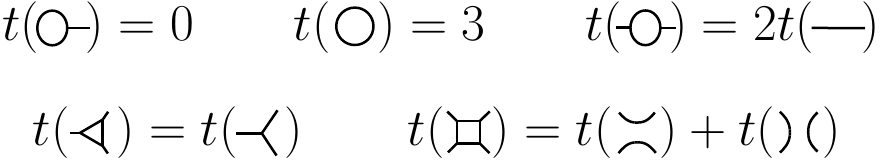}}
		\caption{Skein relations on $t(\Gamma)$.}
	\end{center}
\end{figure}
These relations mirror the above direct sum decompositions for the state spaces of webs. Induction on the complexity of a web $\Gamma$ yields the following result. 
\begin{theorem}\label{taitth}
If web $\Gamma$ is reducible then $\langle \Gamma \rangle$ is a free graded $R$-module and 
\begin{equation*}
    \mathsf{rk}_R \langle\Gamma\rangle = t(\Gamma).
\end{equation*}
\end{theorem}
For a reducible web $\Gamma$, the graded rank $\mathsf{grk}\brak{\Gamma}$ of the free graded $R$-module $\brak{\Gamma}$ takes values in $\Z_{+}[q,q^{-1}]$. It can be considered a \emph{quantization} of the number of Tait colorings of $\Gamma$, with
\begin{equation}\label{quantization}
t(\Gamma)=\mathsf{grk}\brak{\Gamma}|_{q=1}.
\end{equation}

\vspace{0.07in}

Now let us consider state spaces of some simplest planar webs: the empty web, the circle and the $\Theta$-web, see the pictures below. The number of their Tait colorings is 1, 3, and 6, respectively, while the graded ranks are given by 
\[
1, \ \ [3]=q^2+1+q^{-2}, \ \ [3][2]=(q^2+1+q^{-2})(q+q^{-1}),
\]
which specialize to $1,3,6$ upon setting $q=1$. This provides a hint that the state spaces of these webs can be interpreted via $U(3)$-equivariant cohomology of a point, a complex projective space $\mathbb{CP}^2$, and the variety $\mathsf{Fl}_3$ of full flags in $\C^3$ respectively. The field $\Bbbk$ of coefficients must have characteristic $\mathrm{char}(\Bbbk)=2$, since for unoriented foams, respectively webs, evaluations, respectively state spaces, are defined in characteristic $2$ only. Bullet $\bullet$ denotes a one-point topological space. The isomorphisms for these state spaces are 
\begin{figure}[H]
	\begin{center}
		{\includegraphics[width=220pt]{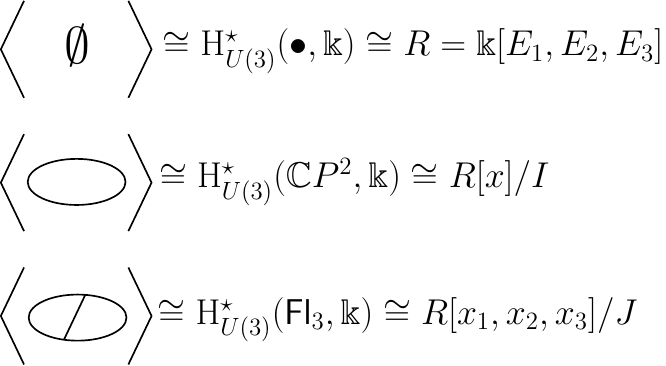}}
		\label{table}
	\end{center}
\end{figure}
where the ideals 
\begin{eqnarray*}
    I&=& (x^3+x^2E_1+x E_2+E_3)\subset R[x],\nonumber\\
    J&=&(x_1+x_2+x_3+E_1, x_1x_2+x_1x_3+x_2x_3+E_2, x_1x_2x_3+E_3)\subset R[x_1,x_2,x_3].
\end{eqnarray*}
That the state spaces of these webs are rings is due to the existence of axis of symmetry in each web. The ring structure is explained in detail below in Section~\ref{subsec_Kuperberg} for the $\Theta$-web in the orientable foam case (in particular, see Figure~\ref{3_07} there). 

\begin{figure}[H]
	\begin{center}
		{\includegraphics[width=300pt]{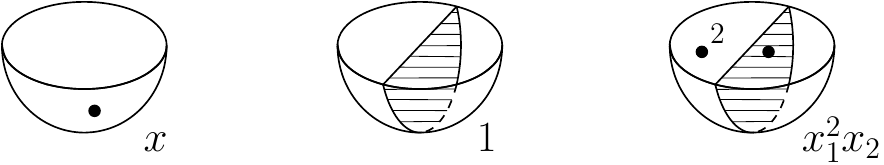}}
        \caption{Left: foam for the generator $x$ of the state space $\brak{O}$.  Center and right: the identity element and monomial $x_1^2x_2$ of $\brak{\Theta}$.}
		\label{gens}
	\end{center}
\end{figure}

Generator $x$ for the state space of the circle $\brak{O}$ is represented by a disk with a dot, see Figure~\ref{gens}. This foam naturally appears in the proof of Lemma~\ref{skein3} and generates the state space as an $R$-module. The ring structure on the state space comes from the \emph{pants} cobordism from the union of two circles to a circle.  Monomials in $x_1,x_2,x_3$ for the state space of the theta-foam $\brak{\Theta}$ are given by taking a singular cup in Figure~\ref{gens} and adding some number (possibly none) of dots to each of its facets. 

The state spaces are free graded modules over the ground ring $R$ of graded ranks $1, [3]=q^2+1+q^{-2}$ and $[3][2]=(q^2+1+q^{-2})(q+q^{-1})$, respectively. These ranks are also the Poincare polynomials of cohomology rings of a point, $\mathbb{CP}^2$ and $\mathsf{Fl}_3$, respectively, scaled by a power of $q$ so that they become Laurent polynomials invariant under the symmetry $q\leftrightarrow q^{-1}$.

Beyond these webs state spaces $\brak{\Gamma}$ quickly become complicated to describe, with some exceptions. For instance, if $\Gamma$ is obtained from a $\Theta$-web by repeatedly creating digon regions, then $\brak{\Gamma}$ can be identified with $U(3)$-equivariant cohomology of an iterated flag variety in $\C^3$. For reducible webs, the above relations allow an inductive construction of state spaces and free $R$-module bases in them given by explicit foams with boundary $\Gamma$. 

\vspace{0.07in}

\begin{figure}[H]
	\begin{center}
		{\includegraphics[width=100pt]{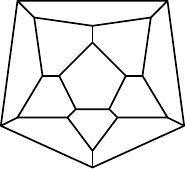}}
		\caption{Graph of the dodecahedron is not reducible.}
		\label{red2}
	\end{center}
\end{figure}

There is an abundance of examples of non-reducible webs. Any connected web in $\R^2$ such that every region, including the outside region, has at least five edges is not reducible. Figure~\ref{red2} shows the smallest example of a non-reducible graph (a non-reducible graph with the smallest number of vertices), also known as the dodecahedron graph.
Determining state spaces of non-reducible graphs is an open problem, see Boozer~\cites{B1,B2}. One can propose the following problem. 

\begin{problem}
    Is the state space  $\brak{\Gamma}$ of any planar web $\Gamma$ a free $R$-module $\langle \Gamma \rangle$ of rank $t(\Gamma)$? 
\end{problem}
This property of the state space $\brak{\Gamma}$ holds if $\Gamma$ is reducible, see Theorem~\ref{taitth}. There is rather limited evidence in favor of this conjecture for non-reducible $\Gamma$, see~\cites{B1,B2}. 

\vspace{0.07in} 

Now let us consider the state space of the disjoint union of webs. An element of the state space $\brak{\Gamma}$ can be represented by a finite $R$-linear combination of foams $F$ with boundary $\Gamma$. Given two webs $\Gamma_1,\Gamma_2$, position them in the plane away from each other. Taking the disjoint union of foams $F_1\sqcup F_2$ where $\partial F_1=\Gamma_1$ and $\partial F_2=\Gamma_2$ results in a web with a boundary $\Gamma_1\sqcup \Gamma_2$. The disjoint union is compatible with evaluations and induces a map of graded $R$-modules
\begin{equation}\label{eq_disjoint}
\brak{\Gamma_1}\otimes_R \brak{\Gamma_2} \to \brak{\Gamma_1 \sqcup \Gamma_2}.
\end{equation}
It is unknown whether this map is an isomorphism in general, since there might exist a foam $F$ with $\partial F=\Gamma_1\sqcup \Gamma_2$, which is not representable in $\brak{\Gamma_1\sqcup \Gamma_2}$ by an $R$-linear combination of unions of foams with boundaries $\Gamma_1$ and $\Gamma_2$ (also see an earlier discussion after Proposition~\ref{prop_eval_extend}). However, the state space of the disjoint union has good properties when at least one of the two webs is reducible. 

\begin{prop}
    If web $\Gamma_1$ or $\Gamma_2$ is reducible, then the map \eqref{eq_disjoint} is an isomorphism.
\end{prop}
This result implies that the evaluation of foams and the state space functor $\brak{\ast}$ is a monoidal functor (a TQFT) from the category of foams between reducible webs to the category of graded $R$-modules. 

\begin{corollary}
    Restricted to reducible webs and all foams (cobordisms) between them, the functor $\brak{\ast}$ is a TQFT.
\end{corollary} 

On the decategorified level we are assigning Tait colorings to webs. On the categorified level, we are dealing with Tait colorings of foams, one dimension up. They are needed for the evaluation of closed foams, which then allows us to define state spaces for planar webs, which are generic cross-sections of foams. Specializing to reducible webs recovers the number of Tait colorings as the rank of the module (state space) associated to the web. 

Tait colorings of webs are closely related to the Four-Color Theorem. The latter can be rephrased as the statement that a connected web has a Tait coloring iff it is bridgeless. A bridge is an edge of the web that separates a connected component into two components. See Section~\ref{KMtheory} for more on this topic. 

\begin{remark}
Unoriented $\SL(4)$ foams are studied in~\cite{KPRS}. An interesting question, raised by the authors of that paper, is whether the surjective homomorphism $S_4\lra S_3$ of symmetric groups can be lifted to a relation between unoriented $\SL(3)$ and $\SL(4)$ web state spaces, where one considers $\SL(4)$ webs with all edges of thickness 2. For the latter webs, their Tait colorings assign a pairs of colors $\{i,j\}\subset \{1,2,3,4\}$ to each edge. These pairs can be combined into three groups
\[
(\{1,2\},\{3,4\}), \ \ (\{1,3\},\{2,4\}), \ \ (\{1,4\},\{2,3\}),
\]
with $S_3$ acting on them by permutations, establishing a connection with Tait colorings of $\SL(3)$ webs. Whether this relation extends to the categorified level and to foam evaluations is an open question. 
\end{remark}
\section{Oriented \texorpdfstring{$\SL(3)$}{SL(3)} foams and link homology}
\label{sec_oriented}
Now we turn to oriented $\SL(3)$ foams. The universal construction starts with a multiplicative evaluation for foams and builds state spaces for cross-sections, which for oriented foams are oriented planar webs. Constructing state spaces in this case is less complicated than in the unoriented case. The graded rank of an oriented $\SL(3)$ web state space is the Kuperberg quantum $\slthree$ invariant of webs and links~\cite{Kup}.

\subsection{State space decomposition for oriented graphs}\label{cow}
\begin{figure}[H]
	\begin{center}
		{\includegraphics[width=150pt]{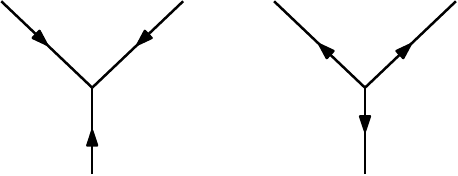}}
		\caption{In-vertex and an out-vertex of an oriented $\SL(3)$ web.}
		\label{3_01}
	\end{center}
\end{figure}

\hfill\vspace{0.02in}\\
In an \emph{oriented} $\SL(3)$ web $\Gamma$ edges and circles (if any) are oriented so that each vertex is either an \emph{in-vertex}, with all edges pointing in, or an \emph{out-vertex}, with all edges pointing out, see Figure~\ref{3_01}. Every edge goes from an out-vertex to an in-vertex. 

\begin{figure}%[H]
	\begin{center}
		{\includegraphics[width=200pt]{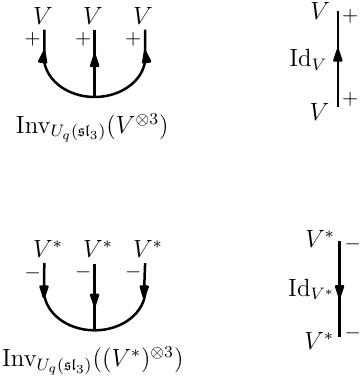}}
		\caption{Left: a trivalent out-vertex, respectively in-vertex, represents a vector in the one-dimensional space of quantum $\mathfrak{sl}_3$ invariants of $V^{\otimes 3}$, respectively in $\textrm{Inv}_{\Uqthree}((V^{\ast})^{\otimes 3})$. Right: vertical lines denote the identity endomorphism of $V$ and $V^{\ast}$, depending on line's orientation.}
		\label{3_03}
	\end{center}
\end{figure}

\begin{figure}%[H]
	\begin{center}
		{\includegraphics[width=280pt]{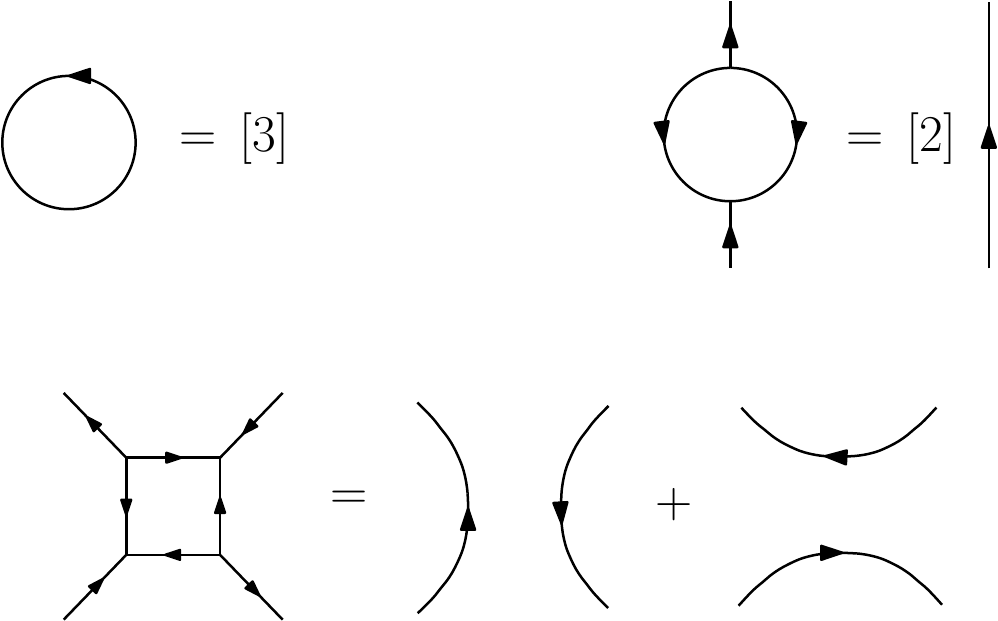}}
		\caption{Kuperberg skein relations for the $\SL(3)$ web invariant. We use a shortcut writing  $\Gamma$ in place of  $P(\Gamma)$.}
		\label{3_02}
	\end{center}
\end{figure}

Webs were first considered by Kuperberg, who used them to develop a uniform approach to quantum invariants for Lie algebras of rank two \cite{Kup}. For the root system $A_2$ Kuperberg webs come from representation theory of the quantum group $\Uqthree$, for a generic $q$. A web describes a planar convolution of generating intertwiners between tensor products of the fundamental 3-dimensional representation $V$ of $\Uqthree$ and its dual $V^{\ast}$. Edges of a web represent the identity intertwiner $\textrm{id}:V\to V$. An \emph{out-vertex} corresponds to the generator of the invariant space $\textrm{Inv}_{\Uqthree}(V^{\otimes 3})$, and an \emph{in-vertex} to the generator of $\textrm{Inv}_{\Uqthree}((V^{\ast})^{\otimes 3})$, as shown in Figure~\ref{3_03}. 
Kuperberg defines a $\Z_+[q,q^{-1}]$-valued invariant of $A_2$-webs inductively on the complexity of a web, via Figure~\ref{3_02} skein relations. We denote his invariant of $\Gamma$ by $P(\Gamma)$. In Figure \ref{3_02} and from now on we omit writing $P(\Gamma)$ in skein relations and simply draw part of the web $\Gamma$ being simplified in the relation.
Here 
\begin{equation}\label{eq_quantum_n}
[n]:=\frac{q^n-q^{-n}}{q-q^{-1}}= q^{n-1}+q^{n-3}+\ldots+q^{-(n-1)}
\end{equation}
is the $q$-\emph{deformed} integer $n$. 

These relations are enough to evaluate the invariant on any planar web and they are consistent. Define the complexity of a web $\Gamma$ as the pair $(v(\Gamma),e(\Gamma))$, where $v(\Gamma)$ is the number of vertices of $\Gamma$ and $e(\Gamma)$ the number of its edges, with a circle being counted as a single edge. These pairs of non-negative integers are ordered lexicographically. Relations in Figure~\ref{3_02} reduce the complexity of a web, with webs on the right-hand side of a relation having lower complexity than the web on the left-hand side.  

A region of a web has an even number of sides, since orientations of its edges alternate as one walks along its boundary. The standard Euler characteristic argument, applied to planar trivalent graphs, proves that such a web always has a region with at most 5 sides. Due to the parity constraint, such a region has $0$, $2$ or $4$ sides, with a $0$-sided region being an innermost circle as shown in the top left of Figure~\ref{3_02}. Relations in the figure show how to reduce a region with $0$, $2$, $4$ sides, with the resulting webs of lower complexity than the original. Kuperberg provides an easy consistency argument~\cite{Kup} to show that his invariant is well-defined and does not depend on the induction rules. Alternatively, one can check the invariance via the state sum formula for $P(\Gamma)$ coming from the intertwiners for the quantum group $\Uqthree$~\cites{Kup,KhK}. 

%{\it \color{blue} Remove?
%Any oriented $\SL(3)$ web is reducible, since it necessarily has a region of valency $4$, $2$ or $0$. A neighbourhood of such a region can be reduced using skein relations analogous to the ones derived previously, as depicted in Figure \ref{3_00}
%}
%\MK{I think we can remove this figure:}
%\begin{figure}[H]
%	\begin{center}
%		{\includegraphics[width=320pt]{3_00.pdf}}
%		\caption{Kuperberg reduction rules for oriented $\SL(3)$ webs.}
%		\label{3_00}
%	\end{center}
%\end{figure}

Hence, to an oriented $\SL(3)$ web $\Gamma$ there is assigned a Laurent polynomial $P(\Gamma)\in\mathbb{Z}_{+}[q,q^{-1}]$ with positive integer coefficients, the Kuperberg invariant of $\Gamma$. Polynomial $P(\Gamma)$ is a quantization of the number of Tait colorings $t(\Gamma)$, so that $t(\Gamma)=P(\Gamma)|_{q=1}$. 

\begin{figure}%[H]
	\begin{center}
		{\includegraphics[width=230pt]{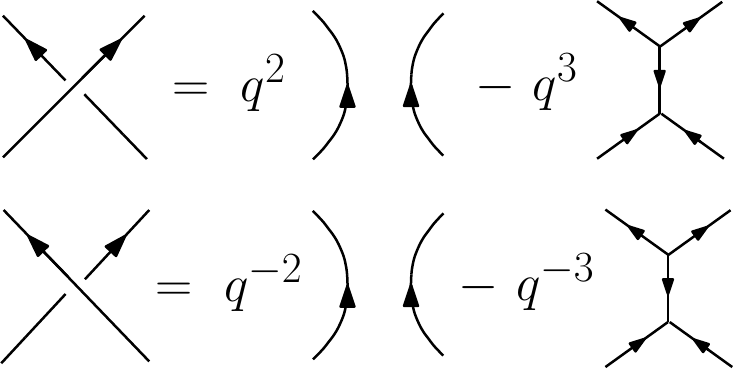}}
		\caption{Resolving crossings into linear combinations of webs in the Kuperberg invariant.}
		\label{3_04}
	\end{center}
\end{figure}

Kuperberg web invariant extends to the quantum $\SL(3)$ invariant of links. The crossings of a link diagram are resolved into linear combinations of web diagrams as depicted in Figure \ref{3_04}.
The invariant $P(D)\in \Z[q,q^{-1}]$ of a link diagram $D$ is defined as the linear combination of the corresponding web invariants. 
One can easily check that these relations are invariant under the Reidemeister moves, and the invariant $P(D)$ depends only on the underlying oriented link $L$ and can be denoted $P(L)$. 

While the integrality of $P$ is retained when passing from web to links, coefficients of $P(L)$ are no longer non-negative, unlike those of $P(\Gamma)$, and the positivity property of the coefficients does not extend from webs to links. 
%However, $P(L)\in\mathbb{Z}[q,q^{-1}]$ for a link $L$, and the link invariant no longer has the positivity property. 

\vspace{0.07in}

The original motivation for considering foams was to categorify the Kuperberg invariant for links, that is, to lift it to a homology theory for links. At the time a categorification of the Jones polynomial was already known~\cite{Kh2}.  The Jones polynomial~\cite{J} can be described via the representation theory of quantum $\slt$~\cites{KR,Tu,RT,CFS}. Kauffman's approach to the Jones polynomial uses a version of Figure~\ref{3_04} relations for $\slt$, where instead of a trivalent graph as the third term of each relation there is just a pair of arcs~\cite{K}. Consequently, no graphs are needed and the invariant reduces to that for a collection of circles in the plane. Each circle evaluates to $[2]=q+q^{-1}$. After categorification, $[2]$ becomes a graded Frobenius algebra $A$ of rank $2$ over the ground ring, and the Frobenius algebra structure is enough to build a complex to categorify the Jones polynomial. This results in a bigraded link homology theory, functorial for link cobordisms. 

For $\slthree$, collections of circles in the plane are replaced by oriented $\SL(3)$ webs. Their invariant $P(\Gamma)$ lies in $\Z_+[q,q^{-1}]$, indicating that its categorification should have a single grading, corresponding to powers of $q$, so that $P(\Gamma)$ is just the graded dimension, a polynomial with non-negative integer coefficients. These graded vector spaces for various $\Gamma$'s should then be assembled into complexes, and Figure~\ref{3_04} relations should lift to realizing the complex of a crossing as the cone of a map between the complexes for the two planar webs. This map may be induced by a natural cobordism between the graph in Figure~\ref{3_04} on the right with two vertices and the graph shown in the middle of Figure~\ref{3_04} which consists of two arcs. This cobordism will be given by a foam with a seam arc that connects the two vertices of the graph, see Figure~\ref{3_13}.

\subsection{Categorification of the Kuperberg invariant for planar webs}\label{subsec_Kuperberg}\hfill\vspace{0.02in}\\
We are looking to categorify $P(\Gamma)\in\mathbb{Z}_{+}[q,q^{-1}]$. More precisely, we want to upgrade $P(\Gamma)$ to homology groups $\mmH^\star(\Gamma)$ carrying a single $\mathbb{Z}$-grading, so that 
\begin{equation*}
	P(\Gamma)=\mathsf{grk}(\mmH^\star(\Gamma)),
\end{equation*}
where for a graded abelian group $V=\oplus_i V^i$ the graded rank 
$$
\mathsf{grk}(V):=\sum_i \mathsf{rk}(V^i)q^i. 
$$

A cobordism between webs $\Gamma$ and $\Gamma^\prime$ should induce a homomorphism between $\mmH^\star(\Gamma)$ and $\mmH^\star(\Gamma^\prime)$, and these homology groups and homomorphisms should give a functor from the category of web cobordisms (foams) to the category of graded abelian groups. 

Ideally, this functor would also be a TQFT on the category of foams, in particular being multiplicative on the disjoint union of webs:
\[
\mmH^{\star}(\Gamma\sqcup \Gamma')\cong \mmH^{\star}(\Gamma)\sqcup \mmH^{\star}(\Gamma'),
\]
extending multiplicativity of the Kuperberg invariant, 
\[
P(\Gamma\sqcup\Gamma')=P(\Gamma) P(\Gamma').
\]
This homology theory for webs and foams can be described by the universal construction. We wil write either $\mmH^{\star}(\Gamma)$ or $\brak{\Gamma}$ for the state space of a web $\Gamma$ that we are looking to construct. 

\vspace{0.07in}

\begin{figure}[H]
	\begin{center}
		{\includegraphics[width=100pt]{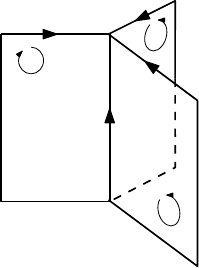}}
		\caption{Compatible orientations of facets near a seam of an oriented $\SL(3)$ foam.}
		\label{3_05}
	\end{center}
\end{figure}

We want to define oriented $\SL(3)$ foams so that their generic cross-sections by planes are oriented planar webs. As was mentioned previously, such webs have in-vertices and out-vertices, and the edges are oriented from an out-vertex to an in-vertex. 
Going one dimension up, we consider 2-dimensional complexes embedded in $\R^3$ such that every facet carries an orientation. Three facets can meet along a seam, so that any pair of facets along a seam has opposite orientations, as shown in Figure~\ref{3_05}. 
The union of any two facets in that figure is a disk, and the orientations of the two halves of the disk are opposite. This orientation condition extends one dimension up the oriented $\SL(3)$ web condition on orientation of edges at a vertex (all in or all out). 

Thus, faces of a foam are oriented in a way which is consistent with the orientation of their cross-sections. We pick the convention for induced orientation of a seam of an oriented $\SL(3)$ foam as shown in Figure~\ref{3_05}.
Next, consider a vertex of an unoriented foam, as in Figure~\ref{foam} on the right. There is no way to orient the six facets at the vertex so that along each of the four seams the orientations are compatible as in Figure~\ref{3_05}. Indeed, the link of the vertex is a trivalent graph on the two-sphere which is the complete graph $K_4$ on the four vertices. It is not bipartite, and carries no orientation of edges so that each vertex is either an in or an out vertex. 

This observation tells us that \emph{oriented} $\SL(3)$ foams, unlike unoriented ones, do not have vertices. It is impossible to introduce an orientation on facets near a vertex consistent with our restriction on orientation of facets along a seam of a foam. 

Absence of vertices greatly simplifies the theory of oriented $\SL(3)$ foams compared to unoriented ones. 
Singularities of oriented $\SL(3)$ foams consist of seams only. If a foam is closed, its set of singular points is a disjoint union of circles. If a foam $F$ has boundary, its singular points (seam points) are a union of circles and edge seams. The latter connect two boundary points of $F$, which are singularities (vertices) of the boundary web of the foam. 
     
For now, restrict to closed foams without boundary and consider a circular seam in a closed foam $F$. There are three facets attached along a seam, and one can check that it is impossible to create a non-trivial permutation of these facets as one goes along a circle~\cite{Kh1}. In other words, a circular seam of a foam has a neighborhood homeomorphic to the direct product of a circle and a tripod (a neighborhood of a vertex in a trivalent graph). 

\begin{figure}[H]
	\begin{center}
		{\includegraphics[width=200pt]{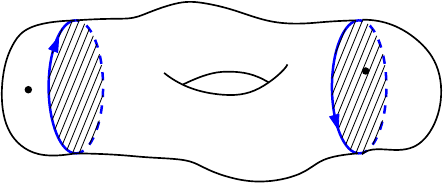}}
		\caption{A closed oriented foam. Singular circles (seams) are highlighted in blue. Two of the facets are shaded for convenience of presentation. Fixing an orientation along a singular circle determines the orientation for all facets and all circles, since this foam is connected.}
		\label{3_06}
	\end{center}
\end{figure}

Example of a closed (oriented $\SL(3)$) foam is depicted in Figure \ref{3_06}. Facets of an oriented foam may carry dots, just as in the unoriented case. 

\vspace{0.07in}

Let us construct a TQFT for such foams. Part of this construction lifts the quantum invariant of a graph to a corresponding graded module. Invariant of the empty graph $\emptyset$ is $1$, so to $\emptyset$ we assign the abelian group $\mathbb{Z}$, placing it in degree $0$.\footnote{Instead of $\Z$ any commutative ring $\Bbbk$ can be used as the ground ring, and we do not need to specialize to $\textrm{char}\Bbbk=0$ as in the unoriented foam case.} Invariant of a circle is given by
\begin{equation*}
   [3]=q^{-2}+1+q^2, 
\end{equation*}
which is, up to an overall grading shift, the graded dimension of
\begin{equation*}
    \mmH^\star(\mathbb{C}P^2)=\mathbb{Z}[x]/(x^3),
\end{equation*}
since
\begin{equation*}
    \textrm{deg}(1)=0,\quad \textrm{deg}(x)=2,\quad \textrm{deg}(x^2)=4.
\end{equation*}
The next web in complexity is the $\Theta$-web. Using relations in Figure \ref{3_02} to reduce a subregion of valency two and then to reduce a circle, we see, that the invariant of this graph is 
$[3] [2]$, 
which is,  up to a grading shift, the graded dimension of the flag variety of $\C^3$:
\begin{equation*}  
    \mmH^\star(\Fl_3)=\mathbb{Z}[x_1,x_2,x_3]/(E_1,E_2,E_3).
\end{equation*}
Now let us indicate why lifting invariants of the circle and the $\Theta$-web to these cohomology rings is a correct guess.

\begin{figure}[!htbp]
	\begin{center}
		{\includegraphics[width=170pt]{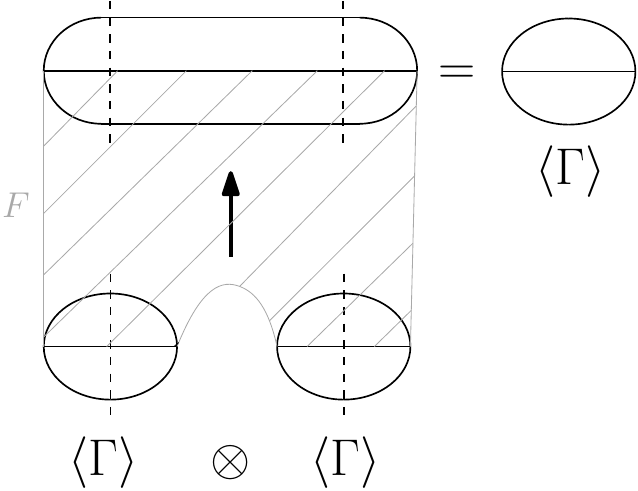}}
		\caption{Cobordism $F$ merging $\Gamma\sqcup \Gamma$ into  $\Gamma$ is indicated in grey. This cobordism is defined for any web $\Gamma$ with a symmetry axis.}
		\label{3_07}
	\end{center}
\end{figure}

\begin{figure}[H]
	\begin{center}
		{\includegraphics[width=80pt]{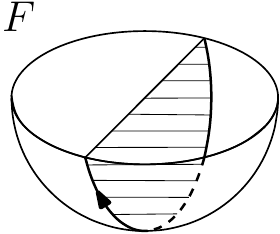}}
		\caption{Foam that represents the unit element in the ring $\langle \Gamma \rangle$, see also Figure~\ref{gens} and the discussion there.}
		\label{3_08}
	\end{center}
\end{figure}

\begin{figure}[H]
	\begin{center}
		{\includegraphics[width=100pt]{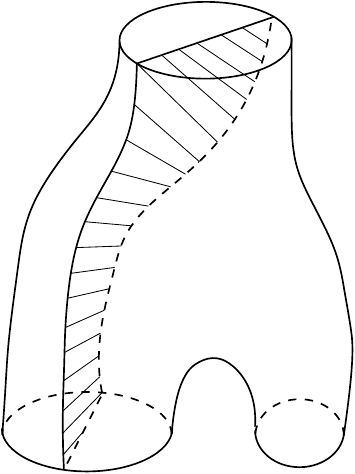}}
		\caption{Merging a circle into one of the three edges of the $\Theta$-graph.}
		\label{3_09}
	\end{center}
\end{figure}

We expect functoriality for the state spaces of webs. That is, a foam with boundary, viewed as a cobordism between two webs, should induce a homomorphism between their homology (state spaces). For a web $\Gamma$ with an axis of symmetry, such as the circle and the $\Theta$-web, part of that functoriality is the structure of a Frobenius algebra on the state space $\brak{\Gamma}$. The foam inducing multiplication on $\brak{\Gamma}$ is shown in Figure \ref{3_07}, where two copies of $\Gamma$ merge into one. 
Thus  the state space $\langle \Gamma \rangle$ of a web $\Gamma$ with a symmetry axis is endowed with a natural associative algebra structure under this multiplication operation
\begin{equation*}
    \langle \Gamma \rangle\otimes \langle \Gamma \rangle\to \langle \Gamma \rangle.
\end{equation*}
 This associative algebra has the unit element, given by foam $F$ depicted in Figure \ref{3_08}.
Furthemore, Figure~\ref{3_08} foam reflected in the horizontal plane gives a cobordism from $\Gamma$ to the empty web and induces a map $\epsilon:\brak{\Gamma}\lra \Z$ that turns $\brak{\Gamma}$ into a Frobenius algebra (over $\Z$). 
Hence, the state space of the circle $S$ and the $\Theta$-web are naturally Frobenius algebras, given a TQFT for foams. 
Moreover, these two Frobenius algebras are related. The cobordism in Figure \ref{3_09} merges a circle into one of the sides of the $\Theta$-web, making $\brak{\Theta}$ a module over $\brak{S}$. 
An example of a pair of graded Frobenius algebras with graded dimensions $[3]$ and $[2][3]$ (up to overall shifts) is given by the cohomology of $\CP^2$ and that of the flag variety $\Fl_3$ of $\C^3$, respectively. 
The latter is the space of full flags
\begin{equation*}
    0\subset L_1 \stackrel{}{\subset} L_2 \stackrel{}{\subset} \mathbb{C}^3,
\end{equation*}
where $\dim L_i =i$, $i=1,2$. 
Pick a Hermitian metric on $\C^3$ and form the triple $(L_1,L_1',L_1'')$ of mutually-orthogonal lines in $\C^3$. Here  $L_1'$ is the orthogonal complement of $L_1$ in $L_2$, and $L_1''$ the complement of $L_2$ in $\C^3$. 
One obtains three different maps of 
${\Fl}_3$ onto ${\CP}^2$ by remembering only one line in the above triple. 

These three maps induce maps on cohomology group 
$\mmH^\star(\mathbb{C}P^2)\lra \mmH^\star(\Fl_3)$,
making $\mmH^\star(\Fl_3)$ an $\mmH^\star(\mathbb{C}P^2)$-module in three different ways and giving 
us three multiplication maps 
\begin{equation}\label{eq_3_maps}
    \mmH^\star(\Fl_3) \otimes \mmH^\star(\mathbb{C}P^2) \to \mmH^\star(\Fl_3).
\end{equation}
Having these maps on the state spaces is consistent with the foam depicted in Figure~\ref{3_09} and its variations of merging the circle into two of the other facets of $\Theta\times [0,1]$.  

On the level of cohomology, maps \eqref{eq_3_maps} send generator $x$ of $\mmH^\star({\CP}^2)$ to one of the generators $x_i$ of $\mmH^\star(\Fl_3)$, via  $1\otimes x\mapsto x_i$. Categorifying the quantum invariant $[3]$ of the  circle web to cohomology ring $\mmH^\star({\CP}^2)$ and quantum invariant $[3][2]$ of the $\Theta$-graph to $\mmH^\star(\Fl_3)$ perfectly matches conditions on these state spaces provided by the cobordisms.\par

So far we have a commutative Frobenius algebra for the state space of a circle, a commutative Frobenius algebra for the state space of a $\Theta$-web, and their interactions in \eqref{eq_3_maps}. 

We define a homogeneous non-degenerate trace on $\mmH^\star({\CP}^2)$ by
\begin{equation}\label{eq_tr_cp}
\varepsilon(1)=0,\quad\varepsilon(x)=0,\quad\varepsilon(x^2)=-1,
\end{equation}
and on $\mmH^\star(\Fl_3)$ by the formula \eqref{eval} below. It turns out that we do not need to consider any other webs or Frobenius algebras to do the universal construction in the oriented $\SL(3)$ foams case. 

\begin{figure}[H]
	\begin{center}
		{\includegraphics[width=250pt]{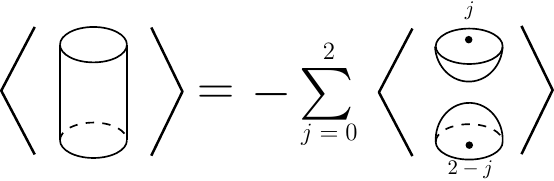}}
		\caption{The neck-cutting relation.}
		\label{3_10}
	\end{center}
\end{figure}

Trace \eqref{eq_tr_cp} leads to a surgery formula, a decomposition of the identity map from a circle to itself into a sum of basis elements (of the circle's state space) times their dual basis elements, depicted in Figure \ref{3_10}. 
There is a similar surgery formula for a 2-dimensional TQFT where the tube, representing the identity map, is decomposed as the sum over cup and cap cobordisms decorated by a basis element, respectively the dual basis element, for a choice of basis in a commutative Frobenius algebra. In the above example, the pair of dual bases are  $(1,x,x^2)$ and $(-x^2,-x,-1)$.

Given a closed foam $F$, possibly with dots, it may have some number of seam circles and no other singularities. We would like to evaluate $F$ to an element $\brak{F}$ of the ground ring $\Z$. The idea is to use the above surgery formula to split off each seam circle from the rest of the foam. A seam circle has three facets attached to it. On each of them choose a circle, consisting of smooth points in the interior of the facet and given by pushing off the seam circle into the facet. Take the foam out of $\R^3$ and do a surgery on it along each of these three circles. Repeat for all seam circles, with the total number of surgeries $3s$, where $s$ is the number of seam circles. 

The result is a sum over many terms, where each term is, up to a sign, a foam which is the disjoint union of connected foams of two types: 
\begin{itemize}
    \item A connected surface $S_{g,n}$ of some genus $g\ge 0$ with some number $n$ of dots,  
    \item A $\Theta$-foam with a singular seam, three disk facets attached to it, and $a,b,c$ dots placed on its three facets, see Figure~\ref{evaluate} on the right. 
    \begin{figure}[!htbp]
	\begin{center}
		{\includegraphics[width=250pt]{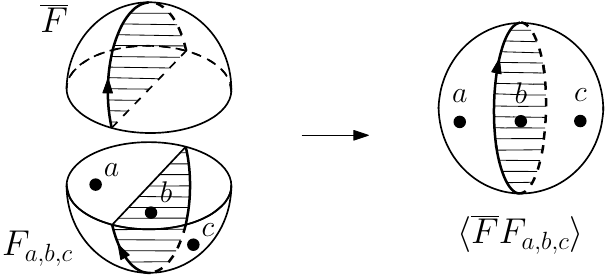}}
		\caption{Foam on bottom left represents a monomial $x_1^ax_2^bx_3^c\in \mmH^{\ast}(\Fl_3)$. Foam in on the top left describes the trace map $\varepsilon:\mmH^{\ast}(\Fl_3)\lra\Z$. Closed foam $\overline{F}F_{a,b,c}$ on the right evaluates to $\varepsilon(x_1^ax_2^bx_3^c)\in \Z$.}
        \label{evaluate}
	\end{center}
\end{figure}
\end{itemize}
A foam of the first type (connected surface, possibly with dots) is evaluated using the Frobenius algebra $\mmH^\star(\mathbb{C}P^2)$ with the trace~\eqref{eq_tr_cp}. 
Namely, 
\[
\brak{S_{g,n}} = 
\begin{cases}
    3 & \mathrm{if} \ g=1,n=0,  \\
    -1 & \mathrm{if} \ g=0,n=2, \\
    0 & \mathrm{otherwise}.
\end{cases}
\]
This 2D TQFT evaluates a dotless torus to $3$ (the rank of free abelian group $\mmH^\star(\mathbb{C}P^2)$), a sphere with two dots to $-1$, with all other evaluations $0$.  

$\Theta$-foam evaluation can be interpreted via the trace on the cohomology ring of the full flag variety $\mmH^\star(\Fl_3)$, see Figure~\ref{evaluate}. The trace is homogeneous and gives an isomorphism $\varepsilon:\mmH^6(\Fl_3)\cong\Z$. On monomials in the generators $x_1,x_2,x_3$ it is defined by 
\begin{eqnarray}\label{eval}
\varepsilon(x_1^a x_2^b x_3^c)=
\begin{cases}
1,\quad\textrm{for $(a,b,c)= (0,1,2)$ up to cyclic permutation},\\
-1,\quad\textrm{for $(a,b,c)= (0,2,1)$ up to cyclic permutation},\\
0,\quad\textrm{for $(a,b,c)\neq (2,1,0)$}.
\end{cases}
\end{eqnarray}
In the top degree $6$  monomials that are nontrivial in the cohomology ring  
are $x_1^2x_2$ and its permutations $x_i^2x_j$ (recall that $\deg(x_i)=2$). Furthermore, $x_i^2x_j=-x_j^2x_i$ in $\mmH^6(\Fl_3)$, for $i,j\in\{1,2,3\},i\not=j$,  so an odd permutation of the indices adds a sign to the value of the trace on a monomial. 

Foam $F$ with boundary depicted in Figure \ref{3_08} represents the unit element of the cohomology ring $\mmH^\star(\Fl_3)$. Putting dots on facets of $F$  gives monomial elements in this ring, see Figure~\ref{evaluate} bottom left. Denote this decorated foam by $F_{a,b,c}$. Closing up the diagram via the reflected foam $\overline{F}$ produces a closed foam $\overline{F}F_{a,b,c}$ that evaluates to 
\[
\brak{\overline{F}F_{a,b,c}}:= \varepsilon(x_1^ax_2^bx_3^c).
\]
In this way, we define the evaluation of a $\Theta$-foam, possibly with dots, via the trace on the cohomology ring $\mmH^\star(\Fl_3)$.

As explained earlier, an evaluation of dotted $\Theta$-foams, dotted spheres, and the surgery formula depicted in Figure \ref{3_10} allows to evaluate any closed foam and apply the universal construction. 
 
\begin{example}
Consider the foam in Figure~\ref{3_06} and also copied below.
\begin{figure}[H]
	\begin{center}
		{\includegraphics[width=200pt]{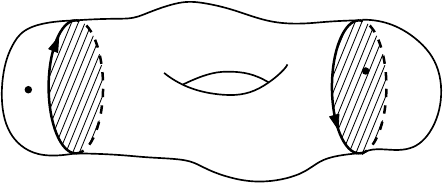}}
		\caption{}
		\label{3_11_1}
	\end{center}
\end{figure}
Let us compute its evaluation, where we omit writing the brackets for convenience. In the algorithm given earlier, a foam is surgered along three circles given by pushing off a seam circle into the three surrounding facets. Instead, we can be thriftier, without changing the evaluation of the foam. First,  use the neck-cutting relation to reduce the genus of each facet to $0$, reducing Figure~\ref{3_11_1} foam to 
\begin{figure}[H]
	\begin{center}
		{\includegraphics[width=250pt]{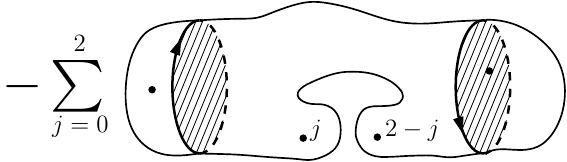}}
		%\caption{}
		%\label{3_11_2}
	\end{center}
\end{figure}
Two dots with weights $w_1,w_2$ floating on the same facet can be combined into a dot with the weight $w_1+w_2$  on the same facet, simplifying the sum to 
\begin{figure}[H]
	\begin{center}
		{\includegraphics[width=160pt]{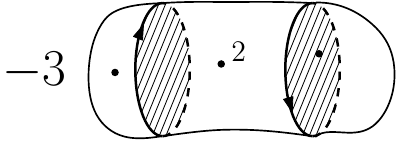}}
		%\caption{}
		%\label{3_11_3}
	\end{center}
\end{figure}
Applying the neck-cutting relation to the annulus facet gives a sum of $\Theta$-foams with dots:  
\begin{figure}[H]
	\begin{center}
		{\includegraphics[width=200pt]{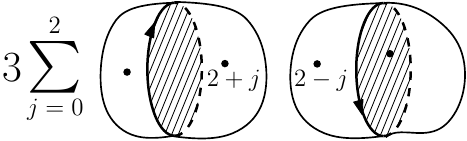}}
		%\caption{}
		%\label{3_11_4}
	\end{center}
\end{figure}
The evaluation is 
$$ \brak{F}=
3\,\sum\limits_{j=0}^{2}\varepsilon(x_2^{2+j}x_3)\varepsilon(x_1 x_2^{2-j}) =3\,\varepsilon(x_2^{2}x_3)\varepsilon(x_1 x_2^{2}) =-3,
$$
where, via (\ref{eval}), only the term $j=0$ survives.
\end{example}

\vspace{0.07in}

Using the rules above, an arbitrary closed foam $F$ in $\R^3$ can be evaluated to an integer $\brak{F}$. We apply the universal construction, described in Section~\ref{subsec_uc}, to this evaluation. The universal construction assigns an abelian group $\brak{\Gamma}$ to an oriented $\SL(3)$ web $\Gamma$. A foam $F$ with boundary, viewed as a cobordism from a web $\Gamma_0$ to a web $\Gamma_1$, induces a homomorphism of groups $\brak{F}:\brak{\Gamma_0}\lra\brak{\Gamma_1}$. 

Abelian group $\brak{\Gamma}$ is generated by foams $F$ with top boundary only, $\partial F =\Gamma$. Foam $F$ has degree given by the same formula \eqref{foam_degree} as in the unoriented case, where $\sum_{i<j}\chi_{ij}(F)$ is determined by any Tait coloring of $F$. Proposition~\ref{prop_graded} and Corollary~\ref{cor_graded} hold in the oriented case as well, with the ring $R$ replaced by $\Z$ (see~\cite{MV} for the construction of oriented $\SL(3)$ homology over the ground ring $\Z[E_1,E_2,E_3]$ of symmetric functions in 3 variables rather than over $\Z$). 

It is straightforward to derive direct sum decompositions of oriented webs that categorify Kuperberg's skein rules in Figure~\ref{3_02}. These are oriented analogues of Lemmas~\ref{skein2},~\ref{skein3} and~\ref{skein5}, see~\cite[Section 3.4]{Kh1}, for instance.   

The state spaces $\brak{\Gamma}$ of oriented webs are graded abelian groups. Decompositions of state spaces lift Kuperberg skein relations to isomorphisms of graded abelian groups. This construction lifts $P(\Gamma)\in\mathbb{Z}_{+}[q,q^{-1}]$ to a TQFT for webs and foams. 
\begin{theorem}\label{thm_or_TQFT}
$\langle \Gamma \rangle$ is a free graded abelian group with the graded rank 
\begin{equation*}
    \mathsf{grk}\, \langle\Gamma\rangle \ =\  P(\Gamma),
\end{equation*}
where $P(\Gamma)$ is the Kuperberg invariant of  oriented $\SL(3)$ web $\Gamma$. Functor $\brak{\ast}$ is a TQFT 
\[
 \mathsf{OFoam} \lra \Z\text{-}\mathsf{gmod} 
\]
from the category $\mathsf{OFoam}$ of oriented $\SL(3)$ webs and foams to the category of graded abelian groups. 
\end{theorem}\par
Objects of $\mathsf{OFoam}$ are oriented webs (in $\R^2$) and morphisms are oriented foams (in $\R^2\times [0,1]$) with oriented webs as  boundaries. Objects of $\Z\text{-}\mathsf{gmod}$ are graded abelian groups and morphisms are homogeneous homomorphisms of graded abelian groups. 

\vspace{0.07in}

Recall the state space $\brak{\Gamma}$ of an \emph{unoriented} web $\Gamma$ defined at the end of Section~\ref{subsec_uc}. Given an \emph{oriented} web $\Gamma$, denote by $\Gamma_{\mathsf{un}}$ this web considered as unoriented. Webs of the form $\Gamma_{\mathsf{un}}$ for some $\Gamma$ are exactly unoriented webs that are bipartite, or, equivalently, webs where every planar region has an even number of sides. For such $\Gamma$, there are two state spaces: 
\begin{itemize}
    \item State space $\brak{\Gamma}$, defined in the present section, which is a free graded $\Z$-module, of graded rank equal to the Kuperberg invariant $P(\Gamma)$. 
    \item State space $\brak{\Gamma_{\mathsf{un}}}$ defined in Section~\ref{subsec_uc}, which is a module over the ring $R$ of symmetric polynomials, see \eqref{eq_ring_R}, over a field $\Bbbk$ of characteristic 2. 
\end{itemize}
State space $\brak{\Gamma_{\mathsf{un}}}$ is a free graded 
$R$-module, of the same graded rank $P(\Gamma)$. The two state spaces are closely related, and, informally speaking, almost the same. 

The key point is that vertices are absent from oriented $\SL(3)$ foams but create high complexity for general unoriented $\SL(3)$ foams, making it hard to understand state spaces of unoriented $\SL(3)$ webs. One type of exception is a bipartite unoriented $\SL(3)$ web $\Gamma_{\mathsf{un}}$, coming from an oriented $\SL(3)$ web $\Gamma$. 
Any element in the state space $\brak{\Gamma_{\mathsf{un}}}$ can be realized by a linear combination of foams $F$ \emph{without vertices}, with $\partial F =\Gamma_{\mathsf{un}}$. This is due to the direct sum decompositions in Lemmas~\ref{skein2},~\ref{skein3},~\ref{skein5}. These direct sum decompositions mirror those in the oriented case. Thus, for bipartite unoriented webs, foams with vertices can be avoided when building their state spaces. 
%For oriented foams $F$ reduce evaluation defined in the present section (shortly before Theorem~\ref{thm_or_TQFT}) modulo two, so the ground ring of the theory is $\mathbf{F}_2$. Reduce evaluation in Section {\bf section number}.

\subsection{Categorification of the Kuperberg invariant for links}\label{subsec_cat_Ku}
\begin{figure}[!htbp]
	\begin{center}
		{\includegraphics[width=290pt]{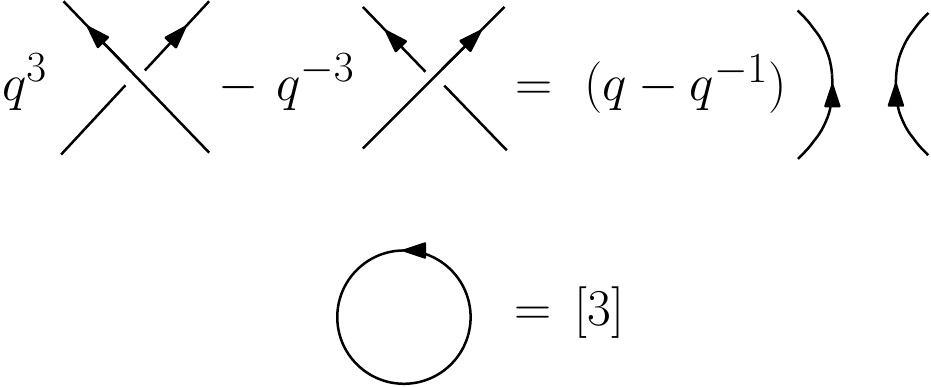}}
		\caption{Skein relation for the Reshetikhin-Turaev quantum $\SL(3)$ invariant, a.k.a. the Kuperberg link invariant, and its normalization. This invariant is also a  specialization of the HOMFLYPT 2-variable invariant, via $(t,q)\mapsto (q^3,q)$.}
		\label{3_12}
	\end{center}
\end{figure}

\begin{figure}[!htbp]
    \centering 
    \includegraphics[width=340pt]{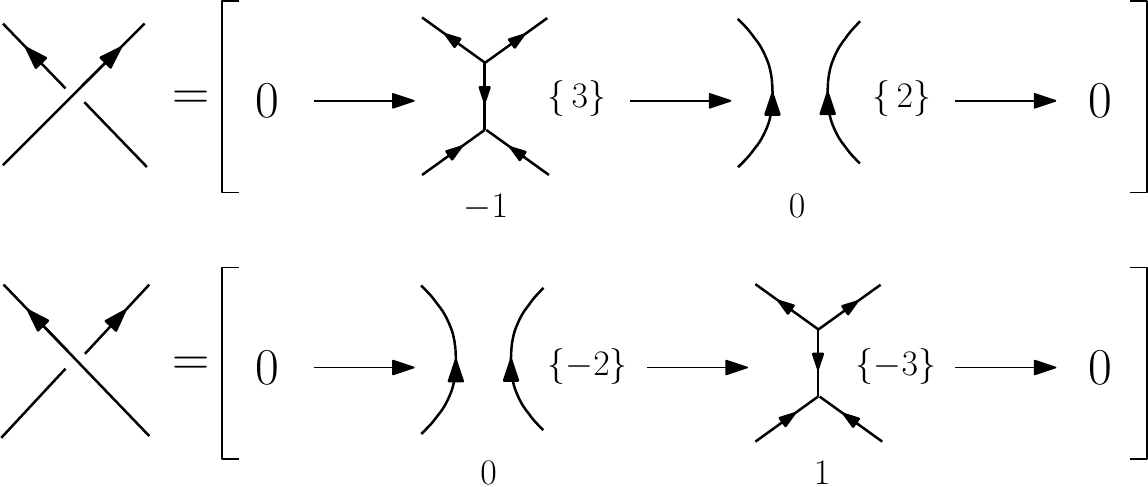} 
    \caption{Converting the crossing into its web resolutions and forming the cone of the map between corresponding complexes induced by the foam $F$ or $\overline{F}$.} 
    \label{cones_fig}
\end{figure}

\hfill\vspace{0.02in}\\
The Kuperberg invariant of a link $P(L)\in\mathbb{Z}[q,q^{-1}]$ is described by the skein relations and the normalization, as in Figure~\ref{3_12}.
\noindent 
One can check that these relations follow from the ones in Figure~\ref{3_04} and web skein relations in Figure~\ref{3_02}. Figure~\ref{3_04} relations can be converted into long exact sequences of homology groups coming from the iterated cone construction in Figure~\ref{cones_fig}, as follows. 

\begin{figure}[H]
	\begin{center}
		{\includegraphics[width=380pt]{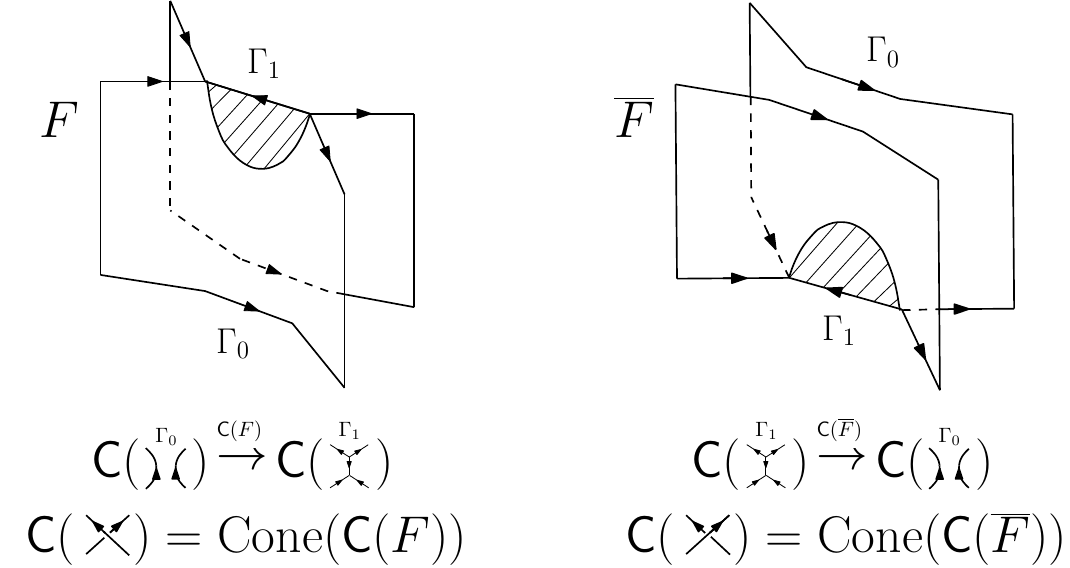}}
		\caption{Foam $F$ between two adjacent resolutions $\Gamma_0,\Gamma_1$ of $D$ induces the map $\brak{F}:\brak{\Gamma_0}\lra \brak{\Gamma_1}$. Foam $\overline{F}$ is the reflection of $F$ in a horizontal plane, inducing the map $\brak{\overline{F}}:\brak{\Gamma_1}\lra\brak{\Gamma_0}$. Summing over all resolutions of the remaining crossings gives corresponding maps between the complexes.}
		\label{3_13}
	\end{center}
\end{figure}

Pick a generic planar projection $D$ of an oriented link $L$ with $n$ crossings. Each of the $n$ crossings admits two resolutions into webs with four boundary points. One of these resolutions consists of two arcs, the other is a web with two trivalent vertices. There are $2^n$ possible resolutions of $D$ into webs. There are canonical cobordisms between pairs of resolutions $\Gamma_0,\Gamma_1$ that differ near one crossing of $D$, via the foam $F$ shown in Figure~\ref{3_13} and its reflection $\overline{F}$. Foam $F$ is a cobordism from $\Gamma_0$ to $\Gamma_1$, foam $\overline{F}$ is a cobordism from $\Gamma_1$ to $\Gamma_0$. These cobordisms induce maps between state spaces $\brak{\Gamma_0}$ and $\brak{\Gamma_1}$ in both directions. 

Place state spaces $\brak{\Gamma}$ of all $2^n$ resolutions $\Gamma$ of $D$ into vertices of an $n$-dimensional cube $\mathsf{Cu}(D)$. Relate state spaces $\brak{\Gamma_0},\brak{\Gamma_1}$ in adjacent vertices by the map induced by $F$ if the crossing is negative and by $\overline{F}$ if the crossing is positive. Among the two crossings on the left of  Figure~\ref{cones_fig}, the top crossing is positive and the bottom one is negative. 

Internal grading (also called $q$-grading) of state spaces can be shifted so that the induced maps have degree $0$. Choice of resolution of a crossing contributes $\pm 2$ or $\pm 3$ to the total shift as shown in Figure~\ref{cones_fig}.

The resulting cube $\mathsf{Cu}(D)$ carries a graded abelian group in each vertex. Each oriented edge carries a grading-preserving homomorphism of groups, and each square facet is commutative. Adding minus signs so that each square anticommutes, collapse this cube into the complex $\mathsf{C}(D)$ of graded abelian groups and position this complex starting in the negative homological degree $-n_+$, where $n_+$ is the number of positive crossings of $D$. Contributions of local resolutions to homological degrees of individual terms are also shown in Figure~\ref{cones_fig}, by the numbers $0,\pm 1$ under the diagrams. The rightmost term of the complex is in positive homological degree $n_-$, the number of negative crossings of $D$. 

Complex $\mathsf{C}(D)$ can also be defined inductively on the number of crossings for planar projections that contain both crossings and trivalent vertices (projections of spatial webs). Induction step is encoded by the rules in Figure~\ref{cones_fig}, where the maps of complexes are induced by the foam $F$ or its reflection $\overline{F}$.  

\begin{theorem}\label{thm_invariance}
    Reidemeister moves of oriented link diagrams induce natural homotopy equivalences of complexes $\mathsf{C}(D)$ of graded abelian groups and isomorphisms of their homology groups $\mathsf{H}(D)$.  
\end{theorem}

Define link homology $\mathsf{H}(L)$ as $\mathsf{H}(D)$ for a diagram $D$ of an oriented link $L$. Isomorphism classes of bigraded groups are invariants of $L$, and with more work one can establish that the groups themselves (and not just their isomorphism classes) are naturally associated to links.

%{\color{teal} ORIGINAL VERSION: Assume that there is already a chain complex assigned to any diagram with $n-1$ crossings and any number of trivalent vertices. There is a canonical cobordism between two diagrams on the right-hand side of these relations, which induces a map between state spaces of such diagrams in a similar fashion to how the unoriented case was worked out. This, in turn, induces a map between chain complexes associated to such diagrams. The chain complex assigned to the diagram on the left-hand side of relations \ref{3_04} is a cone of a map given by the foam. This setting is shown in the Figure \ref{3_13}. Similarly for the opposite crossing in the left-hand side of relations \ref{3_04}, but the order is reversed and the cobordism is turned upside down.}

%{\color{teal} One takes a corresponding cone of a map induced by this foam and for multiple crossings inductively builds a complex as an iterated cone of maps of complexes for simpler diagrams with fewer crossings. The homology of a resultant iterated cone gives the desired link homology $\mmH^\star(L)$. }

In the bigrading decomposition 
\begin{equation}\label{eq_SL3_homology}
    \mathsf{H}(L)=\bigoplus_{i,j}\mathsf{H}^{i,j}(L)
\end{equation}
 $i$ corresponds to the homological degree and $j$ to the $q$ degree. This homology is a categorification of $P(L)$ in the sense that 
\begin{equation*}
	P(L) =\sum_{i,j}(-1)^i q^j\textrm{rk}(\mathsf{H}^{i,j}(L)),
\end{equation*}
where $P(L)$ is the Kuperberg invariant of links, which is also the quantum $\SL(3)$ Reshetikhin-Turaev invariant for the fundamental representation. Replacing $q^{\pm 3}$ by $q^{\pm N}$ in the Figure~\ref{3_12} skein relation and suitably modifying the other skein relations results in the quantum $\SL(N)$ Reshetikhin-Turaev quantum invariant~\cite{RT} and a one-variable specialization of the HOMFLYPT polynomial~\cites{HOMFLY,PT}. For $N=2$ none of the difficulties related to constructing graded state spaces of webs and dealing with foams appear, since in this case planar diagrams reduce to collections of circles. 

More information and references on the $\SL(3)$ link homology theory $\mathsf{H}(L)$ can be found at the end of Section~\ref{subsec_sl3web}.

\begin{remark}
$\SL(3)$ foams in the oriented case are easier to deal with compared to the unoriented $\SL(3)$ foams. Closed oriented $\SL(3)$ foams do not have vertices, only singular seams that close up into singular circles. This allows for the evaluation of such foams using a system of two commutative Frobenius $R$-algebras to evaluate dot-decorated surfaces and dot-decorated $\Theta$-foams. The two algebras are the cohomology rings $\mmH^\star_{U(3)}(\mathbb{C}P^2)$ and $\mmH^\star_{U(3)}(\Fl_3)$, respectively, or their non-equivariant counterparts  $\mmH^\star(\mathbb{C}P^2)$ and $\mmH^\star(\Fl_3)$. In Section~\ref{sec_oriented} we explained foam evaluation and the construction of link homology in the latter, non-equivariant case. In the equivariant case~\cite{MV} state spaces of webs and cohomology groups of links are modules over the ring $\Z[x_1,x_2,x_3]^{S_3}$ of symmetric functions in 3 variables, which is the integral version of the ring $R$ in \eqref{eq_ring_R}.    
\end{remark}
\section{\texorpdfstring{$\SL(3)$}{SL(3)}-web algebra} In the previous section we used oriented $\SL(3)$ webs and foam cobordisms between them to lift the Kuperberg invariant to a homology theory for links. This construction can be extended from links to tangles. For that, one needs to extend the state sum invariant and its categorification from trivalent oriented webs and foam cobordisms between them to trivalent oriented webs with boundary and their cobordisms, which are foams with boundary and corners. 

\subsection{$\SL(3)$ webs with boundary}\label{websbnd}\hfill\vspace{0.02in}\\
In this section we consider the category of webs with boundary modulo Kuperberg's relations. In this category hom spaces are free $\Z[q,q^{-1}]$-modules over the basis of non-elliptic webs, and the category is braided monoidal. Tensoring hom spaces with $\Q(q)$ results in a braided monoidal category which is equivalent (even isomorphic) to the category of intertwiners between tensor products of representations $V$ and $V^{\ast}$ of quantum $\slthree$.

\vspace{0.07in}

Let us define a category $\Web$ of $\SL(3)$ webs with boundary. Objects in $\Web$ are sequences of orientations 
$$
\unde=(\eps_1,\ldots,\eps_n),\quad\eps_j\in\{\pm 1\},
$$ 
and the set of such sequences, including the empty sequence $\varnothing$, is denoted by $\Seq$. We may alternatively just write a sequence of signs, omitting $1$'s from the list.  To any $\unde\in\Seq$ one can assign a sequence of dots on a horizontal line, thinking of it as part of the boundary of an oriented $\SL(3)$ web. Each dot is signed by "$+$" and "$-$" in accordance with $+1$ and $-1$ in $\unde$. 

Closed oriented $\SL(3)$ webs were defined in Section \ref{cow}. Now we consider these webs properly embedded in the strip $\R\times [0,1]$ with boundary sets defined as $\partial_j\Gamma:=\Gamma\cap (\mathbb{R}\times \{j\})$, $j=0,1$ and signs on points of the boundary set induced by orientations of edges of the web near the boundary in the standard way. We refer to $\partial_0\Gamma$, $\partial_1\Gamma$ as the lower and the upper boundary of $\Gamma$. Our webs are {\it properly embedded}, so that each boundary point is an endpoint of a single strand. Hence, we encode web boundaries by signed sequences $\unde$.

Next, we define morphisms in the category $\Web$. 
For two sequences $\unde,\undep\in\Seq$ the set of morphisms 
$$
\Web^{\undep}_{\unde}:=\Hom_{\Web}({\unde},{\undep}),
$$
consists of $\Z[q,q^{-1}]$-linear combinations of isotopy classes of webs $\Gamma$ with $\partial_0 \Gamma =\unde$ and $\partial_1 \Gamma=\undep$, modulo skein relations depicted in Figure \ref{3_02}. Let us also define $\Web^{\unde}:=\Web_\varnothing^{\unde}$. This is a $\Z[q,q^{-1}]$-module spanned by webs with the empty lower boundary and the upper boundary $\unde$. 

\begin{figure}[H]
	\begin{center}
		{\includegraphics[width=450pt]{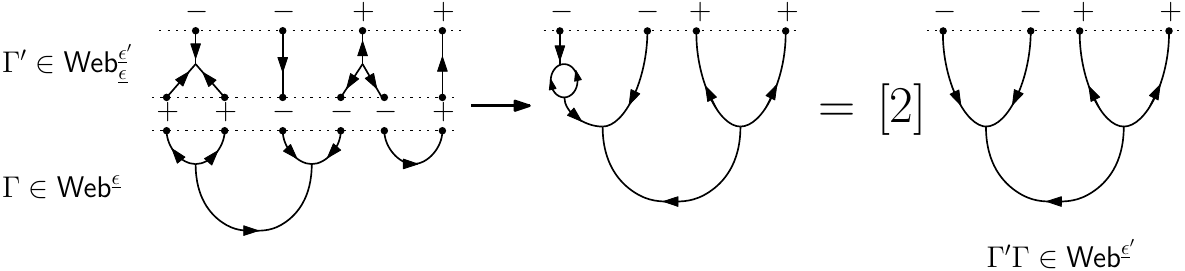}}
		\caption{An example of $\Gamma^\prime\in\Web^{\undep}_{\unde}$ and its composition with $\Gamma\in\Web^{\unde}$, which results in $\Gamma^\prime\Gamma\in\mathsf{Web}^{\undep}$. Here $\unde=(+,+,-,-,-,+)$ and $\undep=(-,-,+,+)$. The composition web is simplified using the second skein relation in Figure \ref{3_02}.}
		\label{morex}
	\end{center}
\end{figure}

Webs $\Gamma,\Gamma^\prime$ with $\partial_0 \Gamma^\prime=\partial_1 \Gamma$ can be composed along the common boundary to the web $\Gamma^\prime \Gamma$. Passing to $\Z[q,q^{-1}]$-linear combinations defines the composition of morphisms in $\Web$. Example of a web $\Gamma^\prime$ corresponding to an element in $\Web^{\undep}_{\unde}$ and its composition with $\Gamma$ in $\Web^{\unde}$ is given in Figure \ref{morex}.
A web decomposes $\R\times [0,1]$ into regions, by cutting the plane strip along all edges of the web (or taking  connected components of the complement $\R\times [0,1]\setminus \Gamma$). 
A region of a web $\Gamma\in\Web^{\undep}_{\unde}$ is called \emph{inner} if it does not intersect the boundary $\mathbb{R}\times \{0,1\}$. For example, none of the regions of webs $\Gamma,\Gamma'$ in Figure~\ref{morex} is inner. We say that $\Gamma$ is \emph{non-elliptic} if all inner regions of $\Gamma$ have at least 6 sides. (Note that each inner region of $\Gamma$ has an even number of sides, due to the bipartite property of $\Gamma$.) We refer to~\cite{Kup} for a conceptual explanation of this terminology. Web depicted in Figure \ref{boundeqiv} is an example of a non-elliptic web in $\Web^{\unde}$ for $\unde=(+,-,+,+,+,+,-)$. It has one inner region, with 6 sides. Webs $\Gamma,\Gamma'$ in Figure~\ref{morex} are non-elliptic. 

\begin{figure}[H]
	\begin{center}
		{\includegraphics[width=250pt]{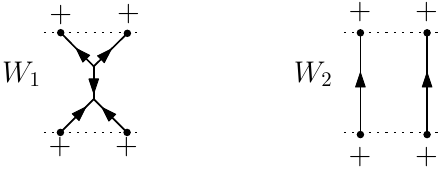}}
		\caption{For $\unde=\undep=(+,+)$, $\Web^{\undep}_{\unde}$  is a free $\Z[q,q^{-1}]$-module with a two-element basis $\B^{++}_{++}=\{W_1,W_2\}$.}
		\label{basisspan}
	\end{center}
\end{figure}

Kuperberg~\cite{Kup} establishes the following result. 
\begin{prop}
    $\Z[q,q^{-1}]$-module $\Web^{\undep}_{\unde}$ is free, with a basis of all non-elliptic webs with upper and lower boundary given by $\undep,\unde$ respectively. 
\end{prop}
Denote the basis of non-elliptic webs in $\Web^{\undep}_{\unde}$ by $\B^{\undep}_{\unde}$. Figure~\ref{basisspan} shows the basis $\B^{++}_{++}$. 

\begin{figure}[H]
	\begin{center}
		{\includegraphics[width=320pt]{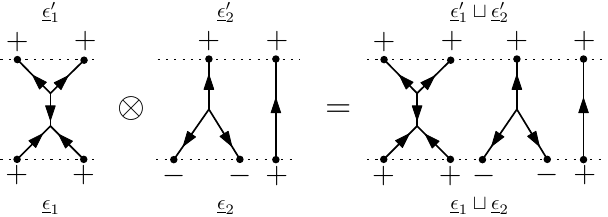}}
		\caption{An example illustrating the monoidal structure for webs $\otimes: \Web^{\undep_1}_{\unde_1}\otimes \Web^{\undep_2}_{\undep_2} \to \Web^{\undep_1\sqcup\undep_2}_{\unde_1\sqcup\unde_2}$.}
        \label{monoidal}
	\end{center}
\end{figure}

Category $\Web$ is monoidal. The monoidal structure is given on objects by concatenating two sign sequences, $\unde_1\otimes \unde_2:=\unde_1 \sqcup \unde_2$. The tensor product of webs is given by placing them in parallel, as depicted in Figure \ref{monoidal}, and then extending the tensor product to linear combinations of webs. Kuperberg skein relations are local, thus compatible with the tensor product of morphisms.

Category $\Web$ is pivotal. This means that each object has a dual, with the dual of a sequence $\unde$ given by the sequence $\unde^{\ast}$ which is $\unde$ written in the opposite order, with all signs reversed. Duality morphisms are given by diagrams of parallel half-circles and satisfy standard axioms corresponding to isotopies of oriented arcs in the plane.
\begin{figure}[H]
	\begin{center}
		{\includegraphics[width=300pt]{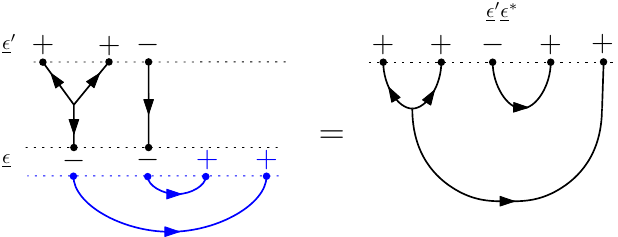}}
		\caption{Action of duality morphisms given by pivotal structure on $\Web^{\undep}_{\unde}$ for $\unde=(-,-)$, $\undep=(+,+,-)$. Here duality morphisms are emphasized in blue, and the resulting web gives an element of $\Web^{\undep\unde^\ast}$. Using pivotal structure one can rotate webs in any direction.}
        \label{pivotal}
	\end{center}
\end{figure}

\begin{remark}\label{morphisms_remark}
Using the pivotal structure, $(\unde,\undep)$-webs in the strip $\R\times [0,1]$ can be converted to $(\varnothing,\undep\unde^{\ast})$-webs in the lower half-plane. Soon we will see that the latter are central for a categorification of the Kuperberg tangle invariant. Due to a homeomorphism between a closed disk minus a boundary point and the lower half-plane, it is also possible to consider webs $\Web^{\unde}$ in a disk $D^2$, see Figure \ref{boundeqiv}. A category of such webs satisfies a set of axioms of a \textit{spider}, which was originally defined in \cite{Kup} as a graphical calculus for representation theory of $\Uqthree$.
\begin{figure}[H]
	\begin{center}
		{\includegraphics[width=200pt]{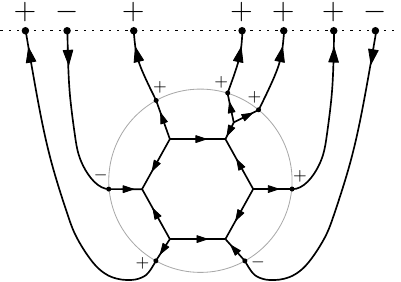}}
		\caption{An example of a web in the lower-plane $\mathbb{R}^2_{-}$ with boundary $\unde=(+,-,+,+,+,-,+)$.  This is a non-elliptic web with an inner region. Essentially the same web in the disk $D^2$ is bounded by the gray circle.}
		\label{boundeqiv}
	\end{center}
\end{figure}
\end{remark}

Category $\Web$ is braided, with the braiding structure on pairs of generating objects (signs) encoded in Figure~\ref{3_04}. That figure defines the braiding on $+\otimes +$, and rotating the figure gives the braiding on all pairs of signs $\pm \otimes \pm$. The braiding then extends to the tensor product of objects in the standard way.

Thus, $\Web$ is a braided pivotal monoidal category. 
This category encodes Kuperberg $\SL(3)$ webs with boundary and skein rules for them, and hom spaces in this category are free $\Z[q,q^{-1}]$-modules with  natural bases given by non-elliptic webs.

\vspace{0.07in} 

\subsection{Functor $\Ku$ from tangles to webs}\hfill\vspace{0.02in}\\
Define the category $\Tan$ of (oriented) tangles to have sign sequences $\unde$ as objects.  An (oriented) tangle $T$ is a proper smooth embedding of a compact oriented 1-manifold, possibly with boundary, into $\R^2\times [0,1]$. Such a manifold is a finite union of intervals and circles. By an $(n,m)$-tangle we mean a tangle with $n$ boundary points on $\R^2\times \{0\}$ and $m$ boundary points on $\R^2\times \{1\}$. We fix these endpoints to be standard in the plane, $\{(0,0),(1,0), \dots, (n,0)\}\subset \R^2$, and likewise for $m$. We consider oriented tangles up to rel boundary isotopy. Orientation of an $(n,m)$-tangle $T$ induces sign sequences $\partial_0 T = \unde$ and $\partial_1 T=\undep$ of length $n$ and $m$, respectively, by inducing orientations on the boundary points. 
\begin{figure}[!htbp]
	\begin{center}
		{\includegraphics[width=150pt]{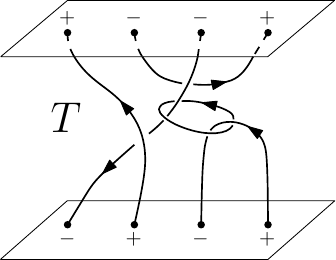}}
		\caption{An example of an oriented $(4,4)$-tangle $T\in\Tan$ with $\partial_0 T =(-,+,-,+)$ and $\partial_1 T =(+,-,-,+)$.}
	\end{center}
\end{figure}

Define the set of morphisms from $\unde$ to $\undep$ in $\Tan$ as oriented tangles $T\in \R^2\times [0,1]$ with $\partial_0 T$ described by $\unde$ and $\partial_1 T$ given by $\undep$. Rel boundary isotopic tangles define the same morphism. In particular, the set (abelian monoid) of endomorphisms of the empty sequence $\emptyset$ in $\Tan$ is the set of isotopy classes of oriented links in $\R^3$. 
Category $\Tan$ is pivotal braided monoidal. 

There is a functor $\Ku: \Tan\lra \Web$ which is the identity on objects. To define it on a morphism (tangle) $T$, project $T$ generically from $\R^2\times [0,1]$ onto the strip $\R\times [0,1]$ to get a diagram $D$ of $T$ with finitely many crossings. Resolve each crossing of $D$ via Figure~\ref{3_04} equations into into linear combinations of $\SL(3)$ webs. The result is a $\Z[q,q^{-1}]$-linear combination of $2^k$ webs, where $k$ is the number of crossings of $D$, which is a morphism in $\Web_{\partial_0 T}^{\partial_1 T}$. Define $\Ku(T)$ to be this morphism. We refer to~\cite{Kup} for details and the following result.  

\begin{prop}
    $\Ku$ is a well-defined functor $\Tan\lra \Web$ which respects pivotal braided monoidal structures of the two categories. 
\end{prop}
Functor $\Ku$ is the first of the three functors shown in \eqref{eq_3_functors} below. 

\vspace{0.07in} 

To relate webs to quantum group representations, it is convenient to extend the ground ring $\Z[q,q^{-1}]$ to the field of rational functions $\Q(q)$.
Define the category $\Web_{\Q(q)}$ to have the same objects as in the category $\Web$. The hom spaces in $\Web_{\Q(q)}$ are given by extending scalars from $\Z[q,q^{-1}]$ to $\Q(q)$: 
\begin{equation}
\Hom_{\Web_{\Q(q)}}(\unde,\undep) := \Hom_{\Web}(\unde,\undep)\otimes_{\Z[q,q^{-1}]} \Q(q).
\end{equation} 
This functor, shown as the middle arrow in \eqref{eq_3_functors}, is just an extension of scalars. 
Since the hom spaces in $\Web$ are free $\Z[q,q^{-1}]$-modules this is a faithful functor. 

\vspace{0.07in} 

Consider the category $\Rep'(\Uqthree)$ of finite-dimensional representations of the quantum universal enveloping algebra $\Uqthree$. The latter is defined over the ground field $\Q(q)$. Recall that $V$ denotes the fundamental 3-dimensional representation of $\Uqthree$ and $V^{\ast}$ denotes its dual. For explicit formulas for the quantum $\slthree$, these representations and generating intertwiners between their tensor products see~\cite{KhK}, for instance.  

For a sequence $\unde$ define $V(\unde)$ inductively by: 
\[
V(+)=V, \ \ V(-)=V^{\ast}, \ \ V(\unde_1\sqcup \unde_2)=V(\unde_1)\otimes V(\unde_2). 
\]
That is, $V(\unde)$ is the tensor product of copies of $V$ and its dual representation $V^{\ast}$ in the order of signs in $\unde$. 

Define $\Rep(\Uqthree)$ to be the full subcategory of $\Rep'(\Uqthree)$ with objects $V(\unde)$ over all sequences $\unde\in \Seq$. In other words, $\mathsf{Rep}(\Uqthree)$ is a full subcategory of the category of finite-dimensional representations of quantum $\slthree$ with objects given by tensor products of $V$ and $V^{\ast}$.
To a web $\Gamma$ there is a natural way to assign an intertwiner (a homomorphism of representations)  
\[
V(\Gamma): V(\partial_0 T)\lra V(\partial_1 T),
\]
as explained in~\cites{Kup,KhK}, for example. In particular, to a trivalent vertex one assigns a particular vector in the one-dimensional space of $\Uqthree$-invariants of $V^{\otimes 3}$. 

This assignment is compatible with the composition of intertwiners and extends to a functor 
\begin{equation}
     \Web_{\Q(q)} \lra \Rep(\Uqthree), 
\end{equation}
shown as the right arrow in \eqref{eq_3_functors}. 
Moreover, this functor is an equivalence of categories (even an isomorphism, since the category $\Rep(\Uqthree)$ has very few objects.) In particular, there are natural isomorphisms of $\Q(q)$-vector spaces
$$
\mathsf{Web}^{\undep}_{\unde}\otimes_{\Z[q,q^{-1}]}\Q(q)\cong\mathrm{Hom}_{\Uqthree}(V^{\otimes \unde}, V^{\otimes \undep }),
$$

\begin{corollary}
    Dimension of the space of intertwiners $\mathrm{Hom}_{\Uqthree}(V^{\otimes \unde}, V^{\otimes \undep })$ equals the number of $(\undep,\unde)$ non-elliptic webs. 
\end{corollary}

These constructions are summarized in the diagram of four categories and three functors between them: 
\begin{equation}\label{eq_3_functors}
\Tan \stackrel{\Ku}{\lra} 
\Web   \lra \Web_{\Q(q)} \stackrel{\cong}{\lra} \mathsf{Rep}(\Uqthree). 
\end{equation} 
All four categories are pivotal braided monoidal and the functors respect these structures. 
%\MK{Importance of invariants and homs from the empty sequence to $\unde$. Balanced sequences and invariant spaces (balanced sequences are discussed below). Need for categorification.}

\subsection{$\SL(3)$-web algebra}\label{subsec_sl3web}\hfill\vspace{0.02in}\\
Now let us extend the homology theory for links developed in earlier sections to a homology theory for tangles. For $\unde\in\Seq$ define the weight of a sequence as
$$
w(\unde):=\sum\limits_{j=1}^n \eps_j.
$$
We say that $\unde$ is \emph{balanced} if $w(\unde)\equiv 0\,(\mathrm{mod}\,3)$ and denote the set of balanced sequences by $\Bal$. In this section, we restrict our considerations to webs $\Gamma$ such that their boundaries are balanced, $\partial_0\Gamma, \partial_1 \Gamma\in \Bal$. If a closed web is cut in half by a generic line in the plane, the sequence of signs at the cut is balanced. Vice versa, any balanced sequence is the boundary of some web in a half-plane. Thus, a sequence $\unde$ is balanced iff there exists a web $\Gamma$ with $\partial_1\Gamma=\unde$ and $\partial_0\Gamma=\varnothing$. For a web $\Gamma$, its lower boundary sequence is balanced iff its upper boundary sequence is balanced. 

%Each oriented $\SL(3)$ web is constructed using elementary graphs: oriented arcs and trivalent out-vertices and in-vertices. They correspond to contractions of morphisms and intertwiners from Figure \ref{3_03}. Boundary sequences of these elementary graphs satisfy the condition $\sum\limits_{j=1}^n \eps_j=0\,(\mathrm{mod}\,3)$. Hence, any sequence $\unde$, which encodes boundary points of some oriented $\SL(3)$ web satisfies this condition. 

\vspace{0.07in}

For $\Gamma\in\Web^{\unde}$ define $\overline{\Gamma}\in\Web_{\unde}^{\emptyset}$ as the reflection of $\Gamma$ about a horizontal line, additionally reversing the orientation of arcs.
\begin{figure}[H]
	\begin{center}
		{\includegraphics[width=380pt]{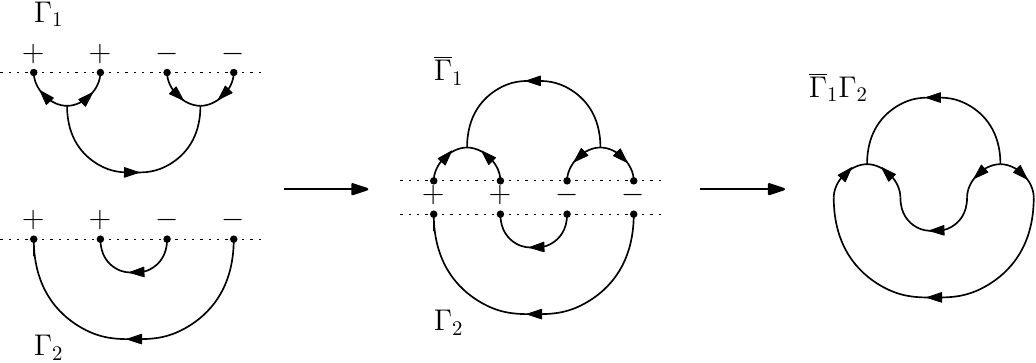}}
		\caption{An example of a graph $\overline{\Gamma}_j\Gamma_i$, which is a gluing of two objects in $\B^{\unde}$, here $\unde=(+,+,-,-)$.}
	\end{center}
\end{figure}

\vspace{0.07in}

As an intermediate object, for a balanced sequence $\unde$ introduce the category $\mcH^{\unde}$. Objects of $\mcH^{\unde}$ are non-elliptic webs $\Gamma$ with boundary conditions $\partial_1\Gamma=\unde$ and $\partial_0\Gamma=\varnothing$. We denote the set of such webs by $\B^{\unde}:= \B^{\unde}_{\varnothing}$ and think of them as non-elliptic webs in the lower half-plane with the boundary $\unde$. This set consists of finitely many elements, and we write $\Gamma_i,\Gamma_j$, etc. for its elements.
\begin{example}
	Consider the sign sequence $\unde=(+,+,-,-)$. There are two non-elliptic webs, $\B^{\unde}=\{\Gamma_1$, $\Gamma_2\}$ with boundary $\unde$, which are the generating objects. They can also be obtained by bending the two elements of the basis $\B^{++}_{++}$ in Figure~\ref{basisspan}. 
\begin{figure}[H]
	\begin{center}
		{\includegraphics[width=280pt]{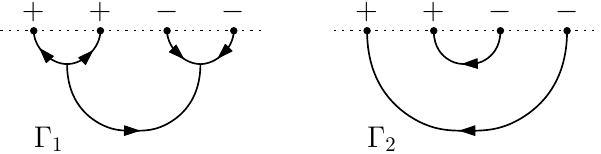}}
		\caption{Elements $\Gamma_1$, $\Gamma_2$ of $\B^{\unde}$, where $\unde=(+,+,-,-)$.}
		\label{4_00}
	\end{center}
\end{figure}
It is straightforward to convert webs in the lower half-plane to webs in a disk, see Remark~\ref{morphisms_remark} and Figure~\ref{boundeqiv}. 
 As part of the pivotal structure of the category, non-elliptic webs $W_1$, $W_2$  in Figure \ref{basisspan} and non-elliptic webs $\Gamma_1$, $\Gamma_2$  in Figure \ref{4_00} give rise to the same webs in a disk.
\end{example}

Morphisms from $\Gamma_i$ to $\Gamma_j$ in $\mcH^{\unde}$ are $\Z$-linear combinations of foams $F$ with boundary given by the union of $\Gamma_i$, $\Gamma_j$ and vertical lines, modulo the relations on foams given by the universal construction in Section~\ref{subsec_Kuperberg}.
\begin{figure}[H]
	\begin{center}
		{\includegraphics[width=150pt]{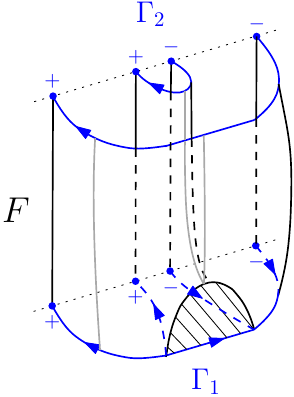}}
		\caption{Foam $F\in\textrm{Hom}_{\mcH^{\unde}}(\Gamma_1,\Gamma_2)$, where $\Gamma_1$, $\Gamma_2$ are highlighted in blue, singular saddle is highlighted in gray. This foam is homeomorphic rel boundary to the foam considered in Figure \ref{3_13}.}
		\label{4_01}
	\end{center}
\end{figure}
An example of such foam is shown in Figure~\ref{4_01}, with 
the boundary $\partial F$ depicted in Figure~\ref{4_02} on the left. In general, for a foam $F$ representing a morphism from $\Gamma_1$ to $\Gamma_2$,  
$$
\partial F = (-\Gamma_1)\times \{0\} \cup \Gamma_2\times \{1\} \cup \unde \times [0,1]. 
$$
An equivalent definition of the graded abelian group which is the hom space $\textrm{Hom}_{\mcH^{\unde}}(\Gamma_1,\Gamma_2)$
is  
\begin{equation}\label{eq_H_hom}
\Hom_{\mcH^{\unde}} (\Gamma_i,\Gamma_j) \ := \ \brak{\overline{\Gamma}_i\Gamma_j}\{|\unde|\}.  
\end{equation}
Here a grading shift is added, where $|\unde|$ is the length of the sequence $\unde$. 
The correspondence between the two definitions, see also Figure~\ref{4_02}, is given by contracting vertical intervals of the boundary of $F$ to points and then fanning out the wedge of the space between the half-planes of $\Gamma_1$ and $\Gamma_2$ into the upper half-space $\R^3_+$.   
\begin{example}
An example of this identification for $\unde=(+,+,-,-)$ is shown in Figure~\ref{4_02}, where on the left we depict a boundary of a foam $F\in \Hom_{\mcH^{\unde}} (\Gamma_1,\Gamma_2)$ and on the right we flatten it into the planar web $\overline{\Gamma}_1\Gamma_2$. The space of morphisms is given by the state space of this web. 
\begin{figure}[H]
	\begin{center}
		{\includegraphics[width=230pt]{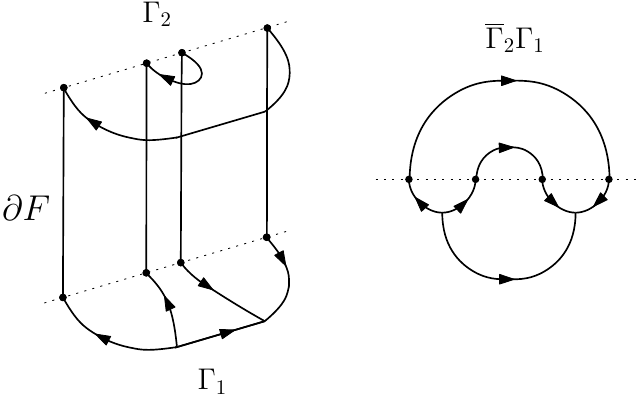}}
		\caption{Identification of $\partial F$ with $\overline{\Gamma}_2\Gamma_1$.}
		\label{4_02}
	\end{center}
\end{figure}
For this sequence the category has only two objects $\Gamma_1$, $\Gamma_2$, see Figure~\ref{4_00}. There are four webs describing spaces of morphisms in $\mcH^{\unde}$. For each of them, we assign corresponding state spaces of graphs, as depicted in Figure \ref{4_03}.
\begin{figure}%[H]
	\begin{center}
		{\includegraphics[width=220pt]{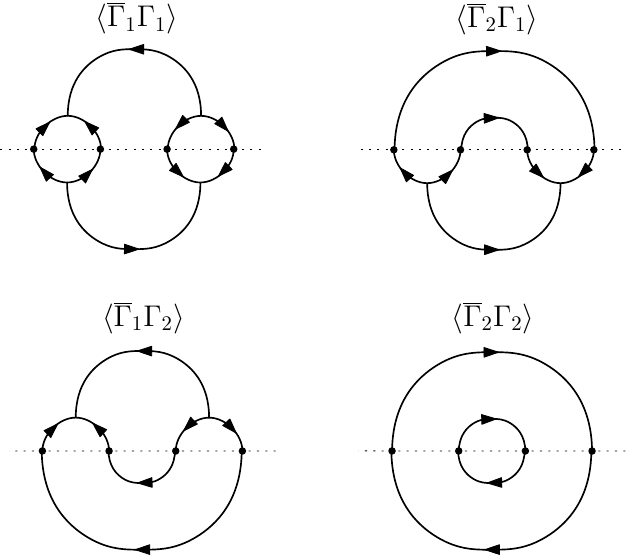}}
		\caption{State spaces assigned to pairs of generating objects $\Gamma_1,\Gamma_2$ for $\unde=(+,+,-,-)$, all with the shift $\{4\}$ from \eqref{eq_H_hom}.}
		\label{4_03}
	\end{center}
\end{figure}
\end{example}

Composition of morphisms in the category $\mcH^{\unde}$ is given by stacking foams on top of each other and extending by linearity. We see that the composition is easy to define using presentation of the boundary on the left of Figure~\ref{4_01}. 

\begin{figure}[H]
	\begin{center}
		{\includegraphics[width=400pt]{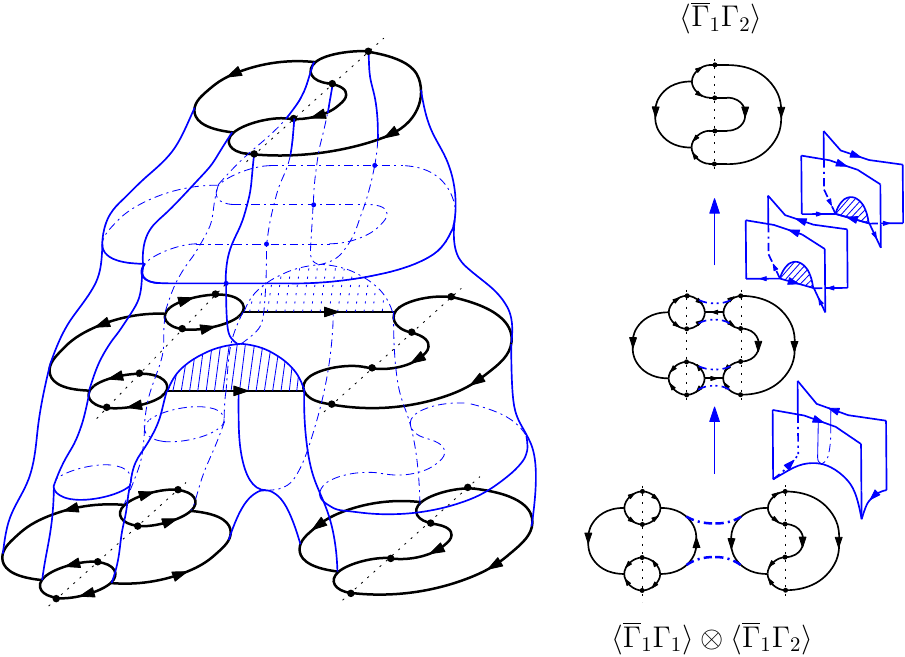}}
		\caption{Multiplication of elements $\langle \overline{\Gamma}_1\Gamma_1\rangle$ and $\langle \overline{\Gamma}_1\Gamma_2\rangle$.}
		\label{4_04}
	\end{center}
\end{figure}

If one instead uses its presentation as $\overline{\Gamma}_1\Gamma_2$, on the right of Figure~\ref{4_01}, and definition \eqref{eq_H_hom}, composition should be defined as follows. The composition of morphisms from $\Gamma_i$ to $\Gamma_j$ with morphisms from $\Gamma_j$ to $\Gamma_k$ is induced by a cobordism from $\overline{\Gamma}_k\Gamma_j \overline{\Gamma}_{j}\Gamma_i$ to $\overline{\Gamma}_k\Gamma_i$, which couples $\Gamma_j$ to its reflection $\overline{\Gamma}_{j}$ converting $\Gamma_j\overline{\Gamma}_{j}$ to the identity cobordism. This cobordism induces a linear map of state spaces 
\begin{equation}
	\langle \overline{\Gamma}_k\Gamma_j\rangle \otimes \langle \overline{\Gamma}_{j}\Gamma_i\rangle\  \mapsto \ \langle \overline{\Gamma}_k\Gamma_i\rangle 
\end{equation}
and gives the desired composition of morphisms. 
The map is homogeneous of degree $|\unde|$, so with the shift as in \eqref{eq_H_hom} the multiplication is degree-preserving. 
An example of the composition of morphisms is shown in Figure~\ref{4_04}.

\vspace{0.07in}

Category $\mcH^{\unde}$ has finitely many objects, and it is a pre-additive category. The sum of its hom spaces, over all possible pairs of objects, is an associative ring. 
Namely, define the graded ring $\mmH^{\unde}$ by 
\begin{equation}
    \mmH^{\unde} \ := \ \bigoplus\limits_{\Gamma_i,\Gamma_j\in \B^{\unde}}\Hom_{\mcH^{\unde}}(\Gamma_i,\Gamma_j) = 
    \bigoplus\limits_{\Gamma_i,\Gamma_j\in \B^{\unde}} \brak{\overline{\Gamma}_i\Gamma_j}\{|\unde|\}.
\end{equation}
That is, we form the sum of spaces of morphisms over all pairs of objects.

Multiplication
\begin{eqnarray*}
	\mmH^{\unde} \otimes \mmH^{\unde}& \to& \mmH^{\unde}\\
	\langle \overline{\Gamma}_k\Gamma_j\rangle \otimes \langle \overline{\Gamma}_{j^\prime}\Gamma_i\rangle& \mapsto& \delta_{j,j^\prime}\langle \overline{\Gamma}_k\Gamma_i\rangle
\end{eqnarray*}
works similarly to the cobordism discussed in Figure \ref{3_07}. An example of multiplication of elements $\langle \overline{\Gamma}_1\Gamma_1\rangle$ and $\langle \overline{\Gamma}_1\Gamma_2\rangle$ is given in Figure \ref{4_04}.

For any $\Gamma_i$ there is an idempotent $1_{\Gamma_i}\in \langle \overline{\Gamma}_i\Gamma_i\rangle$, such that
$$
1_{\Gamma_i} x = x, \ x 1_{\Gamma_j}=x, \ \ \forall x\in   \langle \overline{\Gamma}_i\Gamma_j\rangle.
$$
This idempotent can be described by a canonical cobordism in $\mathrm{Hom}_{\mcH^{\unde}}(\emptyset,\overline{\Gamma}_i\Gamma_i)$. Imagine $\Gamma_i$ to be a side of a fan and then fully expand the fan so that its opposite sides are on the same plane. The two opposite sides form the web $\overline{\Gamma}_i\Gamma_i$, and the fully opened fan is a cobordism from $\emptyset$ to that web representing $1_{\Gamma_i}$. 
\begin{example}
Recall that for $\unde=(+,+,-,-)$ the objects of $\mcH^{\unde}$ are given in Figure \ref{4_00} and the state spaces are those of closed webs in Figure \ref{4_03}. An idempotent for $1_{\Gamma_1}$ is described by a canonical cobordism from the empty web to the web $\overline{\Gamma}_1\Gamma_1$ depicted in Figure \ref{4_05}.
\begin{figure}[H]
	\begin{center}
		{\includegraphics[width=250pt]{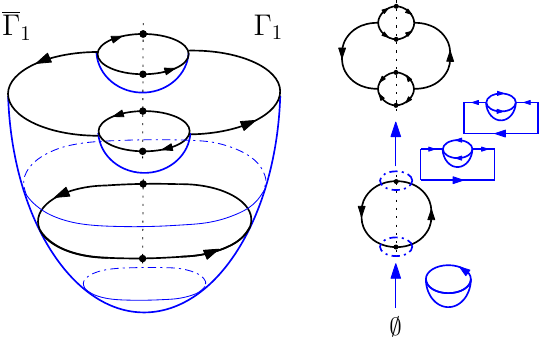}}
		\caption{Example of the idempotent $1_{\Gamma_1}$ for $\Gamma_1\in \B^{\unde}$.}
		\label{4_05}
	\end{center}
\end{figure}
\end{example}

Multiplication by such an idempotent is a projection. Summing over all these mutually-orthogonal idempotents gives the unit element of $\mmH^{\unde}$:
$$
1_{\mmH^{\unde}} \ = \ \sum_{\Gamma_i\in \B^{\unde}} 1_{\Gamma_i}. 
$$
This allows us to decompose $\mmH^{\unde}$ into a matrix-like graded ring
\begin{equation}\label{eq_H}
\mmH^{\unde}=\bigoplus_{\Gamma_i,\Gamma_j\in \B^{\unde}} 1_{\Gamma_j} \mmH^{\unde} 1_{\Gamma_i}=\bigoplus_{\Gamma_i,\Gamma_j\in \B^{\unde}} \langle \overline{\Gamma}_j\Gamma_i\rangle \{|\unde|\}.
\end{equation}

\vspace{0.07in}

For sequences $\unde,\undep\in \Seq$ denote by $\Webst_{\unde}^{\undep}$ the set of webs $U\in \R\times [0,1]$ with the top boundary $\undep$ and bottom boundary $\unde$: 
\[\partial U = \undep \times \{1\} \cup (-\unde)\times \{0\}.
\]
Elements in this set are considered up to a rel boundary isotopy. The concatenation $U_2 \circ U_1$ for composable webs $U_2,U_1$ is associative. Define a monoidal category $\Webst$ with sequences $\unde$ as objects and the set of homs from $\unde$ to $\undep$ given by $\Webst_{\unde}^{\undep}$:
\begin{equation}
    \Hom_{\Webst}(\unde,\undep) \ := \ \Webst_{\unde}^{\undep}.
\end{equation}
Composition of morphisms is given by concatenating webs with boundary. Category $\Web$, introduced in Section~\ref{websbnd}, can be viewed as a linearization of $\Webst$. 
\begin{prop}
    There is a natural functor $\Webst\lra \Web$ which is the identity on objects and takes each web to itself, viewed as an element of the Kuperberg skein module of webs with a given boundary. 
\end{prop}
In particular, there is a natural  map $\Webst_{\unde}^{\undep}\lra \Web_{\unde}^{\undep}$ taking a web to the corresponding element of the $\Z[q,q^{-1}]$-module generated by all webs with that boundary modulo the Kuperberg relations.
Note that this map is not injective: pick a web $\Gamma$ with at least two regions and place a circle in one of the two regions producing webs $\Gamma',\Gamma''$. The two webs define different morphisms in $\Webst$, since they are not rel boundary isotopic, but the same morphism in $\Web$, since a circle reduces to the constant $[3]$ in the ground ring. 

\vspace{0.07in}
     
Now, pick two balanced sequences $\unde$, $\undep\in \Bal$ and a web $U\in \Webst_{\unde}^{\undep}$. 
Define a graded $(\mmH^{\undep},\mmH^{\unde})$-bimodule
\begin{equation}\label{eq_T_to_closed}
\mmH(U):=\bigoplus_{\Gamma_j\in \B^{\undep},\Gamma_i\in \B^{\unde}} \langle \overline{\Gamma}_j U \Gamma_i\rangle \{|\unde|\}.
\end{equation}
\begin{figure}[H]
	\begin{center}
		{\includegraphics[width=400pt]{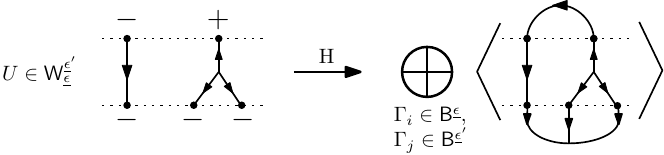}}
		\caption{An example of $\mmH(U)$ for $\unde=(-,-,-)$ and $\undep=(-,+)$. For these sequences diagrams $\Gamma_i,\Gamma_j$ are unique, and $\mmH(U)$ is the state space of the single web shown on the right and called \emph{the $\Theta$-web}, with the grading shifted up by $3$.}
	\end{center}
\end{figure}
In other words, we complete $U$ to a closed web by attaching a non-elliptic web $\Gamma_i$ at its lower endpoints, attaching a non-elliptic web $\overline{\Gamma}_j$ at its upper endpoints, forming the state space of the resulting closed web $\overline{\Gamma}_j U \Gamma_i$ and summing over all such $\Gamma_i,\Gamma_j$. At the end, apply the overall grading shift $\{|\unde|\}$ by the length of the lower sequence $\unde$. 

To see the right $\mmH^{\undep}$-module structure, write 
\begin{equation}\label{eq_H_2}
\mmH^{\undep}=\bigoplus_{\Gamma_i,\Gamma_k\in \B^{\undep}} \langle \overline{\Gamma}_i\Gamma_k\rangle \{|\undep|\}.
\end{equation}
The cobordism from $\overline{\Gamma}_j U \Gamma_i \sqcup \overline{\Gamma}_i\Gamma_k$ to $\overline{\Gamma}_j U \Gamma_k$ that couples $\Gamma_i$ to its reflection $\overline{\Gamma}_i$ induces a map of state spaces 
\[
\brak{\overline{\Gamma}_j U \Gamma_i} \otimes \brak{\overline{\Gamma}_i\Gamma_k} \lra \brak{\overline{\Gamma}_j U \Gamma_k}.
\]
Summing over all $i,j,k$ as above results in a right action of $\mmH^{\undep}$ on $\mmH(U)$. A left action of $\mmH^{\unde}$ on $\mmH(U)$ is constructed likewise, and the two actions commute since the coupling cobordisms for the two actions are nontrivial in two disjoint regions of $\R^2\times [0,1]$. Consequently, $\mmH(T)$ is a bimodule over these two rings. 

Consider composable webs $U_1\in\Webst^{\undep}_{\unde}$ and $U_2\in\Web^{\unde^{\prime\prime}}_{\undep}$. One can check that there is a natural isomorphism of bimodules 
$$
\mmH(U_2\circ U_1)\cong \mmH(U_2)\otimes_{\mmH^{\undep}} \mmH(U_1)
$$
by reducing to the case when sequences $\unde,\unde^{\prime\prime}$ are empty, see~\cite{MPT} and also~\cite{Kh2} for the analogous result in the simpler, $\mathsf{SL}(2)$, case. 

Next, pick webs $U_0,U_1\in\Webst^{\undep}_{\unde}$. 
Consider a foam cobordism $S$ from $U_0$ to $U_1$. Such a cobordism is an oriented $\mathsf{SL}(3)$ foam $S$ with boundary and corners in $\R\times [0,1]^2$ such that $S\cap (\R\times[0,1]\times \{i\})=U_i$ for $i=0,1$. On the other two boundary sides of $\R\times [0,1]^2$ foam $S$ has the simplest possible structure and is given by the direct products of sets of points on $\R$ representing sequences $\unde,\undep$ times $[0,1]$.
\begin{figure}[H]
	\begin{center}
		{\includegraphics[width=350pt]{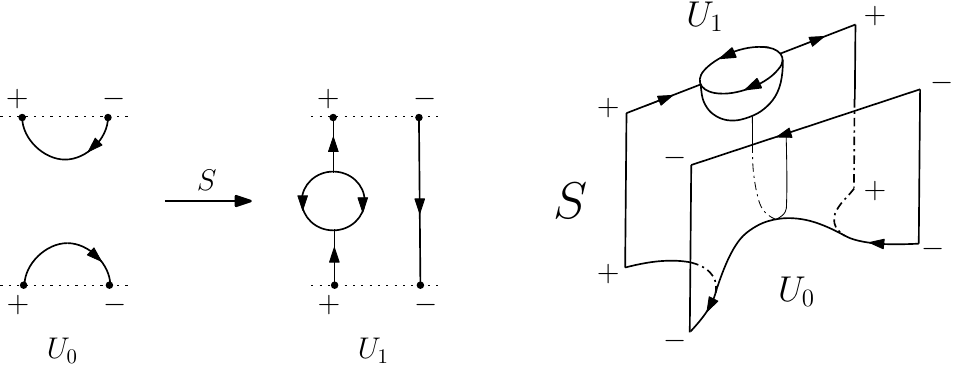}}
		\caption{Example of a cobordism $S$ from $U_0$ to $U_1$, where $U_0,U_1\in\Webst^{\undep}_{\unde}$ for $\unde=\undep=(+,-)$.}
	\end{center}
\end{figure}
To $S$ we assign a homomorphism $\mmH(S)$ of $(\mmH^{\undep},\mmH^{\unde})$-bimodules. To define it, recall that a web with boundary $U$ is converted into closed webs $\overline{\Gamma}_j U \Gamma_i$ to get the corresponding bimodule, see (\ref{eq_T_to_closed}). A similar operation needs to be done on $S$. Foam $S$ can be turned into a foam cobordism between closed webs by capping it offf with   $\Gamma_i\times[0,1]$ on one side and with $\overline{\Gamma}_j\times[0,1]$  on the other. The resulting foam 
\[
(\overline{\Gamma}_j\times[0,1])\,S \, (\Gamma_i\times[0,1])
\]
induces a homomorphism of state spaces from its boundary web $\overline{\Gamma}_j U_0\Gamma_i$  to the other boundary web, $\overline{\Gamma}_j U_1\Gamma_i$. 
Summing these homomorphisms over all $\Gamma_i\in \B^{\unde}$ and $\Gamma_j\in \B^{\undep}$ results in a linear map 
\begin{equation}
\mmH(S)\ : \ \mmH(U_0)\to \mmH(U_1)
\end{equation}
which is a homomorphism of bimodules (and it is homogeneous of certain degree).

%defined by closing up $S$ in all possible ways with $\Gamma_i\times[0,1]$, $\overline{\Gamma}_j\times[0,1]$ for $\Gamma_i\in \B^{\unde},\Gamma_j\in \B^{\undep}$ and forming the direct sum 

%$$\mmH(S):=\bigoplus_{\Gamma_i\in \B^{\unde},\Gamma_j\in \B^{\undep}} \mmH((\overline{\Gamma}_j\times[0,1])\,S \, (\Gamma_i\times[0,1])), $$
%where
%$$\mmH((\overline{\Gamma}_j\times[0,1])\,S \,( \Gamma_i\times[0,1])):\mmH(T_1)\to \mmH(T_2).$$

\vspace{0.07in}

Our invariant assigns the bimodule $\mmH(U)$ to a web with boundary $U\in \Webst_{\unde}^{\undep}$. We can now extend it from webs to tangles. Given an oriented tangle $T\in \Tan_{\unde}^{\undep}$ pick its generic projection $D$ onto the strip $\R\times [0,1]$. Suppose first that $D$ has a single crossing $c$. Mimicking the rules in Figure~\ref{cones_fig}, consider webs $D_0,D_1$ given by resolving the crossing as shown in that figure. There is a canonical cobordism $S_c$ from $D_0$ to $D_1$ which defines a homomorphism of bimodules $\mmH(D_0)\stackrel{\mmH(S_c)}{\lra}\mmH(D_1)$, see Figure~\ref{3_13}. Now, to a single-crossing tangle projection $D$ associate the complex of bimodules
\[
\mmH(D) \ : \ 0 \lra \mmH(D_0) \stackrel{\mmH(S_c)}{\lra}\mmH(D_1) \lra 0,
\]
with the homological degree and $q$-grading shifts as shown in Figure~\ref{cones_fig}. Given an arbitrary generic projection $D$ of a tangle $T$, factorize it into the composition of projections with at most one crossing each, $D=D^m\circ \dots D^2\circ D^1$, and assign to $D$ the tensor product of complexes of bimodules for the projections with at most one crossing:
\[
\mmH(D)  \ := \ \mmH(D^m)\otimes \dots \otimes  \mmH(D^1).
\]
The tensor product is taken over the intermediate rings $\mmH^{\unde_{m-1}}, \dots, \mmH^{\unde_{1}}$, where the diagram $D^i$ has the bottom sign sequence $\unde_{i-1}$ and top sign sequence $\unde_i$.
If $D^i$ has no crossings, it can also be viewed as a web without vertices and has the associated bimodule $\mmH(D^i)$, which corresponds to a complex with the unique nontrivial term in homological degree $0$. 

One can check that the Reidemeister moves of tangle diagrams induce homotopy equivalences between complexes $\mmH(D)$ of graded bimodules, see~\cite[Chapter 2]{Ro}, for instance. 

\begin{theorem}
    Homotopy equivalence class of $\mmH(D)$ in the category of complexes of graded $(\mmH^{\undep},\mmH^{\unde})$-bimodules is an invariant of an oriented tangle $T$, denoted $\mmH(T)$.
\end{theorem}

Checking higher-level moves (movie moves for tangles~\cite{CS}), one can establish that $\mmH(T)$ and not just its equivalence class is an invariant of $T$. Furthermore, given a tangle cobordism $S$ from tangle $T_0$ to tangle $T_1$, there is a well-defined homomorphism 
\[
\mmH(S) \ : \ \mmH(T_0) \lra \mmH(T_1) 
\]
of complexes of bimodules, in the category of complexes modulo chain homotopies. 

Together, these homomorphisms define a 2-functor from the 2-category of tangle cobordisms to the 2-category of complexes of bimodules up to chain homotopies. Let us provide details. 

Define $\Tcob$ to be the following 2-category:  
\begin{itemize}
    \item Objects: balanced sign sequences $\unde$, 
    \item 1-morphisms from $\unde$ to $\undep$: tangles $T$ with the lower, respectively upper boundary given by the sign sequence $\unde$, respectively $\undep$, 
    \item 2-morphisms from $T_0$ to $T_1$ (where the two tangles have matching lower and upper boundaries): oriented surfaces $S\in \R^2\times [0,1]^2$ up to rel boundary isotopies with $\partial_0 S = -T_0$, $\partial_1 S = T_1$. The two side boundaries of $S$ are the standard products $\unde\times [0,1]$, $\undep\times [0,1]$ in $\R^2\times [0,1]$.  
\end{itemize}
(Since one restricts to balanced sign sequences, we may alternatively call $\Tcob$ the 2-category of \emph{balanced} tangle cobordisms.) 

The 2-category $\mmH$ has 
\begin{itemize}
    \item Objects: sign sequences $\unde$, 
    \item 1-morphisms from $\unde$ to $\undep$: complexes of graded $(\mmH^{\undep},\mmH^{\unde})$-bimodules with degree-preserving differential, 
    \item 2-morphisms are degree-preserving homomorphisms of complexes of graded bimodules up to chain homotopies. 
\end{itemize}

\begin{theorem}\label{thm-clark}
    The above construction is a 2-functor 
\begin{equation}
    \mmH \ : \ \Tcob \lra \mmH
\end{equation} 
from the 2-category of tangle cobordisms $\Tcob$ to the 2-category $\mmH$ of complexes of bimodules up to chain homotopies. 
\end{theorem}

This theorem follows at once from Clark~\cite{Cl}.  
Clark works in the Bar-Natan framework of \emph{canopolies}~\cite{BN}, where tangles  live in a 3-ball rather than in $\R^2\times [0,1]$, having a single boundary sign sequence rather than a pair of sequences (and likewise for webs, foams and cobordisms). Canopolies framework also allows one to work with all tangles and tangle cobordisms rather than restricting to balanced ones. 

In the easier case of categorified $\mathsf{SL}(2)$ invariants, the analogue of the above theorem is established in~\cites{Kh5,BN}. In those papers the invariant $\mmH(S)$, for a tangle cobordism $S$, is shown to be well-defined up to an overall sign. There are several ways to fix the sign problem~\cites{CMW,Bl,Ca,Vog,Sa}. The sign issue does not appear in the $\SL(N)$ link homology functoriality for $N\ge 3$, see~\cite{ETW}, and can be traced, on the decategorified level, to the skew-commutativity of the self-duality isomorphism for the fundamental representation of Lie group $\SL(2)$, see~\cite{CMW}. 

\vspace{0.07in}

To summarize, a sequence $\unde$ gives rise to the category $\mcH^{\unde}$ or, equivalently, to the graded associative unital ring $\mmH^{\unde}$. This ring is also Frobenius, equipped with a nondegenerate trace. Rings $\mcH^{\unde}$ and their bimodules associated to webs with boundary give rise to complexes of bimodules for tangle diagram and to the 2-functor $\mmH$ in Theorem~\ref{thm-clark}. This 2-functor specializes to the bigraded link homology theory in~\cite{Kh4}, functorial for link cobordisms.  

Rings $\mmH^{\unde}$ are studied in depth in~\cite{Ro} and~\cite{MPT}. Every web with boundary $\unde$ gives rise to a projective $\mmH^{\unde}$-module. 
Non-elliptic webs with few endpoints give rise to indecomposable projective $\mmH^{\unde}$-modules, but eventually projective modules for non-elliptic web stop being indecomposable~\cite{MN}, which is a categorified version of the observation that the basis of non-elliptic webs eventually diverges from the Lusztig dual canonical basis in the space of $U_q(\slthree)$ invariants of tensor products of fundamental representations~\cite{KhK}.  Classification of non-elliptic webs giving rise to indecomposable projective modules was achieved in~\cite{Ro}. 

Analogues of rings $\mmH^{\unde}$ for unoriented $\SL(3)$ webs and foams are not understood, in general. Variations on these rings for up to 8 boundary points are studied in~\cite{Th}.

$\SL(2)$ analogues of web rings, defined earlier in~\cite{Kh3}, are easier to understand. A short introduction to these rings and (likely incomplete) list of their applications in low-dimensional topology and representation theory, with many references, can be found in~\cite{KhL}, see also the recent work~\cite{BVDHMS}. 

A computer program to compute $\SL(3)$ homology of knots and links has been developed by Lewark~\cites{Le1,Le2}. Mackaay and Vaz~\cite{MaVa} proved that $\SL(3)$ link homology constructed via foams and matrix factorizations are naturally isomorphic. 

Khovanov (or $\SL(2)$) link homology~\cite{Zh} gives rise to the Rasmussen invariant~\cite{Ra}, which is a homomorphism from the knot concordance group onto $2\Z$ with important properties. Using his invariant, Rasmussen gave the first algebraic proof of the Kronheiner-Mrowka-Milnor theorem on the slice genus of positive knots~\cite{Ra}, without the use of gauge theory. Similar concordance invariants exist for the $\SL(3)$ and $\SL(N)$ link homology theories, and we refer the reader to~\cites{Lo,Le3,LL,Sch} for their definition, properties, and independence from the original Rasmussen invariant. 

See Manolescu~\cite{Man} for 
a very recent review of applications of combinatorial link homology to $4$-manifold topology, including to smooth structures on 4-manifolds.  

Mackaay, Pan and Tubbenhauer~\cite{MPT} show that the center of $\mmH^{\unde}$ is isomorphic to the cohomology ring of a suitable Spaltenstein variety. One expects that this isomorphism can be upgraded to an equivalence of derived categories of graded projective $\mmH^{\unde}$-modules and a suitable subcategory of the category of equivariant coherent sheaves on the quiver variety thickening of that Spaltenstein variety. The $\SL(2)$ case of that correspondence can be found in~\cite{AN}, see also~\cites{CK1,CK2}.

Ring $\mmH^{\unde}$, their Morita and derived Morita equivalent rings, and their generalizations from $N=2,3$ to any $\SL(N)$ appear throughout representation theory and categorification in many ways. We point the reader to the papers~\cites{LQR,MS,Mc1,McY,QR,McW} to explore these connections and for more information on structural and representation-theoretical meaning of webs and foams. Annular $\SL(N)$ foams are developed in~\cites{QRS,GoWe,QR2} and \cites{AkKh,Ak} for non-equivariant, respectively equivariant, annular link homology. See also \cite{AkKhZ} for $\SL(3)$ annular web rings in the equivariant setting. 
$\SL(3)$ homology is also reviewed in online notes~\cite{KhQ}.

\subsection{Frobenius rings from \texorpdfstring{$n$}{n}-dimensional topological quantum field theories}\hfill\vspace{0.02in}\\
Categories $\mcH^{\unde}$ and rings $\mmH^{\unde} $ defined in Section~\ref{subsec_sl3web} exist in far greater generality, as we now explain, although the usefulness of the general setup is unclear. 
 Suppose we have a $n$-dimensional TQFT (or a lax TQFT), that is, a (lax) tensor functor
$$
\mathcal{Z}:\textbf{Cob}_n\to \textbf{Vect}_\Bbbk 
$$
from some category of cobordisms to the category of vector spaces. In the example below, we consider the category of $n$-dimensional oriented cobordisms between closed oriented $(n-1)$-manifolds. 

The functor, in particular, takes a closed $n$-dimensional manifold $M^n$  to a number $\mathcal{Z}(M^n)$:
\begin{equation*}
    M^n,\,\partial M^n=\emptyset \mapsto  \mathcal{Z}(M^n)\in\Bbbk.
\end{equation*}
%\begin{eqnarray*}
%    \mathcal{Z}:  \textrm{Hom}_{\textbf{Cob}_n} & \to & \textrm{Hom}_{\textbf{Vect}_\Bbbk}\\
%     M,\,\partial M=\emptyset & \mapsto & \mathcal{Z}(M)\in\Bbbk
%\end{eqnarray*}
It also takes a closed $(n-1)$-dimensional manifold $N^{n-1}$ and assigns to it the finite-dimensional vector space $\mathcal{Z}(N^{n-1})$.
\begin{eqnarray*}
    \mathcal{Z}:  Ob(\textbf{Cob}_n) & \to & Ob(\textbf{Vect}_\Bbbk)\\
     N & \mapsto & \mathcal{Z}(N)
\end{eqnarray*}
Given a closed $(n-2)$-dimensional manifold $K^{n-2}$, one can build a category $\mathcal{Z}(K^{n-2})$ as follows. Objects $Ob(\mathcal{Z}(K))$ are $(n-1)$-manifolds $N^{n-1}$ such that $\partial N=K$. To describe morphisms in $\mathcal{Z}(K)$ from $N_1$ to $N_2$, form the $(n-1)$-manifold
$$
\overline{N}_2 N_1=N_1\underset{K}{\cup} \overline{N}_2
$$
given by the two manifolds $N_1$, $N_2$ glued along the boundary $\partial N_1=\partial N_2=K$, one with the opposite orientation. Set 
\[
\textrm{Hom}_{\mathcal{Z}(K)}(N_1,N_2) \ := \mathcal{Z}(\overline{N}_2 N_1).
\]
Composition of morphisms is given by the coupling cobordism between the composition $N_2\overline{N}_2$ and the product $K\times [0,1]$, similar to that in the Figure~\ref{4_04}. 

%will be the set of all possible cobordisms $\textrm{Hom}(N_1,N_2)\in \textrm{Hom}_{\textbf{Cob}_n}$, where $N_1,N_2\in Ob(\mathcal{Z}(K))$. 

\begin{figure}[H]
	\begin{center}
		{\includegraphics[width=330pt]{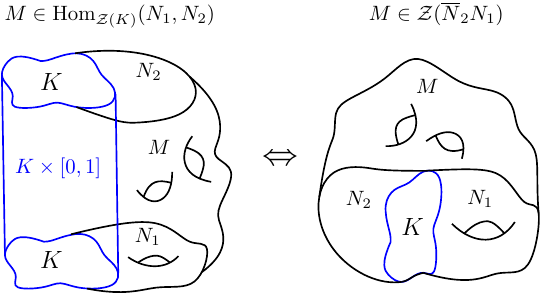}}
		\caption{Left-hand side represents manifold $M$ as a cobordism in $\textrm{Hom}_{\mathcal{Z}(K)}(N_1,N_2)$, right-hand side represents manifold $M$ as an element of the $\Bbbk$-vector space $\mathcal{Z}(\overline{N}_2 N_1)$.}
		\label{4_06}
	\end{center}
\end{figure}

An $n$-manifold $M$ with boundary and corners, which is a cobordism from $N_1$ to $N_2$, as shown in Figure~\ref{4_06}, 
gives an element of $\textrm{Hom}_{\mathcal{Z}(K)}(N_1,N_2)$, by converting it to (essentially the same) manifold with the boundary $\overline{N}_2 N_1$. More precisely, $M^n$ satisfies 
$$
\partial M\cong N_1\cup(\overline{N}_2)\cup (K\times [0,1]).
$$
% One can apply the functor $\mathcal{Z}$ to obtain a $\Bbbk$-vector space $\mathcal{Z}(\overline{N}_2 N_1)$, and from Figure \ref{4_06} we see that $\textrm{Hom}_{\mathcal{Z}(K)}(N_1,N_2)=\mathcal{Z}(\overline{N}_2 N_1)$.
In the category $\mathcal{Z}(K)$, it seems natural to further restrict to objects $N$ with $\partial N\cong K$ such that each connected component of $N$ has nonempty boundary, and pick one representative from each isomorphism class (rel boundary) of such $N$. A variation of this construction, using foams embedded in $\R^3$ as cobordisms, recovers the category $\mcH^{\unde}$ from Section~\ref{subsec_sl3web}. Summing over all the hom spaces results in an associated ring, similar to our conversion from $\mcH^{\unde}$ to $\mmH^{\unde}$. If $\mathcal{Z}(K)$ has infinitely many objects $N$, even with the above restrictions, the resulting ring is non-unital but it is equipped with a system of idempotents, one for each $N$, which serves as a replacement for the unit element.  

For us, categories $\mcH^{\unde}$ and rings $\mmH^{\unde}$ were  an intermediate step in the extension of link homology to tangles. It would be interesting to find other uses of categories $\mathcal{Z}(K)$ for suitable cobordism categories and tensor functors $\mathcal{Z}$.  

%Manifolds $M^n$ such that
%$$\partial M=N_1\cup(\overline{N}_2)\cup (K\times [0,1])$$
%give a map
%$$N_2\overset{[M]}{\to} N_1,$$
%which defines an element of $\mathcal{Z}(\overline{N}_2 N_1)$. This shows that
%$$
%\mathcal{Z}(K)=\bigoplus_{N_i,N_j\in Ob(\mathcal{Z}(K))}\mathcal{Z}(\overline{N}_i N_j)
%$$
%is a category with an algebra structure.
\section{Oriented \texorpdfstring{$\GL(N)$}{GL(N)} foams and a categorification of Reshetikhin-Turaev link invariants}\label{sec_GLN}

The foam evaluation described in Section~\ref{sec-one} followed~\cite{KhR}. It is a specialization of the unoriented version of the Robert-Wagner foam evaluation formula \cite{RW1} from $N$ colors to three. Motivation for this formula came from the work of Kronheimer and Mrowka \cites{KM15, KM16, KM17}. They defined a more general homology theory for trivalent graphs embedded into oriented $3$-manifolds, which comes from the $\SO(3)$ gauge theory for $3$-orbifolds. 
Kronheimer and Mrowka also conjectured the existence of a combinatorial topological theory for planar trivalent graphs and foam cobordisms between them, and evaluation formula in~\cite{KhR} allowed to prove their conjecture. 

In the first appearance of foams in relation to link homology, only oriented bipartite planar graphs and oriented $\SL(3)$ foam cobordisms between them were needed~\cite{Kh1}. In this case foams have no vertices and their evaluation can be written down in a simpler way, as explained in Section~\ref{subsec_Kuperberg}. This leads to a link homology theory which categorifies the Kuperberg bracket~\cites{Kh1,MV,MN} while avoiding complexities related to having foams with singular vertices. 

Once vertices are present one cannot localize and reduce the evaluation to that near singular seams as described in Section~\ref{subsec_Kuperberg}. 
In oriented $\SL(N)$ foams, where $N\ge 4$, vertices are unavoidable (unlike the oriented $\SL(3)$ foams). A combinatorial approach to $\SL(N)$ foams for $N\ge 4$ via the universal construction was absent until relatively recently~\cite{RW1}, and more categorical approaches played the main role in constructing $\SL(N)$ link homology theories~\cites{KhR,S,MS,W2}. Some of these earlier and categorical approaches utilize categorifications of quantum groups and their representations~\cites{ChRo,La1,KhLa1,Rou,La2,Br,W2}, combined with the categorified Howe duality~\cites{LQR,QR} and result in the theory of $\SL(N)$ foams and the associated link homology theories. One also expects that the Kapustin-Li formula for foams~\cites{MSV,KhRz2} can be linked to the Robert-Wagner evaluation. 

\vspace{0.07in}

\begin{figure}[H]
	\begin{center}
		{\includegraphics[width=280pt]{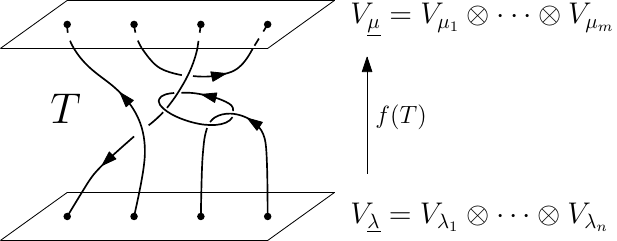}}
		\caption{Reshetikhin-Turaev invariant $f(T)$.}
		\label{5_00}
	\end{center}
\end{figure}

Starting at the decategorified level, let us first recall the setup of quantum link and tangle invariants, via representations of quantum groups. 
(Tangles and tangle cobordisms may need to be framed, but we do not explicitly mention this below.) 
 One considers an oriented tangle $T$, where components of the tangle are colored by irreducible representations $V_{\lambda}$ of the quantum group $U_q(\mathfrak{g})$. The invariant of $T$ is a map between tensor products of these representations, as shown in Figure \ref{5_00}. 
This is the Reshetikhin-Turaev construction \cite{RT}, which is a (braided monoidal) functor from the category $\mathsf{OFTan}_{\mfg}$ of oriented framed tangles, with components labelled by irreducible representations $V_{\lambda}$ of $\mfg$, to the category $\Rep(U_q(\mathfrak{g}))$,
\begin{equation}\label{eq_f_V}
    f \ : \ \mathsf{OFTan}_{\mfg} \lra \Rep(U_q(\mathfrak{g})). 
\end{equation}
When the tangle is a link $L$, the Reshetikhin-Turaev invariant is a homomorphism from the trivial representation to itself, given by a scalar $P(L)$, which lives in $\Z[q,q^{-1}]$ (more precisely, the invariant of links lives in $\Z[q^{\pm 1/D}]$, where $D$ depends on $\mathfrak{g}$, see~\cite{TLe}, and we then relabel $q^{1/D}$ to $q$, for simplicity).

In the extension of the Reshetikhin-Turaev invariants from links to tangles, tensor products $V_{\underline{\lambda}}$ may be replaced by their invariant subspaces 
\[
\Inv(V_{\undlambda}) \ :=  \ \Hom_{U_q(\mfg)}(V_0,V_{\undlambda}),
\]
where $V_0$ is the trivial representation, with the map $f(T)$ restricting to the corresponding map, also denoted $f(T)$, between the invariant spaces:
\[
f(T) \ : \ \Inv(V_{\undlambda}) \lra \Inv(V_{\underline{\mu}}).
\]
This gives a variation on the functor $f$ in \eqref{eq_f_V}, which in this case takes values in the category of vector spaces over the ground field $\Q(q)$: 
\begin{equation}\label{eq_f_V2}
    f \ : \ \mathsf{OFTan}_{\mfg} \lra \Q(q)\mathsf{-vect}. 
\end{equation}

\vspace{0.07in}

Let us describe a perfect model for a categorification of Reshetikhin-Turaev invariants $P(L)$, for oriented links $L$, and, more generally, $f(T)$ for tangles $T$. 

A categorification of this invariant requires replacing a tensor product $V_{\undlambda}$ of representations or the invariant subspace $\Inv(V_{\undlambda})$ in it by a triangulated category $\mathcal{C}_{\underline{\lambda}}$ so that its Grothendieck group $K_0(\mathcal{C}_{\underline{\lambda}})$ is an integral version of $V_{\undlambda}$ or its invariant subspace. 
Namely, $\mathcal{C}_{\underline{\lambda}}$ must come with an automorphism $\{1\}$ that respects the triangulated category structure. This automorphism is often given by a grading shift, when $\mathcal{C}_{\underline{\lambda}}$ is built from the category of suitable graded modules over a graded ring. The automorphism makes the Grothendieck group $K_0(\mathcal{C}_{\underline{\lambda}})$ into a $\Z[q,q^{-1}]$-module, with $\{1\}$ inducing multiplication by $q$ on the Grothendieck group. 

Next, there should be a homomorphism of $\Z[q,q^{-1}]$-modules
\begin{equation}
    K_0(\mathcal{C}_{\underline{\lambda}}) \lra V_{\undlambda} \ \  \mathrm{or} \ \ 
    K_0(\mathcal{C}_{\underline{\lambda}}) \lra 
    \Inv(V_{\undlambda}) 
\end{equation}
that becomes an isomorphism upon tensoring $K_0$ with $\Q(q)$ over $\Z[q,q^{-1}]$: 
\begin{equation}
    K_0(\mathcal{C}_{\underline{\lambda}})\otimes_{\Z[q,q^{-1}]}\Q(q) \stackrel{\cong}{\lra} V_{\undlambda} \ \  \mathrm{or} \ \ 
    K_0(\mathcal{C}_{\underline{\lambda}})\otimes_{\Z[q,q^{-1}]}\Q(q)  \stackrel{\cong}{\lra}  
    \Inv(V_{\undlambda}). 
\end{equation}
Equivalently, this can be phrased as an isomorphism between the Grothendieck group and an integral version (over $\Z[q,q^{-1}]$) of the tensor product module 
\begin{equation}
    K_0(\mathcal{C}_{\undlambda})\cong V^\mathbb{Z}_{\undlambda} \quad \mathrm{or} \quad 
    K_0(\mathcal{C}_{\undlambda})\cong \Inv(V^\mathbb{Z}_{\undlambda}),  \quad 
    V^\mathbb{Z}_{\undlambda} \otimes_{\Z[q,q^{-1}]}\Q(q) \cong V_{\undlambda}. 
\end{equation}
In particular, $V^\mathbb{Z}_{\undlambda}$ are free $\Z[q,q^{-1}]$-modules of rank equal to the dimension of $V_{\undlambda}$ or $\Inv(V_{\undlambda})$, respectively. For all tangles $T$ the map $f(T)$ should restrict to a linear map between the integral versions of tensor products:  
\[
f(T):V^{\mbZ}_{\undlambda}\lra V^{\mbZ}_{\underline{\mu}}. 
\] 
One collection of such integral lattices $V_{\undlambda}^{\mbZ}$ is given by $\Z[q,q^{-1}]$-spans of the duals of Lusztig canonical bases in $V_{\undlambda}$. 
We refer to~\cites{Lu,W2} and~\cite{TLe}
for these and related integrality results. 

\vspace{0.1in} 

Upon categorification, 
intertwiner $f(T)$ for a tangle $T$ should be 
replaced by an exact functor 
\begin{equation}
\mathcal{F}(T): \mathcal{C}_{\underline{\lambda}}\to \mathcal{C}_{\underline{\mu}}
\end{equation}
that on the Grothendieck group level descends to $f(T)$, so that the following diagram commutes: 
\[
\begin{CD}
K_0(\mcC_{\undlambda}) @>[\mathcal{F}(T)]>> K_0(\mcC_{\underline{\mu}})\\
@VV{\cong}V @VV{\cong}V\\
V^{\mbZ}_{\undlambda} @>{f(T)}>> V^{\mbZ}_{\underline{\mu}}
\end{CD}
\]
Composition of tangles should match the composition of functors, with isomorphisms $\mcF(T_2\circ T_1)\cong \mcF(T_2)\circ \mcF(T_1)$ for composable tangles. 

A link $L$ is a tangle between empty sequences of weights. The category $\mcC_{\emptyset}$ is expected to be some base category, such as the category of complexes of graded vector spaces over a field $\Bbbk$, modulo chain homotopies and with finite-dimensional total cohomology groups. Functor $\mcF(L)$ is then given by tensoring with a bigraded vector space, $\mcF(L)(X) \cong \mathsf{H}(L)\otimes X$, and bigraded groups $\mathsf{H}(L)$ define a homology theory for oriented links $L$.  

A lifting $f\mapsto \mcF$ should take us one dimension up. To tangles $T$ we associate functors $\mcF(T)$. Assume that $S$ is a tangle cobordism from a tangle $T_0=\partial_0 S$ to $T_1=\partial_1 S$. A tangle cobordism is an oriented surface in $\R^2\times [0,1]^2$ with boundary and corners, viewed as a cobordism from $\partial_0 S$ to $\partial_1 S$. They were discussed earlier, at the end of Section~\ref{subsec_sl3web}, with the 2-category of tangle cobordisms denoted $\Tcob$ there. 

Since our links have components decorated by irreducible representations of $\mfg$, components of tangle cobordisms should be decorated by such representations as well, and we denote the corresponding 2-category by $\Tcob_{\mfg}$. 

The invariant $\mcF(T)$ should extend to tangle cobordisms and assign a natural transformation 
\begin{equation}
\mcF(S): \mcF(T_0) \lra \mcF(T_1)
\end{equation}
of functors to a tangle cobordism $S$ from $T_0$ to $T_1$. 
These transformations should fit together into a 2-functor 
\begin{equation}\label{eq_Tcob_2f}
\mcF \ : \ \Tcob_{\mfg} \lra \mathsf{NT}
\end{equation}
from the 2-category of tangle cobordisms to a suitable 2-category $\mathsf{NT}$ (natural transformations) with 
\begin{itemize}
    \item Objects: sequences $\undlambda$, 
    \item 1-morphisms: exact functors from $\mcC_{\undlambda}$ to $\mcC_{\underline{\mu}}$, 
    \item 2-morphisms: natural transformations between these exact functors. 
\end{itemize}

Specializing from tangle cobordisms to link cobordisms $S$, the invariant of a latter is a homogeneous homomorphism of bigraded groups (or vector spaces) 
\[\mathsf{H}(S) \ : \ \mathsf{H}(\partial_0 S)\lra \mathsf{H}(\partial_1 S).
\]
Restricting to the category of link cobordisms, the assignment $L \mapsto \mathsf{H}(L)$ and $S\mapsto \mathsf{H}(S)$ is a functor from that category to the category of bigraded abelian groups (or vector spaces) and homogeneous homomorphisms. 

A categorified Reshetikhin-Turaev invariant assigns a functor to a tangle, and  a natural transformation between functors to a tangle cobordism. On the level of quantum invariants we deal with $3$-dimensional constructions, while their categorifications produce invariants of tangle cobordisms, taking us into $4$ dimensions, one dimension up from the quantum invariants, as envisioned by L. Crane and I. B. Frenkel in their foundational work~\cite{CrFr}. 

\vspace{0.1in}

The 2-functor setup in \eqref{eq_Tcob_2f} is the ideal scenario. To date, it has been realized in several related ways for $\mfg=\mathfrak{gl}_N$ and labeling of tangle and tangle cobordism components by fundamental (miniscule) representations $\Lambda^a_q (V)$, for $1\le a\le N-1$. 

Below, we briefly sketch a combinatorial approach to this construction via foam evaluation in~\cite{RW1}, followed by its extension in~\cite{ETW} to tangle cobordisms. This approach works for the most general coefficient ring of symmetric polynomials in $N$ variables, it postpones working with triangulated categories until later steps, and it gives a clean proof of functoriality for tangle cobordisms. 

Other realizations include: 
\begin{itemize}
    \item A matrix factorizations approach in~\cites{KhR,Yn,Wu}, 
    \item An approach via singular-parabolic blocks of highest weight categories for $\mathfrak{sl}_k$, over all $k$, see~\cites{S,MS}. Here $k$ and $N$ are different variables, 
    \item An approach via Webster algebras~\cite{W2},
    \item An approach via coherent sheaves or via Fukaya-Floer categories on suitable quiver varieties~\cites{AS,AN,CK2,M,SS}.
\end{itemize}

Webster's categorification~\cite{W2} of the Reshetikhin-Turaev invariants~\cite{RT} works in full generality, for any simple Lie algebra $\mfg$ and any labeling of components of $T$ or $L$ by irreducible representations of $\mfg$. Webster homology is expected to be functorial under cobordisms when the components are labelled by miniscule representations of $\mfg$, see~\cite[Section 8.4]{W2}. 
Other theories with a good chance for functoriality include Bodish-Elias-Rose spin link homology~\cite{BER} and some finite-dimensional categorifications of the colored Jones polynomial~\cites{Kh5,Me}. 

For non-miniscule representations, Webster~\cite{W2}, Cautis~\cite{Cau}, Cooper-Krushkal~\cite{CoKr}, Stroppel-Sussan~\cite{StSu} homology and related homology theories are unlikely to be functorial. 
This is due to the homology of the unknot being infinite-dimensional over the ground field, as checked in~\cites{CoKr,W2,StSu} for the unknot labeled by 3-dimensional irreducible representation of $\slt$ for several of these theories. It is an important open problem to find functorial modifications of these and other homology theories.  
Several of the these theories require representing a link as the closure of a braid, making their functoriality under cobordisms almost impossible to establish with that definition. 
Categorical approaches to categorified quantum invariants are sophisticated, and they were often done for $\mathfrak{g}=\GL(N)$ or $\SL(N)$ and a specific class of irreducible representations of $\mfg$, such as fundamental representations.
%If we compute functor $\mathcal{F}(T)$ on Grothendieck groups, we get
%$$
%\xymatrix{
%V^\mathbb{Z}_{\underline{\lambda}} \ar[r]^{[\mathcal{F}(T)]} & V^\mathbb{Z}_{\underline{\mu}}
%}
%$$
%$$
%V^\mathbb{Z}_{\underline{\lambda}}\overset{[\mathcal{F}(T)]}%{\longrightarrow} V^\mathbb{Z}_{\underline{\mu}}
%$$

\vspace{0.07in} 

In order to proceed with an alternative, more combinatorial approach, beyond a categorification of the Jones polynomial and the Kuperberg bracket (which are the quantum $\SL(2)$ and $\SL(3)$ invariants for the fundamental representation, respectively)  one needs to come up with a sophisticated foam evaluation formula, which is a bit different from what has been described in previous lectures. This step had to wait several years and was solved in a spectacular fashion by Robert and Wagner~\cite{RW1} for $\GL(N)$ foams in the oriented case (and see~\cites{MSV,QR,RWe,MaVa2,We} for different earlier approaches to $\GL(N)$ and to triply-graded link homology from foams). At the time only oriented $\SL(3)$ closed foam evaluation was available, where vertices were absent. Robert and Wagner's work, together with the Kronheimer-Mrowka orbifold homology~\cite{KM15}, also led to the evaluation formula for unoriented $\SL(3)$ foams~\cite{KhR}, described in Section~\ref{sec-one}. 

\vspace{0.07in}

Here and below we mostly gloss over the differences between $\SL(N)$ and $\GL(N)$ webs and $\SL(N)$ and $\GL(N)$ foams and refer to~\cite{We2}  for more details. On the level of webs, edges of thickness $N$ are present in $\GL(N)$ webs but not in $\SL(N)$ webs, since the top exterior power $\Lambda^N(V)$ is the trivial representation of $\SL(N)$ but not of $\GL(N)$. In an $\SL(N)$ web, three lines of thickness $a,b,c$ with $a+b+c=N$ can merge into a vertex, see Figure~\ref{5_08} on the left, but  not in a $\GL(N)$ web. 
Differences in $\SL(N)$ and $\GL(N)$ foams, and especially in their evaluations are more subtle, and we specialize to $\GL(N)$ foams here. Robert-Wagner evaluation~\cite{RW1} is that for $\GL(N)$ foams, in our terminology.  

\vspace{0.07in}

\begin{figure}[!htbp]
	\begin{center}
		{\includegraphics[width=230pt]{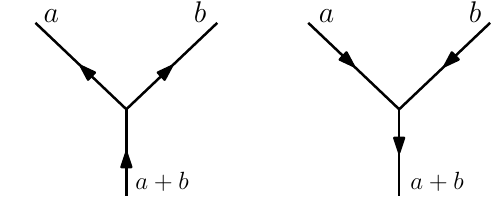}}
		\caption{Vertices of MOY graphs, where thick lines split and merge.}
		\label{5_01}
	\end{center}
\end{figure}

Starting at the decategorified level of webs, one considers Murakami-Ohtsuki-Yamada graphs (called \emph{MOY graphs or webs}), representing compositions of specific intertwiners between tensor products of fundamental representations of quantum $\mathfrak{gl}_n$~\cite{MOY}. In this case we allow vertices near which the lines have thickness $1\leq a,b,a+b\leq N$, respectively, and these lines split and merge as shown in Figure~\ref{5_01}. A boundary point labelled $a$ represents the irreducible representation $\Lambda^a_q(V)$ of quantum $\mathfrak{gl}_n$ which quantizes the $a$-th exterior power $\Lambda^a(V)$ of the fundamental representation of the Lie algebra $\mathfrak{gl}_n$. 
Trivalent vertices in Figure~\ref{5_01}  represent inclusion and projection intertwiners between representations of quantum groups:
\[
\Lambda^{a+b}_q(V) \lra 
\Lambda^{a}_q(V) \otimes  \Lambda^{b}_q(V), \ \ \ 
\Lambda^{a}_q(V) \otimes  \Lambda^{b}_q(V) \lra 
\Lambda^{a+b}_q(V). 
\]

\begin{remark}
To consider $\SL(N)$ webs and foams, where for the fundamental representation we have $\wedge^N_q (V)\simeq\mathbb{C}$ (the top $q$-exterior power of the fundamental representation is trivial), one can add vertices of new types, where weights and orientations are shown in Figure \ref{5_08}.
\begin{figure}[H]
	\begin{center}
		{\includegraphics[width=330pt]{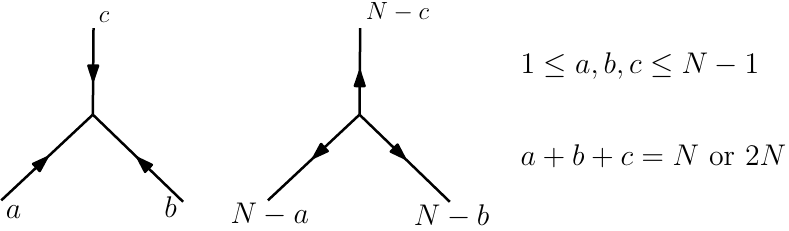}}
		\caption{Two types of vertices in $\SL(N)$ MOY webs, with all orientations either in or out.}
		\label{5_08}
	\end{center}
\end{figure}
The case $N=3$, restricted to webs with all edges of weight 1, reproduces webs for the Kuperberg invariant~\cite{Kup} and foam evaluation for its categorification in Section~\ref{sec_oriented}.
\end{remark}

\begin{figure}[H]
	\begin{center}
		{\includegraphics[width=400pt]{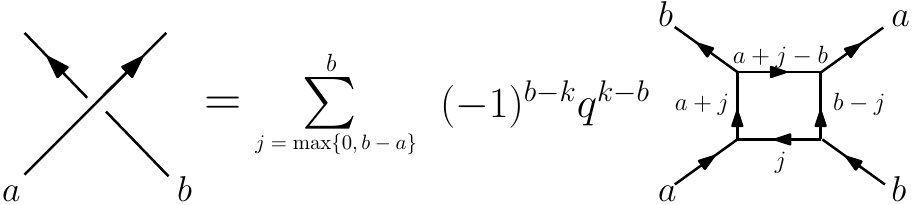}}
		\caption{Decomposition of an $(a,b)$-crossing into a $q$-linear combination of planar MOY webs.}
		\label{5_02}
	\end{center}
\end{figure}

\begin{figure}[H]
	\begin{center}
		{\includegraphics[width=200pt]{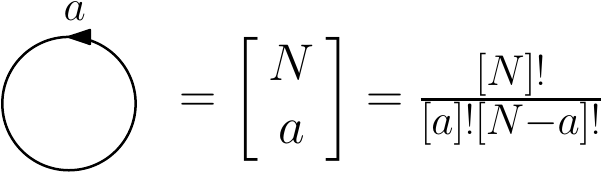}}
		\caption{A circle of thickness $a$ evaluates to the quantum binomial coefficient, also see \eqref{eq_quantum_n} for the definition of quantum integer $[n]$.}
		\label{5_03}
	\end{center}
\end{figure}

A crossing decomposes as a linear combination of compositions of these intertwiners, as shown in Figure \ref{5_02}. A similar formula for the opposite crossing is omitted.
Normalization of the unknot is chosen as in Figure \ref{5_03}.
Using the Murakami-Ohtsuki-Yamada graphical calculus~\cite{MOY}, a closed MOY web $\Gamma$ evaluates to a Laurent polynomial in $q$ with non-negative integer coefficients $P_N(\Gamma)\in \Z_+[q,q^{-1}]$. An oriented link $L$ with components labeled by exterior powers of the fundamental representation, i.e., by numbers from $1$ to $N$, evaluates to a Laurent polynomial with integer coefficients $P_N(L)\in \Z[q,q^{-1}]$. 

A full set of skein relations on MOY webs with boundary was derived by Cautis, Kamnitzer and Morrison~\cite{CKM}.

\vspace{0.07in}

\begin{figure}[!htbp]
	\begin{center}
		{\includegraphics[width=350pt]{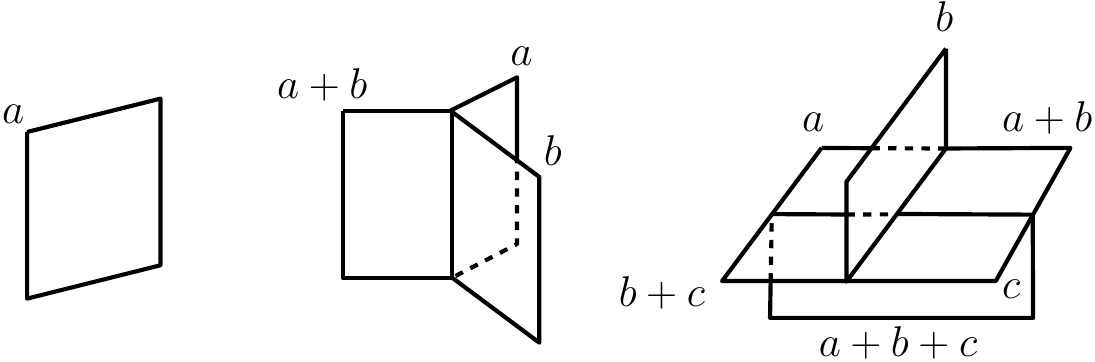}}
		\caption{Three types of points of a $\GL(N)$ foam, left to right: regular points, points along a seam, and vertices.}
		\label{5_04}
	\end{center}
\end{figure}

\begin{figure}[!htbp]
	\begin{center}
		{\includegraphics[width=350pt]{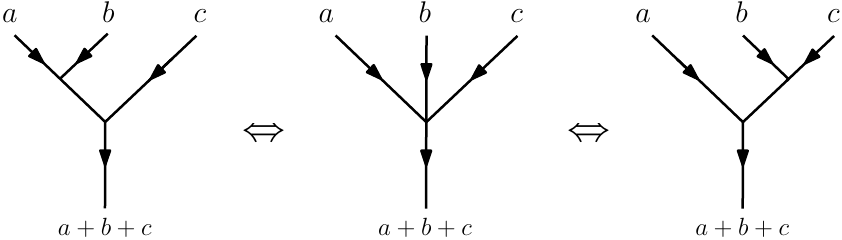}}
		\caption{Associativity diagram for merging three strands of thickness $a$, $b$, $c$ into a single strand of thickness $a+b+c$.}
		\label{5_007}
	\end{center}
\end{figure}

\begin{figure}[!htbp]
	\begin{center}
		{\includegraphics[width=120pt]{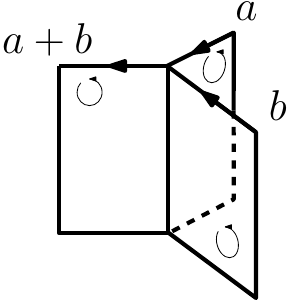}}
		\caption{Orientations of three facets along a seam of a $\GL(N)$ foam.}
		\label{5_05}
	\end{center}
\end{figure}

To categorify this $\GL(N)$ invariant of webs, links and tangles one needs to go one dimension up, work with $\GL(N)$ foams and find a suitable evaluation for them. Generalizing from $\SL(3)$ to $\GL(N)$, we allow facets of a foam to be of thickness $1\leq a\leq N$. Possible types of points on a foam (or singularities of a foam) are depicted in Figure~\ref{5_04}.
Along a seam, facets of thickness $a$ and $b$ merge into a facet of thickness $a+b$, see Figure~\ref{5_05} and the middle picture in Figure~\ref{5_04}. 
A neighborhood of a  singular vertex is shown on the right side of Figure \ref{5_04}. At the vertex facets of thickness $a$, $b$, $c$ merge into facets of thickness $a+b$, $b+c$, and eventually, into the facet of thickness $a+b+c$. Taking parallel cross-sections of this foam near the vertex produces the associativity diagram for the web merge of three strands, as in Figure~\ref{5_007}.
Every facet is oriented and facet orientation is preserved when moving from a facet of thickness $a$ or $b$ to the adjacent facet of thickness $a+b$.
Along a seam, moving between facets of thickness $a$ and $b$ reverses orientation, see Figure~\ref{5_05}. 

Similar to $\SL(3)$ foams, dots on facets of a $\GL(N)$ foam are allowed.  
A dot $d$ on a facet $f$ of thickness $a$ carries a homogeneous  symmetric functions in $a$ variables
\begin{equation}\label{eq_P_d_a}
P_d\in \mathbb{Z}[y_1,\ldots,y_{a}]^{S_{a}}.
\end{equation}
Homogeneity condition is imposed so that the degree of a foam is well-defined. 
Let $\Dots(F)$ be the set of dots on a foam $F$. By $f(d)$ denote the facet that contains the dot and by $a(d)$ the thickness of the facet, so that $a=a(d)$ in \eqref{eq_P_d_a}.

A coloring of $F$ is a map
$$
c:f(F)\to 2^{I_N},
$$
where $I_N=\{1,\ldots, N\}$ is the set of $N$ colors. A coloring $c$ assigns a subset $c(f)\subset I_N$ of cardinality $a$, $|c(f)|=a$, to a facet $f$ of $F$ of thickness $a$. These subsets should satisfy the following compatibility condition along seams. 

Denote thin facets of thickness $a,b$ along a seam by $f_1,f_2$ and the thick facet by $f_3$. We require that subsets $c(f_1),c(f_2)$ are disjoint and their union is $c(f_3)$: 
\begin{equation}\label{eq_disjoint1}
c(f_3) =c(f_1) \sqcup c(f_2), \ \ |c(f_1)|=a, \ |c(f_2)|=b, \ c(f_1)\cap c(f_2)= \emptyset. 
\end{equation}
The following \emph{flow rule} restates the above relation: each color of a thin facet, say $i\in c(f_1)$, flows into the same color of the thick facet along the seam, $i\in c(f_3)$, and that is how all colors of the thick facet are obtained. Adjacent thin facets along a seam cannot contain the same color.   
\begin{figure}[!htbp]
	\begin{center}
		{\includegraphics[width=300pt]{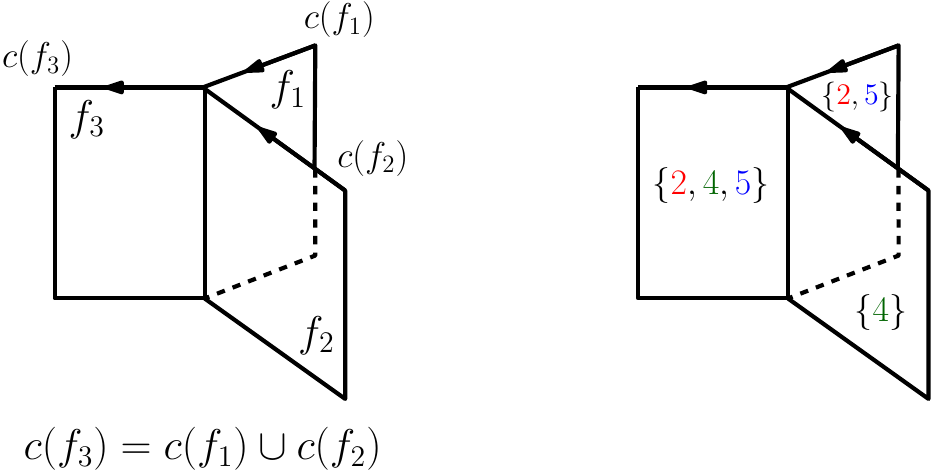}}
		\caption{Left: the flow rule along the seam. Right: an example of a flow along the seam, where $a=2$, $c(f_1)=\{2,5\}$ and $b=1$, $c(f_2)=\{4\}$.}
		\label{5_005}
	\end{center}
\end{figure} 

This rule also tells us about possible colorings of facets near a vertex, where four seams merge. At a vertex shown in Figure~\ref{5_04} thin facets (of thickness $a,b,c$) are colored by pairwise-disjoint subsets $c(f_1),c(f_2),c(f_3)$ of $I_N$, and the remaining three facets are colored by suitable unions of these subsets.

Taking the union of closed facets that contain a given color $i$ yields a closed surface in a foam, denoted $F_i(c)\subset F$ and called a \emph{unicolored surface} of $F$ and $c$. It inherits orientation from that of facets of $F$. 
The flow rule above implies that $F_j(c)$ is an oriented closed surface. It is embedded in $\R^3$ since $F\subset \R^3$. 
Being a closed surface in $F\subset \R^3$, it has even Euler characteristic $\chi(F_i(c))$. 

Define a \emph{bicolored surface} $F_{ij}(c)$ for $1\le i<j\le N$, as the closure of the symmetric difference 
\[
F_{ij}(c) \ := \ \overline{F_i(c)\Delta F_j(c)}. 
\]
It consists of the union of closed facets $f$ such that $c(f)$ contains exactly one element of the set $\{i,j\}$. That is, the set $c(f)$ contains either $i$ but not $j$ or $j$ but not $i$. Taking the closure means that whenever an open facet $f$ belongs to $F_{ij}(c)$, its boundary, a union of seams and vertices, also does. The following proposition is easy to verify by checking the surface property near a vertex of $F$. 
\begin{prop}
    $F_{ij}(c)$ is a closed surface embedded in $\R^3$. 
\end{prop}
Consequently, the Euler characteristic $\chi(F_{ij}(c))$ is an even integer. 

\vspace{0.07in}

Given a dot $d$ on a facet $f$ of $F$ carrying a symmetric polynomial $P_d$ and a coloring $c$ of $F$, define the 
polynomial $P_d(c)\in \Z[x_1,\dots, x_N]$  by relabeling variables $y_1,\dots, y_a$, see \eqref{eq_P_d_a}, into variables $x_k$, where $k$ runs over all $a$ colors in $c(f)$, for the facet $f=f(d)$. Since $P_d$ is symmetric, the choice of a bijection from $y_1,\dots, y_a$ to $\{x_k\}_{k\in c(f)}$  to set up a homomorphism 
\begin{equation}\label{eq_homo_S}
\Z[y_1,\dots, y_a]^{S_a}\ \lra \ \Z[x_1,\dots, x_N]
\end{equation}
taking $P_d$ to $P_d(c)$ does not matter.  
For example, if $a=2$, $P_f=y_1^2+y_2^2-y_1 y_2$ and $c(f)=\{3,5\}$, then $P_f(c)=x_3^2+x_5^2-x_3 x_5$. 

\vspace{0.07in} 

Define the evaluation of a foam $F$ with a given coloring $c$ by  
\begin{equation}
\label{eq_eval_F_c}
\langle F, c\rangle \ :=\ (-1)^{s(F,c)}\frac{P(F,c)}{Q(F,c)}
\end{equation}
where
\begin{eqnarray}
\label{eq_first_fla}
P(F,c) & = & \prod\limits_{d\in \Dots(F)} P_d(c) \ \in \ \Z[x_1,\dots, x_N], \\
\label{eq_second_fla}
Q(F,c) & = & \prod\limits_{1\leq i<j\leq N}(x_i-x_j)^{\chi(F_{ij}(c))/2},
\end{eqnarray}
and
\begin{eqnarray}
\label{eq_three_fla}
s(F,c) & = & \theta^+(c) +\sum\limits_{j=1}^N j\frac{\chi(F_j(c))}{2}, \\
\label{eq_four_fla}
\theta^+(c) & = & \sum\limits_{i<j}\theta^+_{ij}(c).
\end{eqnarray}
The numerator $P(F,c)$ is the product of labels of dots on $F$, converted into a homogeneous polynomial in $x_1,\dots, x_N$ via homomorphisms \eqref{eq_homo_S}. The denominator is a product of powers of $x_i-x_j$'s that depend on the Euler characteristics of bicolored surfaces $F_{ij}(c)$. The exponents are integers that may be either positive, negative or $0$. 

Expression $s(F,c)$ encodes the sign in \eqref{eq_eval_F_c} and it incorporates, as one summand, the sum of scaled Euler characteristics of unicolored surfaces $F_j(c)$. 

The other summand involves $\theta^{+}_{ij}(c)$, defined as the number of \emph{positive} singular circles on the bicolored surface $F_{ij}(c)$. Namely, for a given $c$, consider the union $F(i,j;c)$ of facets $f$ of $F$ that contain at least one of the colors $i$ and $j$ in $c(f)$. Subspace $F(i,j;c)$ contains the bicolored surface $F_{ij}(c)$ along with the facets $f$ such that $\{i,j\}\subset c(f)$. One can view $F(i,j;c)$ as a $\GL(2)$ foam with $F_{ij}(c)$ the union of 1-facets (colored by a single color $i$ or $j$) and the complement given by the union of 2-facets or \emph{double facets} (colored by $\{i,j\}$). 

Bicolored surface $F_{ij}(c)$ may contain \emph{singular circles}, along which the double facets of $F(i,j;c)$ are attached. Thus, along a singular circle, there are two 1-facets and a double facet attached, colored by $i$,$j$, and $\{i,j\}$, respectively. The circle inherits orientation from that of adjacent single and double facets. 
 
\begin{figure}[!htbp]
	\begin{center}
		{\includegraphics[width=300pt]{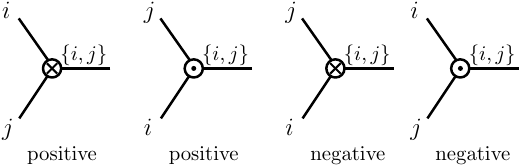}}
		\caption{Positive and negative types of seam circles on the 2-colored subfoam $F_i(c)\cup F_j(c)\subset F$, where $i<j$. Central dot denotes orientation of the singular circle pointing in the direction outside of the plane of the figure, while a cross denotes the opposite orientation.}
		\label{5_06}
	\end{center}
\end{figure}

Define a singular circle $C$ to be \emph{positive} if  by looking in the direction of orientation of $C$ and rotating counterclockwise in the orthogonal plane, the coloring of attached facets is in the cyclic order $(i,j,\{i,j\})$ for $i<j$ as shown on the leftmost diagram in Figure~\ref{5_06}. The second diagram in Figure~\ref{5_06} shows the same cross-section viewed from the opposite direction, against the orientation of a circle. The number of positive singular circles on $F_{ij}(c)$ is denoted $\theta_{ij}^+(c)$, and it enters formula \eqref{eq_four_fla} in the definition of $\langle F,c\rangle$. 
The evaluation of a foam $F$ is given by the sum of evaluations over all colorings $c$ of $F$: 
\begin{equation}
\label{eq_eval_F}
\langle F\rangle \ :=\ \sum_{c}\langle F, c\rangle.
\end{equation} 
\begin{prop}~\cite{RW1}\label{prop_N_integral}
    For any closed foam $F$ its evaluation $\langle F\rangle$ belongs to the subring of symmetric polynomials 
    \begin{equation}\label{eq_ground_ring}
        R_N \ := \ \Z[x_1,\dots,x_N]^{S_N}.
    \end{equation} 
\end{prop}
Theorem~\ref{thm_in_R} is an analogous integrality result, for unoriented $\SL(3)$ foams. The above proposition and its proof came out first, and Theorem~\ref{thm_in_R} followed in~\cite{KhR}.  The proof consists of looking at 2-spheres in the bicolored surface $F_{ij}(c)$ and showing that Kempe moves of colorings allow to cancel out potential denominators $x_i-x_j$ in the evaluation. Part of the proof requires carefully checking changes in signs $s(F,c)$ in \eqref{eq_three_fla} under a Kempe move, which is not needed for  Theorem~\ref{thm_in_R}, where the ground ring has characteristic 2.  

\vspace{0.07in} 

Given Proposition~\ref{prop_N_integral} integrality result, Robert and Wagner~\cite{RW1} define state spaces $\langle\Gamma\rangle$ using the universal construction, for a $\GL(N)$ web (or MOY graph) $\Gamma$. These state spaces are graded $R$-modules. They then derive inductive direct sum decompositions of state spaces $\langle\Gamma\rangle$ analogous to those derived in Section~\ref{section2} for unoriented $\SL(3)$ webs but technically much more involved~\cite{RW1}. 
Their detailed study of $\GL(N)$ foams and associated state spaces leads them to the following result~\cite{RW1}. 
\begin{prop} \label{prop_rank_MOY}
The state space of any $\GL(N)$ MOY web $\Gamma$ is a free graded $R$-module of graded rank equal to the quantum MOY invariant: 
\begin{equation}\label{eq_rank_MOY}
\mathsf{grk}_R (\langle \Gamma \rangle) = P_N(\Gamma).
\end{equation} 
\end{prop}
Here $P_N(\Gamma)$ is the Murakami-Ohtsuki-Yamada~\cite{MOY} invariant of planar $\GL(N)$ webs $\Gamma$ discussed earlier in this section. 

\begin{remark}
For the simplest MOY graphs, their state spaces can be interpreted via equivariant cohomology groups of partial or iterated flag varieties. Earlier in the lectures we have discussed the relation between unoriented and oriented $\SL(3)$ webs and foams and cohomology of projective spaces and flag varieties for $\C^3$, see the discussion after Theorem~\ref{taitth} and in Section~\ref{subsec_Kuperberg}. 
For the $\GL(N)$ webs, from the formula in Figure~\ref{5_03} we can guess that the state space associated to the planar circle of thickness $a$ should be isomorphic to the $U(N)$-equivariant cohomology of the complex Grassmannian, $\textrm{Gr}(a,N)$ of $a$-dimensional planes in $\C^N$. This indeed turns out to be the case~\cite{RW1}. The state space of the $\GL(N)$ $\Theta$-web, a planar web with 2 vertices and three edges of thickness $a,b$ and $a+b$ between them is isomorphic to the $U(N)$-equivariant cohomology of the variety of partial flags $\{L\subset L'\subset \C^N\}$ with $\dim(L)=a,\dim(L'/L)=b$. 
\end{remark}

Graded $R$-module invariant $\langle \Gamma\rangle$ of MOY webs is functorial: a MOY foam $F$ with boundary induces a homogeneous homomorphism of graded $R$-modules
\[
\brak{F} \ : \ \brak{\partial_0 F} \lra \brak{\partial_1 F} .  
\]
These homomorphisms come from the universal construction setup, and their counterparts for oriented and unoriented $\SL(3)$ foams were discussed earlier in the lectures. 

\begin{prop}~\cite{RW1}
    State spaces $\brak{\Gamma}$ for MOY webs $\Gamma$ and homomorphisms $\brak{F}$ for MOY foams $F$ with boundary constitute a TQFT $\brak{\ast}$ on the category of MOY $\GL(N)$ foams taking values in the category of (free) graded $R$-modules 
    \[
    \brak{\ast} \ : \ \mathsf{Foams}_{\GL(N)} \lra R\mathsf{-gmod}.
    \]
\end{prop}

With this proposition in hand, one can now mimic the construction of link homology groups in Section~\ref{subsec_cat_Ku}. Recall the expression for the $(a,b)$-crossing in Figure~\ref{5_02}  as an alternating sum of web diagrams (squares with spikes) with powers of $q$ with coefficients. Compare this decomposition to essentially its special case in Figure~\ref{3_04}, when $a=b=1$, $N=3$ and one replaces $\GL(3)$ with $\SL(3)$. The latter two-term difference is categorified by the complexes in Figure~\ref{cones_fig}, with the differentials induced by foams in Figure~\ref{3_13}.  

A generalization of this complex to categorify an arbitrary $(a,b)$-crossing in Figure~\ref{5_02} on the left uses the webs in Figure~\ref{5_02} on the right and differentials between their state spaces for $j$ and $j+1$ induced by appropriate foams~\cite{ETW}, also see earlier work~\cites{QR,MSV} and papers~\cites{Wu,Yn} where these complexes appear in the context of matrix factorizations. Given a planar projection $D$ of a link $L$ with $n$ crossings, one can set up an $n$-dimensional hyperrectangle of suitable size $k_1\times \dots \times k_n$ with state spaces of suitable webs in integer vertices and maps induced by foams connecting state spaces in adjacent vertices. Passing to the total complex of this commutative hyperrectangle and taking cohomology groups results in a bigraded $R_N$-module $\mathsf{H}_N(L)$, for an oriented link $L$ with components colored by numbers from $1$ to $N-1$ (fundamental weights). This results in a bigraded homology theory that categorifies the MOY quantum invariant. Several techniques, including the canopolis framework~\cite{BN} and related approaches discussed in these lectures, extend this construction to tangles and tangle cobordisms, giving a 2-functor as in \eqref{eq_Tcob_2f}, with the above restriction on labels of components.   

An extension of functoriality from link cobordisms in $\R^3\times [0,1]$ to those in $\S^3\times [0,1]$ together with its applications to lasagna skein modules of 4-manifolds was done in~\cite{MWW}. 

\vspace{0.07in}

The construction sketched above is specific to the exterior powers of the fundamental representation of $\GL(N)$ and MOY intertwiner networks. As an intermediate step one uses that quantum invariants of planar MOY webs have positivity and integrality property $P_N(\Gamma)\in\mathbb{Z}_{+}[q,q^{-1}]$. This property is essential for building homology groups in a manageable way using the universal construction, forming state spaces for webs and defining chain complexes for link projections via these state spaces and maps between them induced by appropriate foams.  

\vspace{0.07in}

Taking the smallest simple Lie algebra $\slt$, we observe that 
intertwiners between arbitrary representations of $U_q(\slt)$  do not possess positivity and integrality properties. In particular, integrality and positivity properties are absent from the theory of planar $q$-spin networks~\cites{KR, KL, CFS}. Lack of integrality and positivity for closed intertwiner networks  is a major obstacle to a possible functorial categorification of quantum $\slt$ invariants of links ({\it i.e.}, of the colored Jones polynomial) and to a  categorification  of $q$-spin networks embedded in $\R^3$. It is not known whether various finite-dimensional categorifications of the colored Jones polynomial, see~\cite{Me} and references therein, are functorial under link cobordisms. 

In a limiting case of this construction, when $N\to\infty$, foam evaluation for braid-like graphs and foams provides a combinatorial description of the category of Soergel and singular Soergel bimodules. For braid-like MOY graphs $\Gamma_0,\Gamma_1$ the space of all braid-like foams from $\Gamma_0$ to $\Gamma_1$  modulo relations is isomorphic to the space of homomorphisms between singular Soergel bimodules associated to $\Gamma_0,\Gamma_1$, see~\cites{Vaz,MaVa2},~\cites{RWe,We,GoWe} for the original relations between the (singular) Soergel category and foams and \cites{RW2,KhRW} for an  approach via foam evaluation. Individual singular Soergel bimodules can also be interpreted via foams with boundary and corners. We refer to the above papers for the notions of braid-like webs and foams. Soergel bimodules and their homomorphisms are recovered by specializing to braid-like MOY graphs with the lines of thickness $1$ and $2$ only. 

A surprising action of the Lie algebra $\mathfrak{sl}_2$ on $\GL(N)$ web state spaces and associated link homology groups has recently been constructed in~\cites{QRSW1,QRSW2}. Defining the action requires some amount of equivariance, so that the ground ring must have a least one polynomial generator in a positive degree and cannot be a field. It is an open problem to understand these actions and connect them to well-known $\mathfrak{sl}_2$ actions, such as in the Lefschetz theory or on the Asaeda-Przytycki-Sikora annular homology~\cite{GLW}. 

\section{Kronheimer-Mrowka  theory and the Four-Color Theorem}\label{KMtheory}
Kronheimer-Mrowka $\SO(3)$ gauge theory for $3$-orbifolds was developed in \cites{KM15, KM16, KM17} . The orbifolds can have lines with a $\mathbb{Z}_2$-structure, locally represented by a line in $\R^3$ modulo the rotation by $\pi$, and $3$-vertices, which locally correspond to the quotient $\mathbb{R}^3/V_4$, where $V_4\subset \SO(3)$ is the Klein group. This theory gives rise to a homology theory $J^\sharp$ for 3-orbifolds, which is a version of the instanton Floer homology. Kronheimer-Mrowka theory has remarkable connections to the Four-Color Theorem. 

Any trivalent graph $\Gamma$ embedded in a $3$-manifold gives a $3$-orbifold. Restricting to trivalent graphs embedded to $\mathbb{R}^3$, i.e. knots and links with trivalent vertices, one obtains a functorial homology theory for such objects. Further restricting to $\mathbb{R}^2$ and using embedding $\mathbb{R}^2\subset \mathbb{R}^3$, one obtains a homology theory for planar webs. This theory assigns a finite-dimensional $\mathbb{F}_2$-vector space to a web, without a grading.

We say that graph $\Gamma\subset \R^2$ has a {\it bridge} when there exists an embedded circle in $\R^2$, which intersects $\Gamma$ at a single point of an edge. Kronheimer and Mrowka prove the  following non-vanishing theorem:
\begin{theorem}\label{nonvanish}
For a planar trivalent graph $\Gamma$
\begin{equation*}
    J^\sharp(\Gamma)=0\quad\text{if and only if}\quad\text{graph $\Gamma$ has a bridge}.
\end{equation*}
\end{theorem}
The proof of this theorem and its generalization in~\cite{KM15} relies on Floer homology and Gabai's work on sutured manifolds. 
Based on computations of $J^\sharp$ for some webs and its properties, Kronheimer and Mrowka propose the following conjecture~\cite{KM15}. 
\begin{conj}\label{conject}
For a planar trivalent graph $\Gamma$
\begin{equation*}
\textrm{dim}_{\mathbb{F}_2}(J^\sharp(\Gamma))=t(\Gamma),
\end{equation*}
where $t(\Gamma)=|\textrm{Tait}(\Gamma)|$.
\end{conj}
Truth of this conjecture, combined with Theorem \ref{nonvanish}, would imply the Four-Color Theorem, giving a conceptual proof of the latter that bypasses approaches requiring computer-aided case-by-case analysis~\cites{AH1,AH2,RSST,St}. In \cite{KM17} it was shown that
\begin{equation}\label{eqKM}
t(\Gamma)\leq \textrm{dim}_{\mathbb{F}_2}(J^\sharp(\Gamma)).
\end{equation}
%but the status of the Conjecture \ref{conject} is currently unknown. 

Kronheimer and Mrowka conjectured the existence of a combinatorial counterpart of their homology theory, restricted to planar trivalent graphs and their foam cobordisms~\cite{KM15}. 
A specialization of the construction presented in Sections \ref{sec-one} and \ref{section2} led to a proof of their conjecture~\cite{KhR}. 
The idea is to downsize the evaluation from taking values in the ring $R=\Bbbk[E_1,E_2,E_3]$ of symmetric 3-variable polynomials to taking values in the ground field $\Bbbk$ of characteristic two.

Starting with closed foam evaluation $\langle F \rangle\in R$, we can
do a base change given by a homomorphism
\begin{eqnarray*}
f_0:R & \to & \Bbbk\\
E_1,E_2,E_3& \mapsto & 0,
\end{eqnarray*}
where every element of a positive degree goes to $0$. Denote composite evaluation by 
\[
\langle F\rangle_0 \ := \ f_0(\langle F\rangle) \ \in \Bbbk
\]
and apply the universal construction to it to get a functor $\langle\ast\rangle_0$, which gives a state space $\langle\Gamma\rangle_0$ for each planar trivalent web $\Gamma$. It is a finite dimensional $\Bbbk$-vector space and it has a grading, where the degree of a dot on a facet is $2$. One can check that this state space is a free graded $\Bbbk$-vector space of dimension at most $t(\Gamma)$. The definition and the match with the construction in~\cite{KM15} make it clear that $\langle\Gamma\rangle_0$ is a subquotient of $J^\sharp(\Gamma)$ when the base field $\Bbbk=\mathbb{F}_2$. Together with (\ref{eqKM}) it gives the following inequality of dimensions
\begin{equation*}
    \textrm{dim}\langle\Gamma\rangle_0\leq t(\Gamma)\leq\textrm{dim}(J^\sharp(\Gamma)),
\end{equation*}
where $\textrm{dim}(J^\sharp(\Gamma))$ is known to be nonzero for bridgeless graphs.
\begin{conj}
\begin{equation*}
    \textrm{dim}_{\mathbb{F}_2}\langle\Gamma\rangle_0=\textrm{dim}_{\mathbb{F}_2}(J^\sharp(\Gamma)).
\end{equation*}
\end{conj}

%Consider the space state $\langle\Gamma\rangle$ obtained from the evaluation of a foam $\langle F\rangle\in R$, where $R=\Bbbk[E_1,E_2,E_3]$.

This conjecture would also imply the Four-Color Theorem. Theorem~\ref{taitth} shows that this conjecture holds for reducible webs, but otherwise there is only a limited evidence for this conjecture, see~\cites{B1,B2}. Furthermore, Boozer~\cite{B2} proved that the associated topological theories for foams, as cobordisms between unoriented $\SL(3)$ webs, are not isomorphic. 

Understanding state spaces $\langle \Gamma\rangle_0$ for non-reducible webs $\Gamma$ is a very difficult problem. We refer the reader to Boozer~\cites{B1,B2} for more information and results and to Thatte~\cite{Th} for a possibly related study of unoriented $\SL(3)$ web algebras. 

For another very recent approach to the Four-Color Theorem via TQFTs, likely related to the one above, see~\cite{BMc}.  

\begin{remark}
$\textrm{gdim}\langle\Gamma\rangle_0$ is a Laurent polynomial in $q$, specializing at $q=1$ to the number of Tait colorings $t(\Gamma)$. A more familiar deformation of $t(\Gamma)$ comes from studying the $3$-dimensional irreducible representation $V_2$ of $U_q(\slt)$ and networks of intertwiners one builds from the trivalent vertex representing a generator of the one-dimensional space of invariants of $V_2^{\otimes 3}$.  This one-variable polynomial invariant of planar trivalent graphs, which also extends to an invariant of trivalent spatial graphs, is known as the Yamada polynomial~\cites{Ya,AgK}. The deformation of $t(\Gamma)$ discussed above is different from the Yamada polynomial. This deformation is mysterious, and at present it is unknown how to compute it for non-reducible graphs, with the exception of evaluations in~\cites{B1,B2}.  
\end{remark}
\bibliography{refs}
\bibliographystyle{plainnat}
\end{document}